\documentclass[10pt,sort & compress]{elsarticle}
\usepackage{geometry}
\usepackage{amssymb}
\usepackage{latexsym}
\usepackage{amsfonts}
\usepackage{amsthm}
\usepackage{amsbsy}
\usepackage{booktabs}
\usepackage{amsmath}
\usepackage{graphics,graphicx}
\usepackage{subfigure}
\usepackage{multirow,multicol}
\usepackage{bm}
\usepackage{bbm}
\usepackage{color}
\usepackage{float}
\usepackage{calc}
\usepackage{caption}
\usepackage{dsfont}
\usepackage{ifthen}
\usepackage{mathrsfs}
\usepackage[colorlinks,linkcolor=magenta,citecolor=cyan]{hyperref}
\usepackage{epstopdf}
\usepackage{epsfig}
\usepackage{algorithm}
\usepackage{algorithmic}
\usepackage{pifont}
\usepackage{natbib}
\usepackage{overpic}
\usepackage{psfrag}
\usepackage{rotating}
\usepackage{stmaryrd}
\usepackage{verbatim}
\usepackage{setspace}
\usepackage{tikz}

\newdefinition{remark}{Remark}
\newdefinition{method}{Method}
\newdefinition{example}[theorem]{Example}

\numberwithin{theorem}{section}

\newcommand{\orcid}[1]{\href{https://orcid.org/#1}{\includegraphics[width=8pt]{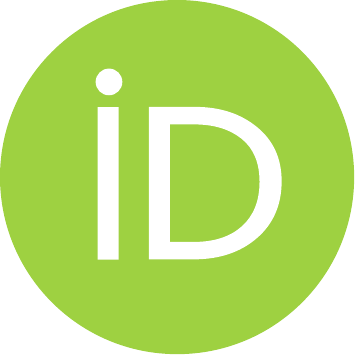}}}
\journal{}
\begin{document}
\begin{frontmatter}
\title{Physical informed neural networks with soft and hard boundary constraints for solving advection-diffusion equations using Fourier expansions}
\author[label2]{Xi'an Li\orcid{0000-0002-1509-9328}\corref{lxa}}
\ead{lixian9131@163.com}
\author[label1]{Jiaxin Deng}
\author[label3]{Jinran Wu\orcid{0000-0002-2388-3614}}
\author[label4]{Shaotong Zhang\orcid{0000-0002-2806-2989}}
\author[label1]{Weide Li\orcid{0000-0003-2288-4751}}
\author[label5]{You-Gan Wang\orcid{0000-0003-0901-4671}}
\ead{ygwanguq2012@gmail.com}

\address[label2]{Ceyear Technologies Co., Ltd., Qingdao 266555, PR China}
\address[label1]{School of Mathematics and Statistics, Lanzhou University, Lanzhou 730000, China, PR China}
\address[label3]{The Institute for Learning Sciences and Teacher Education, Australian Catholic University, Brisbane 4000, Australia}
\address[label4]{Frontiers Science Center for Deep Ocean Multispheres and Earth System, Key Lab of Submarine Geosciences and Prospecting Techniques, MOE and College of Marine Geosciences, Ocean University of China, Qingdao 266100, PR China}
\address[label5]{School of Mathematics and Physics, The University of Queensland, Brisbane, Australia}
\cortext[lxa]{Corresponding author}

\begin{abstract}
Deep learning methods have gained considerable interest in the numerical solution of various partial differential equations (PDEs). One particular focus is physics-informed neural networks (PINN), which integrate physical principles into neural networks. This transforms the process of solving PDEs into optimization problems for neural networks. To address a collection of advection-diffusion equations (ADE) in a range of difficult circumstances, this paper proposes a novel network structure. This architecture integrates the solver, a multi-scale deep neural networks (MscaleDNN) utilized in the PINN method, with a hard constraint technique known as HCPINN. This method introduces a revised formulation of the desired solution for \textcolor{black}{ADE} by utilizing a loss function that incorporates the residuals of the governing equation and penalizes any deviations from the specified boundary and initial constraints. By surpassing the boundary constraints automatically, this method improves the accuracy and efficiency of the PINN technique. To address the ``spectral bias'' phenomenon in neural networks, a subnetwork structure of MscaleDNN and a Fourier-induced activation function are incorporated into the HCPINN, resulting in a hybrid approach called SFHCPINN. The effectiveness of SFHCPINN is demonstrated through various numerical experiments involving \textcolor{black}{ADE} in different dimensions. The numerical results indicate that SFHCPINN outperforms both standard PINN and its subnetwork version with Fourier feature embedding. It achieves remarkable accuracy and efficiency while effectively handling complex boundary conditions and high-frequency scenarios in ADE.

\end{abstract}
\end{frontmatter}

\section{Introduction}\label{sec:01}

\textcolor{black}{The ADE} is a fundamental equation with widespread applications in various scientific and engineering \textcolor{black}{field}. It finds relevance in fields such as physical biology~\citep{FENNEL2015153}, marine science~\citep{WARNER2007355}, earth and atmospheric sciences~\citep{Earth}, as well as mantle dynamics~\citep{RICARD200731}. This article primarily focuses on investigating the dynamics of the unsteady advection-diffusion equation (ADE) under various boundary constraints. The following equation mathematically describes the ADE:
\begin{equation}\label{eq0101}
\frac{\partial u}{\partial t} = \mathbf{div}\big{(}\bm{p} \cdot\nabla u\big{)}-\bm{q}\cdot\nabla u + f.
\end{equation}
This equation captures the interaction between convection and diffusion at different temporal and spatial scales. It involves a scalar variable denoted as $u$, which is transported through advection and diffusion. The constant or vector parameters $\bm{p}$ and $\bm{q}$ represent the advection field's speed and the diffusion coefficient in different directions, respectively. The term $f$ signifies the internal source or sink's capacity, while the concentration gradient is represented by $\nabla u$, where $\nabla$ denotes the gradient operator and $\mathbf{div}$ denotes the divergence operator.
Similar to other types of partial differential equations(PDEs), the analytical solutions for ADE are generally seldom available, so solving these PDEs numerically using approximation methods is necessary. Numerical methods such as finite element method (FEM)~\cite{burman2007continuous, zheng2010note, dhawan2012numerical, bachini2021intrinsic, cier2021automatically}, finite difference method (FDM)~\cite{sari2010high, prieto2011application, shankar2018hyperviscosity} and finite volume method (FVM)~\cite{ollivier2002high} are commonly used to solve ADE. In these approaches, the computational domain of interest is divided into a set of simple regular mesh, and the solution is computed in these mesh patches. Generally, to reduce the numerical error, the size of the mesh is required to be small when solving PDEs by these mesh-dependent methods, it will yield significant computational and storage challenges~\citep{boztosun2002analysis}. Given a specific mesh size, several numerical techniques have been developed to reduce the errors of mesh methods~\cite{ewing2001summary}, such as the upwind scheme~\cite{brooks1982streamline} and Galerkin least squares strategy~\cite{hughes1989new}. However, mesh-dependent methods can be challenging, time-consuming, and computationally expensive when dealing with complex domains of interest and boundary constraints. In contrast, meshless methods that use a set of configuration points without grids have been developed to approximate the solution of ADE, such as Radial Basis Function~\cite{askari2017numerical, zhang2022simulation}, Monte Carlo methods~\cite{erhel2014combined}, and B-spline collocation approaches~\cite{jena2021computational}. While these methods are easy to implement and straightforward, their accuracy may deteriorate compared to grid-based methods.

In the past few years, deep neural networks (DNN) have demonstrated significant potential in solving ordinary and partial differential equations as well as inverse problems. This is due to their ability to handle strong nonlinearity and high-dimensional problems, as highlighted in various studies~\cite{Ritz, BERG201828, RAISSI2019686, sirignano2018dgm, khoo2019switchnet, lyu2022mim, jagtap2020conservative, jagtap2021extended}.
This methodology is preferable because it transforms a PDE problem into an optimization problem, and then approximates the PDE solution through gradient backpropagation and automatic differentiation of the DNN. Furthermore, they are inherently mesh-free and can address high-dimensional and geometrically complex problems more efficiently than mesh-based methods.
 
\textcolor{black}{Physics Informed Neural Networks (PINN) were first proposed by \citet{RAISSI2019686} in 2019 that derived from the concept of physics-constrained learning applied to solve conventional differential equations in the early 1990s~\cite{LEE1990110, IsaacE01, IsaacE02}}, which incorporates physical laws into neural networks by adding the residuals of both the PDEs and the boundary conditions (BCs) as multiple terms in the loss function. \textcolor{black}{After that, lots of efforts have been made to further enhance the performance of PINN, such as Parallel PINN~\cite{shukla2021parallel},  APINN~\cite{hu2023augmented} and Unified scalable PINN~\cite{penwarden2023unified}.} \textcolor{black}{And the PINN is widely applied to address various types of PDEs including material identification~\cite{shukla2021physics}, the optimal location of sensors~\cite{jagtap2022deep} and flows problems~\cite{jagtap2022physics,mao2020physics}.} \textcolor{black}{Some theory of convergence and the analysis of generalization error for PINN are also made by ~\citet{mishra2022estimates, de2023error,hu2022extended}}. In addition, the Deep Ritz Method was proposed by \citet{Ritz} for numerically solving variational problems, and various PINN methods have been proposed for solving PDEs with complex boundaries by researchers such as~\citet{Wang_2020, Gao_2021, Sheng_2021}. Some researchers have also investigated using physical constraints to train PINN~\cite{Rishi11374, Mohammad, Kailai, Jonathan}.
Despite some impressive developments in the field, including DPM~\cite{Jungeun}, PINN still faces significant challenges, as pointed out in the study~\cite{Aditi}. Several theoretical papers~\cite{Sifan1, Sifan2} have highlighted an imbalanced competition between the terms of PDEs and BCs in the loss function, limiting PINN's applicability in complex geometric domains. To address this problem, researchers such as \citet{BERG201828, Sun_2020, Lulu1} have proposed incorporating BCs into the ansatz such that any instance from the ansatz automatically satisfies the BCs, resulting in so-called hard-constraint methods. By satisfying the BCs precisely, the PDE solution becomes an unconstrained problem, which can be more effectively trained in a neural network.

In addition, the standard PINN method may not perform well on problems with high-frequency components because of the low-frequency bias of DNN, which has been documented by~\citet{Xu_2020} and \citet{rahaman2018spectral}. For instance, when using a DNN to fit the function $y=\sin (x) + \sin(2x)$, the DNN output initially approximates $\sin (x)$ and then gradually converges to $\sin (x) + \sin (2x)$. This phenomenon is called the Frequency Principle (F-Principle) or spectral bias, which stems from the inherent divergence in the gradient-based DNN training process. To address this issue, researchers have explored the relationship between DNN and Fourier analysis, which was inspired by the F-Principle~\cite{xu_training_2018, Matthew2020Fourier, wang2020eigenvector}. Recent experimental studies by~\citet{Ellen1, BenNERF} have suggested that a heuristic sinusoidal encoding of input coordinates, termed ``positional encoding'', can enable PINN to capture higher frequency content. To enhance the computational efficiency of the aforementioned PINN, we construct a separate subnetwork for each frequency to capture signals at different frequencies.

Despite many researchers having discovered the F-Principle phenomenon and the unbalanced rivalry between the terms of PDEs and BCs, there is still no universal solution for both issues, which presents an opportunity for further research. Additionally, there is an opportunity to explore the Neumann boundary in more depth, as the majority of research on ADE has focused on the Dirichlet boundary, while the Neumann boundary is less studied and more intricate.

This paper introduces a novel approach called sub-Fourier Hard-constraint PINN (SFHCPINN) to solve a class of ADEs using Fourier analysis and the hard constraint technique. The SFHCPINN approach uses hard-constraint PINN to enforce initial and boundary conditions and incorporates the residual of the governed equation for ADE into the cost function to guide the training process. This transforms the solution of ADE into an unconstrained optimization problem that can be efficiently solved using gradient-based optimization techniques. To further reduce the approximation error of the DNN, a subnetworks framework of DNN with sine and cosine as the activation function is introduced by the F-Principle and Fourier theory. While previous studies primarily focused on simulated data using mean squared error (MSE) loss, the primary contributions of this paper are:

\begin{enumerate}
\item[a.] Our proposed method involves a PINN with a subnetwork architecture and a Fourier-based activation function. This approach aims to address the issue of gradient leakage in DNN parameters by leveraging the F-principle and Fourier theorem.
\item[b.] Our approach involves the use of a structured PINN with hard constraints, which allows for the automatic satisfaction of initial and boundary conditions. This approach enables the accurate resolution of ADEs with complex boundary conditions, including the Dirichlet boundary, Neumann boundary, and mixed form.
\item[c.] We provide compelling evidence for the efficacy of the proposed method by demonstrating the superiority of hard-constraint PINN over soft-constraint PINN in solving a class of ADEs under both Dirichlet and Neumann boundaries.
\end{enumerate}

The structure of this paper is outlined as follows. Section~\ref{sec:Preliminaries} presents an overview of deep neural networks and the standard PINN framework for PDEs. In Section~\ref{sec:SFHCPINN}, we introduce a novel approach to solving ADE using PINN with hard constraints, as well as an activation function based on Fourier analysis. Section~\ref{subsec:algr2SFHCPINN} details the SFHCPINN algorithm for approximating the solution of ADE. In Section~\ref{sec:experiment}, we provide numerical examples to demonstrate the effectiveness of the proposed method for ADE. Finally, we present the conclusions of the paper in Section~\ref{sec:conclusion}.

\section{Preliminaries}\label{sec:Preliminaries}
This section has presented a detailed exposition of the relevant mathematical principles and formulae about DNN and PINN.

\subsection{Deep Neural Networks}\label{sec:0201}
Initially, we present the standard neural cell of Deep Neural Networks (DNN) and the mapping relationship between input $\bm{x}\in\mathbb{R}^{d}$ and an output $\bm{y}\in\mathbb{R}^{m}$, as expressed by \eqref{eq0201}:
\begin{equation}\label{eq0201}
\bm{y}=\sigma(\bm{W} \bm{x}+\bm{b}).
\end{equation}
Here, the activation function $\sigma(\cdot)$ is an element-wise non-linear model, \textcolor{black}{$\bm{W} = (w_{ij}) \in \mathbb{R}^{m\times d}$} and $\bm{b}\in\mathbb{R}^{m}$ are the weight matrix and bias vector, respectively.
The standard unit \eqref{eq0201} is usually known as a hidden layer, and its output is fed into another activation function after modification with a new weight and bias.
Hence, a DNN is constructed with stacked linear and nonlinear activation functions.
The mathematical expression for a DNN with input data $\bm{x}\in\mathbb{R}^{d}$ can be formulated as:
\begin{equation}\label{form2dnn}
	\begin{cases}
		\bm{y}^{[0]} = \bm{x}\\
		\bm{y}^{[\ell]} = \sigma\circ(\bm{W}^{[\ell]}\bm{y}^{[\ell-1]}+\bm{b}^{[\ell]}), ~~\text{for}~~\ell =1, 2, 3, \cdots\cdots, L
	\end{cases}.
\end{equation}
Here, $\bm{W}^{[\ell]} \in \mathbb{R}^{n_{\ell+1}\times n_{\ell}}$ and $\bm{b}^{[\ell]}\in\mathbb{R}^{n_{\ell+1}}$ denote the weights and biases of the $\ell$-th hidden layer, respectively. $n_0=d$ and $n_{L+1}$ is the dimension of output. The notation ``$\circ$'' indicates an element-wise operation. The parameter set of $\bm{W}^{[1]},\cdots \bm{W}^{[L]}, \bm{b}^{[1]},\cdots \bm{b}^{[L]}$ is represented by $\bm{\theta}$, and the DNN output is denoted by $\bm{y}(\bm{x};\bm{\theta})$.

\subsection{Physics Informed Neural Networks}\label{sec:0202}

Let us consider a system of parametrized PDEs given by:
\begin{equation}\label{eq2PPDE}
\begin{aligned}
&\mathcal{N}_{\bm{\lambda}}[\hat{u}(\bm{x},t)]=\hat{f}(\bm{x},t), ~\quad \bm{x} \in \Omega, t \in[t_0, T] \\
&\mathcal{B}\hat{u}\left(\bm{x}, t\right)=\hat{g}(\bm{x}, t), ~\quad\quad\bm{x} \in \partial \Omega, t \in[t_0, T]\\
&\hat{u}(\bm{x}, t_0)=\hat{h}(\bm{x}, t_0), \quad \quad \quad\bm{x} \in \Omega ,
\end{aligned}
\end{equation}
in which $\mathcal{N}_{\bm{\lambda}}$ stands for the linear or nonlinear differential operator with parameters $\bm{\lambda}$, $\mathcal{B}$ is the boundary operators. $\Omega $ and  $\partial \Omega $ respectively \textcolor{black}{illustrate the domain of interest and its boundary}. In general PINN, one can substitute a DNN model for the solution of PDEs \eqref{eq2PPDE}, then obtain the optimal solution by minimizing the following loss function:
\begin{equation}\label{loss2PINN}
L=Loss_{PDE} +\omega_{1}Loss_{IC}+\omega_{2} Loss_{BC}
\end{equation}
\textcolor{black}{with}
\begin{equation}\label{subloss2PINN}
\begin{aligned}
&Loss_{PDE} = \frac{1}{N_P}\sum_{i=1}^{N_P}\left| \mathcal{N}_{\bm{\lambda}}[\hat{u}_{NN}(\bm{x}^i,t^i)]-\hat{f}(\bm{x}^i,t^i)\right|^2 \\
&Loss_{BC} = \frac{1}{N_B}\sum_{i=1}^{N_B}\bigg{|}\mathcal{B}\hat{u}_{NN}\left(\bm{x}^i, t^i\right)-\hat{g}(\bm{x}^i,t^i)\bigg{|}^2\\
&Loss_{IC} = \frac{1}{N_I}\sum_{i=1}^{N_I}\bigg{|}\hat{u}_{NN}(\bm{x}^i, t_0)-\hat{h}(\bm{x}^i,t_0)\bigg{|}^2
\end{aligned}
\end{equation}
where $\omega_{1}$ and $\omega_2$ are the loss weighting coefficients on distinct \textcolor{black}{boundaries}. $Loss_{PDE}$, $Loss_{IC}$, and $Loss_{BC}$ depict the residual of governing equations, the loss of the given initial condition, and the loss of the prescribed BC, respectively. In addition, if some data are available inside the domain, an extra loss term indicating the mismatch between the predictions and the data can be taken into account
\begin{equation}
    Loss_{D} = \frac{1}{N_D}\sum_{i=1}^{N_D}\bigg{|}u_{NN}(\bm{x}^i, t^i)-u_{Data}^i\bigg{|}^2.
\end{equation}
The structure of PINN for solving parametrized PDEs \eqref{eq2PPDE} is depicted in the following Figure~\ref{fig: PINN Structure}.
\begin{figure}[!htbp]
    \centering
    \includegraphics[scale=0.4]{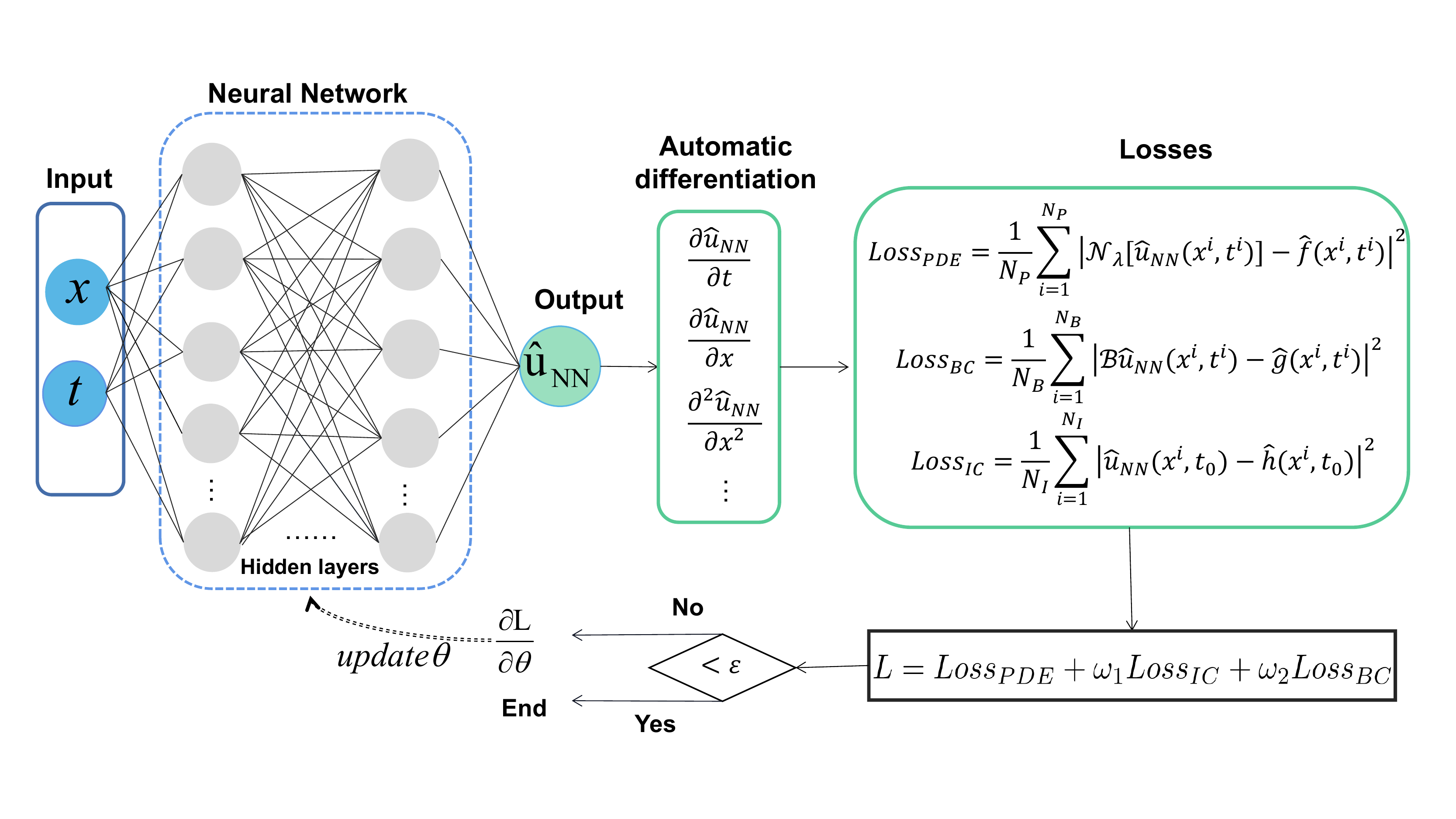}
    \caption {Schematic diagram of physical information neural network (PINN).}
    \label{fig: PINN Structure}
\end{figure}

\section{Fourier induced Subnetworks Hard Constraint PINN to Advection-Diffusion Equation}\label{sec:SFHCPINN}

\subsection{Unified architecture of Hard Constraint PINN to Advection-Diffusion Equation}\label{subsec:HCPINN}
In this section, we now consider the following advection-diffusion equation with prescribed boundary and initial conditions, it is
\begin{equation}\label{eq:ADE_IC_BC}
\begin{aligned}
    &\frac{\partial u(\bm{x}, t)}{\partial t} = \mathbf{div}\big{(}\bm{p}\cdot \nabla u(\bm{x}, t)\big{)}  -\bm{q}\cdot \nabla u(\bm{x}, t) + f(\bm{x}, t), \quad \bm{x} \in \Omega, t \in[t_0, T] \\
    &\mathcal{B}u(\bm{x}, t)=g(\bm{x}, t), \quad ~\quad\quad\bm{x} \in \partial \Omega, t \in[t_0, T]\\
    & u(\bm{x}, t_0) = h(\bm{x}, t_0), \quad~~\quad \quad\bm{x} \in \Omega
\end{aligned}
\end{equation}
where $\Omega$ is a bounded subset of $\mathbb{R}^{d}$ with piecewise Lipschitz boundary which satisfies the interior cone condition and $ \partial \Omega$ represents the boundary of the interested domain. Generally, the boundary is a complicated geometry and composed of essential and natural boundaries, i.e., $\partial \Omega=\Gamma_D \cup \Gamma_N$ and $\Gamma_D \cap \Gamma_N=\emptyset$. The operational item $\mathcal{B}$ indicates the BC, such as the Dirichlet, Neumann, or Robin boundary. 

When a standard PINN trial function $u_{NN}$ is used to approximate the solution of ADE~\eqref{eq:ADE_IC_BC}, it will be firstly differential concerning variable $\bm{x}$ and $t$, respectively, then embed into the residual of the governed equation and constructed the main part of loss function of neural networks. Under the imposed boundary and initial conditions, the neural network solution with parameter $\bm{\theta}$ can be obtained by minimizing
\begin{equation}\label{softloss2ADE}
    \textcolor{black}{Loss(S_R,S_B,S_I;\bm{\theta}) = Loss_{R}(S_R;\bm{\theta}) +  \frac{\gamma}{N_B}\sum^{N_B}_{i=1}\bigg{|}\mathcal{B}u_{NN}(\bm{x}_B^i,t_B^i)-g(\bm{x}_B^i,t_B^i)\bigg{|}^2 + \frac{\omega}{N_I}\sum^{N_I}_{i=1}\bigg{|}\mathcal{I}u_{NN}(\bm{x}_I^i,t_0)-h(\bm{x}_I^i)\bigg{|}^2}
\end{equation}
for $(\bm{x}_B^i, t_B^i)\in S_B$ and $(\bm{x}_I^i,t_0)\in S_I$, as well as  
\begin{equation}\label{loss_in_sr}
  \textcolor{black}{Loss_{R}(S_R;\bm{\theta}) = \frac{1}{N_R}\sum_{i=1}^{N_R}\left| \frac{\partial u_{NN}(\bm{x}_R^i, t_R^i)}{\partial t} - \mathbf{div}\big{(}\bm{p}\cdot \nabla u_{NN}(\bm{x}_R^i, t_R^i)\big{)}  + \bm{q}\cdot \nabla u_{NN}(\bm{x}_R^i, t_R^i) - f(\bm{x}_R^i, t_R^i) \right|^2}
\end{equation}
 for $(\bm{x}^R_i, t^R_i)\in S_R$. Here and hereinafter, $S_R=\{(\bm{x}^R_i, t^R_i)\}_{i=1}^{N_R}$, $S_B=\{(\bm{x}^B_i, t^B_i)\}_{i=1}^{N_B}$ and $S_I=\{(\bm{x}_I^i,t_0)\}_{i=1}^{N_I}$ stand for the sets of distributed sample points on $\Omega\times T$, $\partial \Omega\times T$ and $\Omega\times\{t_0\}$, respectively. In addition, two penalty parameters $\gamma$ and $\omega$ are introduced to control the contributions of boundary and initial for loss function.

\textcolor{black}{Many scholars have studied carefully the choice of the residual term in loss function}, for example, Mixed PINN~\cite{lyu2022mim}, XPINN~\cite{jagtap2021extended}, cPINN~\cite{jagtap2020conservative}, two-stage PINN~\cite{lin2022two} and gPINN~\cite{yu2022gradient}. Another problem to be addressed is how to enforce the initial and boundary conditions (I/BCs). The imposition of I/BCs is crucial for solving PDEs because it allows a unique solution. Considering the optimization nature of the PINN, the primitive way of applying I/BCs is to penalize the discrepancy of initial and boundary constraints for PDEs in a soft manner.

In PINN-based deep collocation methods, the performance of optimization depends on the relevance of each term. However, assigning the weights of each term may be difficult, then the approximations of I/BCs may not be favorable, resulting in an unsatisfactory solution. Consequently, we may apply the boundary constraint in a ``hard'' manner by including particular solutions that satisfy the I/BCs. Consequently, the constraints on the boundaries are gently met. \textcolor{black}{It will improve the capacity of neural networks for dealing with complex geometry issues~\cite{BERG201828, Sun_2020, Lulu1}.} The \textcolor{black}{boundary} operator $\mathcal{B}$ in (\ref{eq:ADE_IC_BC}) will be used in the following to investigate PDEs with Dirichlet BCs. In this instance, our proposed theory for the solution is
\begin{equation}\label{eq:ansatz}
    u_{NN}(\bm{x},t)=G(\bm{x},t)+D(\bm{x},t) NN^{L}(\bm{x},t; \bm{\theta})
\end{equation}
where $NN^L$ is a fully-connected deep neural network, $G(\bm{x},t)$ is \textcolor{black}{an extension function} meeting the I/BCs constraints $\mathcal{B}u(\bm{x}, t)$ and $\mathcal{I}u(\bm{x}, t_0)$, and $D(\bm{x},t)$ is a smooth distance function giving the distance from $(\bm{x},t) \in \Omega \times T$ to $\partial \Omega \times \{t_0\}$. The objective of this concept is to compel the approximate solution to conform to a set of restrictions, notably the Dirichlet BCs.
In other words, while $\bm{x}$ is on the $\partial \Omega \times \{t_0\}$ \textcolor{black}{boundary}, $D(\bm{x},t)$ equals zero, and the value increases as points depart from the I/BCs.

It is worth noting that \textcolor{black}{$u_{NN}(\bm{x},t)$}reaches its I/BCs value at the hypothesis equation's boundary point (\ref{eq:ansatz}). For those ADE problems with simple IC/BC and an easy form of $\partial \Omega$, $D(\bm{x},t)$ could be defined analytically. 
Nevertheless, assuming the geometry is too complicated for an analytic formulation, both the extension boundary function $G(\bm{x},t)$ and the smooth function $D(\bm{x},t)$ may be parameterized using small-scale NNs according to the given I/BCs constraints and a small of configuration points sampled from the interested domain with boundary. Therefore, it will not add any more complexity when optimizing the loss of hard-constraint PINN (HCPINN):
\begin{equation}\label{softlossHCPINNDir}
    {Loss}_{HCPINN}(S_R,S_B,S_I;\bm{\theta})= Loss_R({S}_{R};\bm{\theta}).
\end{equation}
\textcolor{black}{with $Loss_{R}({S}_R;\bm{\theta})$ being defiend as in \eqref{loss_in_sr}.}

Considering $\partial \Omega$ as the Neumann-boundary and $\mathcal{B}$ as the differential operator. 
In contrast to Dirichlet BC, which is stored inside a particular solution, Neumann BC is included throughout the equation loss.
The ansatz solution for the Neumann BCs is the same as (\ref{eq:ansatz}) and $G(\bm{x},t)$ now is an extension meeting the IC constraints $\mathcal{I}u(\bm{x}, t_0)$ and $D(\bm{x},t)$ is a smooth distance function giving the distance from $(\bm{x},t) \in \Omega \times T$  to $\Omega \times \{t_0\}$. The Neumann BCs are encoded into the loss function:
\begin{equation}\label{softlossHCPINNNeu}
    {Loss}_{HCPINN}(S_R,S_B,S_I;\bm{\theta}) = Loss_{R}({S}_R;\bm{\theta}) +  \frac{\gamma}{N_B}\sum^{N_B}_{i=1}\bigg{|}\mathcal{B}u_{NN}(\bm{x}^B_i,t^B_i)-g(\bm{x}^B_i,t^B_i)\bigg{|}^2.
\end{equation}

We then conclude the unified loss function of the HCPINN as follows:
\begin{equation}
    \bm{\theta}^{*} =\underset{\bm{\theta}}{\arg \min } {Loss}_{HCPINN}(\bm{\theta}).
\end{equation}
To obtain the $\bm{\theta}^*$, one can update the parameters $\bm{\theta}$ using gradient descent method overall training samples or a few training samples at each iteration. In particular, Stochastic gradient descent (SGD) \textcolor{black}{and its improved versions, such as Adam~\cite{Adam}, Adagrad~\cite{duchi2011adaptive} and RMSprop~\cite{hinton2012neural}, are} the common optimization technique for deep learning. In the implementation, the SGD method requires only one of $n$ function evaluations at each iteration compared with the gradient descent method. Additionally, instead of picking one term, one can also choose a ``mini-batch'' of terms at each step. In this context, the \textcolor{black}{update scheme of parameters of DNN for vanilla} SGD is given by:
\begin{equation*}
\bm{\theta}^{k+1}=\bm{\theta}^{k}-\alpha_k\nabla_{\bm{\theta}^k}Loss_{HCPINN}(\bm{x};\bm{\theta}^{k}),~~\bm{x}\in S_R~\text{or}~\bm{x}\in S_R~\cup S_B,
\end{equation*}
where the ``learning rate'' $\alpha_k$ decreases with $k$ increasing. 

\begin{remark}
    To compute the smooth distance function $D(\bm{x},t)$ in the Dirichlet condition, we first calculate the non-smooth distance function $d$ and estimate it using a low-capacity NN. At each point $(\bm{x},t)$, we define $d$ as the shortest distance to a boundary point at which a BC must be applied. Indeed,
    \begin{equation}
        d(\bm{x},t)=\min _{(\bm{x},t)^{*} \in \partial\Omega \times \{t_0\}}\left\|(\bm{x},t)-(\bm{x},t)^{*}\right\| .
    \end{equation}
    The exact form of $d$ (and $D$) is not important other than that $D$ is smooth and
    \begin{equation}
        |D(\bm{x},t)|<\epsilon, \quad \forall (\bm{x},t) \in \partial\Omega \times \{t_0\}.
    \end{equation}
    We can use a small subset from $\partial \Omega \times \{t_0\}$ to compute $d$.
\end{remark}

\begin{remark}
    Instead of computing the actual distance function, we could use the more extreme version
    \begin{equation}
        d(\bm{x},t)=\left\{\begin{array}{ll}
          0, & (\bm{x},t) \in \partial\Omega \times \{t_0\}\\
           1, & \text { otherwise }
    \end{array}\right..
    \end{equation}
 Moreover, for issues where the I/BCs are enforced in simple geometry, $D(\bm{x},t)$ and $G(\bm{x},t)$ may be derived analytically~\cite{Lagri,IsaacE01, sukumar2022exact, schiassi2021extreme}. For instance, we can define $D(\bm{x},t)=t/T$ when there is simply initial boundary enforcement on $t=0$, and we can choose $D(\bm{x},t)=(x-a)(b-x)$ or $(1-e^{a-\bm{x}})(1-e^{\bm{x}-b})$ in $\Omega=[a,b]$ when the BCs are only imposed on $\partial\Omega$. For complex cases, it is difficult to identify an analytical formula for $D(\bm{x},t)$, but it is possible to approximate it using spline functions~\cite{Sheng_2021}.
\end{remark}
\begin{remark}
    The ansatz \eqref{eq:ansatz} demands that $G$ be globally defined and smooth, as well as that
    \begin{equation}
        |G(\bm{x},t)-g(\bm{x},t)|<\epsilon, \quad \forall (\bm{x},t) \in \partial\Omega \times \{t_0\}
    \end{equation}
where $g(\bm{x},t)$ is the function satisfying I/BCs \textcolor{black}{of given PDEs}. To compute $G$ we simply train an NN to fit $G(\bm{x},t)$,$\quad \forall (\bm{x},t) \in \partial\Omega \times \{t_0\}$. The loss function used is given by
\begin{equation}\label{lossG}
    \textcolor{black}{{Loss}_{G} = \frac{1}{N_G}\sum^{N_G}_{i=1}\bigg{|}{G}(\bm{x}^i,t^i)-g(\bm{x}^i,t^i)\bigg{|}^2,~~(\bm{x}^i,t^i)\in\Omega\times[t_0, T]}
\end{equation}
and apply SGD as the optimization technique as well.
\end{remark}

\begin{remark}
    In some given BCs, $G(\bm{x},t)$ could be defined directly utilizing the I/BCs. For example, if \eqref{eq:ADE_IC_BC} $\mathcal{B}u(\bm{x}, t)=0$ and $ \mathcal{I}u(\bm{x}, t_0)=0$, then we could directly define $G(\bm{x},t)=0.$
\end{remark}

\subsection{Sub-Fourier PINN and its activation function}\label{subsec:SFPINN}
The activation function is one of the critical issues for designing the architecture of DNN. As a non-linear transformation that bounds the value for given input data, it directly affects the performance of DNN models in practical applications.
Several different types of activation functions have been used in DNN, such as $\text{ReLU}(\bm{z}) = \max\{0,\bm{z}\}$ and $\tanh(\bm{z})$.

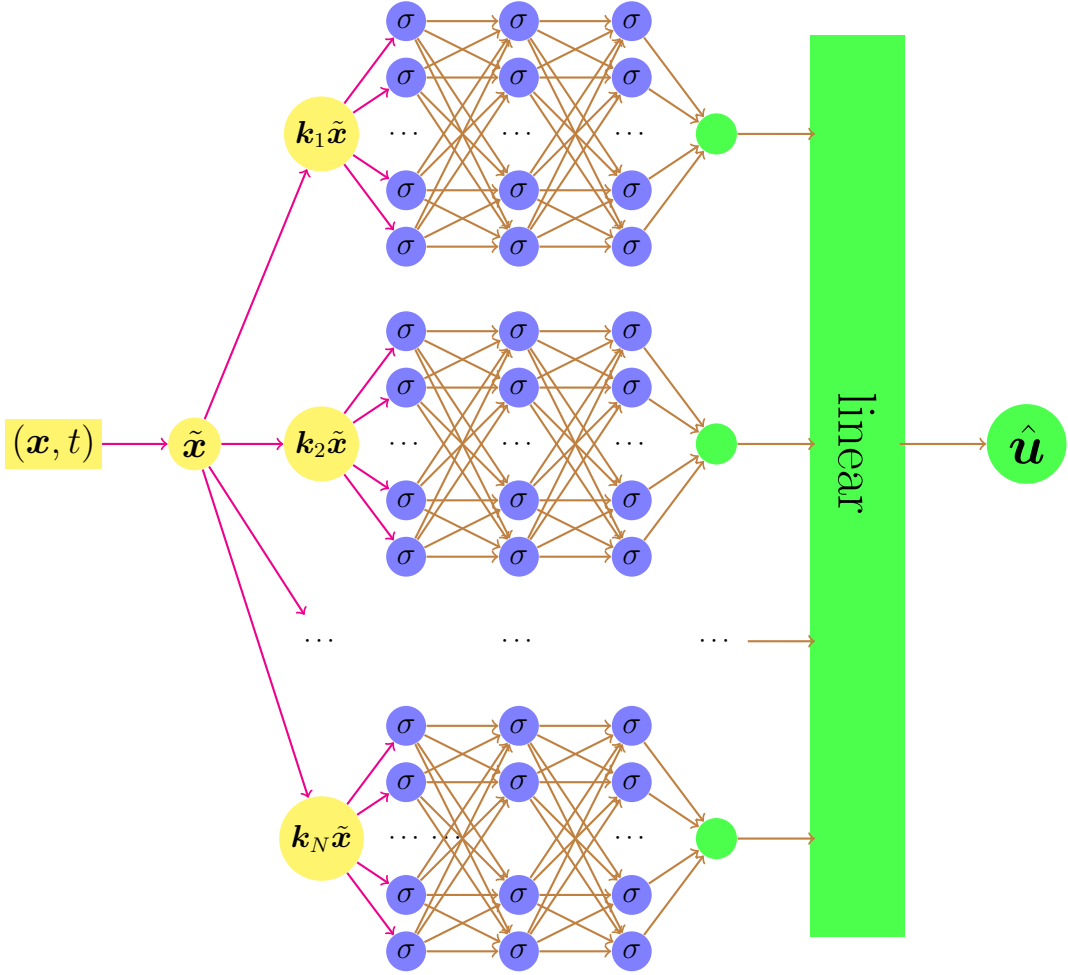
\begin{figure}[!ht]
    \begin{center}
    \begin{tikzpicture}[scale=0.75]	
        \node[rectangle, fill=yellow!70,inner sep=2.25pt](in2pre) at (-4.75,8){\Large$(\bm{x},t)$};
        \node[circle, fill=yellow!70,inner sep=2.25pt](in) at (-2.25,8){\Large$\tilde{\bm{x}}$};
        \node[circle, fill=yellow!70,inner sep=2.5pt] (xh) at  (0.0,1){\large$\bm{k}_N\tilde{\bm{x}}$};
        \draw[line width=0.8pt,color=magenta,->] (in2pre) -- (in);
        \draw[line width=0.8pt,color=magenta,->] (in) -- (xh);
        
        \node[circle, fill=blue!50,inner sep=2.25pt] (h01) at (1.5,-1){\large$\sigma$};
        \node[circle, fill=blue!50,inner sep=2.25pt] (h02) at (1.5,0){\large$\sigma$};
        \node[circle, fill=blue!0] (h03) at (1.5,1){$\cdots$};
        \node[circle, fill=blue!50,inner sep=2.25pt] (h04) at (1.5,2){\large$\sigma$};
        \node[circle, fill=blue!50,inner sep=2.25pt] (h05) at (1.5,3){\large$\sigma$};
        
        \draw[line width=0.8pt,color=magenta,->] (xh) -- (h01);
        \draw[line width=0.8pt,color=magenta,->] (xh) -- (h02);
        \draw[line width=0.8pt,color=magenta,->] (xh) -- (h04);
        \draw[line width=0.8pt,color=magenta,->] (xh) -- (h05);
        
        \node[circle, fill=blue!50,inner sep=2.25pt] (h11) at (3.5,-1.0){\large$\sigma$};
        \node[circle, fill=blue!50,inner sep=2.25pt] (h12) at (3.5,0){\large$\sigma$};
        \node[circle, fill=blue!0] (h13) at (2.25,1){$\cdots$};
        \node[circle, fill=blue!50,inner sep=2.25pt] (h14) at (3.5,2){\large$\sigma$};
        \node[circle, fill=blue!50,inner sep=2.25pt] (h15) at (3.5,3){\large$\sigma$};
        
        \draw[line width=0.8pt,color=brown,->] (h01) -- (h11);
        \draw[line width=0.8pt,color=brown,->] (h01) -- (h12);
        \draw[line width=0.8pt,color=brown,->] (h01) -- (h14);
        \draw[line width=0.8pt,color=brown,->] (h01) -- (h15);
        
        \draw[line width=0.8pt,color=brown,->] (h02) -- (h11);
        \draw[line width=0.8pt,color=brown,->] (h02) -- (h12);
        \draw[line width=0.8pt,color=brown,->] (h02) -- (h14);
        \draw[line width=0.8pt,color=brown,->] (h02) -- (h15);
        
        \draw[line width=0.8pt,color=brown,->] (h04) -- (h11);
        \draw[line width=0.8pt,color=brown,->] (h04) -- (h12);
        \draw[line width=0.8pt,color=brown,->] (h04) -- (h14);
        \draw[line width=0.8pt,color=brown,->] (h04) -- (h15);
        
        \draw[line width=0.8pt,color=brown,->] (h05) -- (h11);
        \draw[line width=0.8pt,color=brown,->] (h05) -- (h12);
        \draw[line width=0.8pt,color=brown,->] (h05) -- (h14);
        \draw[line width=0.8pt,color=brown,->] (h05) -- (h15);
        
        \node[circle, fill=blue!50,inner sep=2.25pt] (h21) at (5.5,-1.0){\large$\sigma$};
        \node[circle, fill=blue!50,inner sep=2.25pt] (h22) at (5.5,0){\large$\sigma$};
        \node[circle, fill=blue!0] (h23) at (5.5,1){$\cdots$};
        \node[circle, fill=blue!50,inner sep=2.25pt] (h24) at (5.5,2){\large$\sigma$};
        \node[circle, fill=blue!50,inner sep=2.25pt] (h25) at (5.5,3){\large$\sigma$};
        \node[circle, fill=green!70,inner sep=5.5pt] (uh) at (7.0,1){};
        
        \draw[line width=0.8pt,color=brown,->] (h11) -- (h21);
        \draw[line width=0.8pt,color=brown,->] (h11) -- (h22);
        \draw[line width=0.8pt,color=brown,->] (h11) -- (h24);
        \draw[line width=0.8pt,color=brown,->] (h11) -- (h25);
        
        \draw[line width=0.8pt,color=brown,->] (h12) -- (h21);
        \draw[line width=0.8pt,color=brown,->] (h12) -- (h22);
        \draw[line width=0.8pt,color=brown,->] (h12) -- (h24);
        \draw[line width=0.8pt,color=brown,->] (h12) -- (h25);
        
        \draw[line width=0.8pt,color=brown,->] (h14) -- (h21);
        \draw[line width=0.8pt,color=brown,->] (h14) -- (h22);
        \draw[line width=0.8pt,color=brown,->] (h14) -- (h24);
        \draw[line width=0.8pt,color=brown,->] (h14) -- (h25);
        
        \draw[line width=0.8pt,color=brown,->] (h15) -- (h21);
        \draw[line width=0.8pt,color=brown,->] (h15) -- (h22);
        \draw[line width=0.8pt,color=brown,->] (h15) -- (h24);
        \draw[line width=0.8pt,color=brown,->] (h15) -- (h25);
        
        \draw[line width=0.8pt,color=brown,->] (h21) -- (uh);
        \draw[line width=0.8pt,color=brown,->] (h22) -- (uh);
        \draw[line width=0.8pt,color=brown,->] (h24) -- (uh);
        \draw[line width=0.8pt,color=brown,->] (h25) -- (uh);

        \node[circle, fill=blue!0,inner sep=3.5pt] (xi) at (0.0,4.5){$\cdots$};
	\node[circle, fill=blue!0,inner sep=3.5pt] (ihidden) at (3.5,4.5){$\cdots$};
	\node[circle, fill=blue!0,inner sep=3.5pt] (ui) at (7.0,4.5){$\cdots$};
        \draw[line width=0.8pt,color=magenta,->] (in) -- (xi);

        \node[circle, fill=yellow!70,inner sep=2.25pt] (xj) at  (0.0,8){\large$\bm{k}_2\tilde{\bm{x}}$};
        \draw[line width=0.8pt,color=magenta,->] (in) -- (xj);

        \node[circle, fill=blue!50,inner sep=2.25pt] (j01) at (1.5,6){\large$\sigma$};
        \node[circle, fill=blue!50,inner sep=2.25pt] (j02) at (1.5,7){\large$\sigma$};
        \node[circle, fill=blue!0] (j03) at (1.5,8){$\cdots$};
        \node[circle, fill=blue!50,inner sep=2.25pt] (j04) at (1.5,9){\large$\sigma$};
        \node[circle, fill=blue!50,inner sep=2.25pt] (j05) at (1.5,10){\large$\sigma$};

        \draw[line width=0.8pt,color=magenta,->] (xj) -- (j01);
        \draw[line width=0.8pt,color=magenta,->] (xj) -- (j02);
        \draw[line width=0.8pt,color=magenta,->] (xj) -- (j04);
        \draw[line width=0.8pt,color=magenta,->] (xj) -- (j05);

        \node[circle, fill=blue!50,inner sep=2.25pt] (j11) at (3.5,6){\large$\sigma$};
        \node[circle, fill=blue!50,inner sep=2.25pt] (j12) at (3.5,7){\large$\sigma$};
        \node[circle, fill=blue!0] (j13) at (3.5,8){$\cdots$};
        \node[circle, fill=blue!50,inner sep=2.25pt] (j14) at (3.5,9){\large$\sigma$};
        \node[circle, fill=blue!50,inner sep=2.25pt] (j15) at (3.5,10){\large$\sigma$};

        \draw[line width=0.8pt,color=brown,->] (j01) -- (j11);
        \draw[line width=0.8pt,color=brown,->] (j01) -- (j12);
        \draw[line width=0.8pt,color=brown,->] (j01) -- (j14);
        \draw[line width=0.8pt,color=brown,->] (j01) -- (j15);
        
        \draw[line width=0.8pt,color=brown,->] (j02) -- (j11);
        \draw[line width=0.8pt,color=brown,->] (j02) -- (j12);
        \draw[line width=0.8pt,color=brown,->] (j02) -- (j14);
        \draw[line width=0.8pt,color=brown,->] (j02) -- (j15);
        
        \draw[line width=0.8pt,color=brown,->] (j04) -- (j11);
        \draw[line width=0.8pt,color=brown,->] (j04) -- (j12);
        \draw[line width=0.8pt,color=brown,->] (j04) -- (j14);
        \draw[line width=0.8pt,color=brown,->] (j04) -- (j15);
        
        \draw[line width=0.8pt,color=brown,->] (j05) -- (j11);
        \draw[line width=0.8pt,color=brown,->] (j05) -- (j12);
        \draw[line width=0.8pt,color=brown,->] (j05) -- (j14);
        \draw[line width=0.8pt,color=brown,->] (j05) -- (j15);

        \node[circle, fill=blue!50,inner sep=2.25pt] (j21) at (5.5,6){\large$\sigma$};
        \node[circle, fill=blue!50,inner sep=2.25pt] (j22) at (5.5,7){\large$\sigma$};
        \node[circle, fill=blue!0] (j23) at (5.5,8){$\cdots$};
        \node[circle, fill=blue!50,inner sep=2.25pt] (j24) at (5.5,9){\large$\sigma$};
        \node[circle, fill=blue!50,inner sep=2.25pt] (j25) at (5.5,10){\large$\sigma$};

        \draw[line width=0.8pt,color=brown,->] (j11) -- (j21);
        \draw[line width=0.8pt,color=brown,->] (j11) -- (j22);
        \draw[line width=0.8pt,color=brown,->] (j11) -- (j24);
        \draw[line width=0.8pt,color=brown,->] (j11) -- (j25);
        
        \draw[line width=0.8pt,color=brown,->] (j12) -- (j21);
        \draw[line width=0.8pt,color=brown,->] (j12) -- (j22);
        \draw[line width=0.8pt,color=brown,->] (j12) -- (j24);
        \draw[line width=0.8pt,color=brown,->] (j12) -- (j25);
        
        \draw[line width=0.8pt,color=brown,->] (j14) -- (j21);
        \draw[line width=0.8pt,color=brown,->] (j14) -- (j22);
        \draw[line width=0.8pt,color=brown,->] (j14) -- (j24);
        \draw[line width=0.8pt,color=brown,->] (j14) -- (j25);
        
        \draw[line width=0.8pt,color=brown,->] (j15) -- (j21);
        \draw[line width=0.8pt,color=brown,->] (j15) -- (j22);
        \draw[line width=0.8pt,color=brown,->] (j15) -- (j24);
        \draw[line width=0.8pt,color=brown,->] (j15) -- (j25);

        \node[circle, fill=green!70,inner sep=5.5pt] (uj) at (7.0,8){};
        \draw[line width=0.8pt,color=brown,->] (j21) -- (uj);
        \draw[line width=0.8pt,color=brown,->] (j22) -- (uj);
        \draw[line width=0.8pt,color=brown,->] (j24) -- (uj);
        \draw[line width=0.8pt,color=brown,->] (j25) -- (uj);

        \node[circle, fill=yellow!70,inner sep=2.25pt] (xk) at  (0.0,13.5){\large$\bm{k}_1\tilde{\bm{x}}$};
        \draw[line width=0.8pt,color=magenta,->] (in) -- (xk);
        
        \node[circle, fill=blue!50,inner sep=2.25pt] (k01) at (1.5,11.5){\large$\sigma$};
        \node[circle, fill=blue!50,inner sep=2.25pt] (k02) at (1.5,12.5){\large$\sigma$};
        \node[circle, fill=blue!0,inner sep=2.25pt] (k03) at (1.5,13.5){$\cdots$};
        \node[circle, fill=blue!50,inner sep=2.25pt] (k04) at (1.5,14.5){\large$\sigma$};
        \node[circle, fill=blue!50,inner sep=2.25pt] (k05) at (1.5,15.5){\large$\sigma$};
        
        \draw[line width=0.8pt,color=magenta,->] (xk) -- (k01);
        \draw[line width=0.8pt,color=magenta,->] (xk) -- (k02);
        \draw[line width=0.8pt,color=magenta,->] (xk) -- (k04);
        \draw[line width=0.8pt,color=magenta,->] (xk) -- (k05);
        
        \node[circle, fill=blue!50,inner sep=2.25pt] (k11) at (3.5,11.5){\large$\sigma$};
        \node[circle, fill=blue!50,inner sep=2.25pt] (k12) at (3.5,12.5){\large$\sigma$};
        \node[circle, fill=blue!0] (k13) at (3.5,13.5){$\cdots$};
        \node[circle, fill=blue!50,inner sep=2.25pt] (k14) at (3.5,14.5){\large$\sigma$};
        \node[circle, fill=blue!50,inner sep=2.25pt] (k15) at (3.5,15.5){\large$\sigma$};
        
        \draw[line width=0.8pt,color=brown,->] (k01) -- (k11);
        \draw[line width=0.8pt,color=brown,->] (k01) -- (k12);
        \draw[line width=0.8pt,color=brown,->] (k01) -- (k14);
        \draw[line width=0.8pt,color=brown,->] (k01) -- (k15);
        
        \draw[line width=0.8pt,color=brown,->] (k02) -- (k11);
        \draw[line width=0.8pt,color=brown,->] (k02) -- (k12);
        \draw[line width=0.8pt,color=brown,->] (k02) -- (k14);
        \draw[line width=0.8pt,color=brown,->] (k02) -- (k15);
        
        \draw[line width=0.8pt,color=brown,->] (k04) -- (k11);
        \draw[line width=0.8pt,color=brown,->] (k04) -- (k12);
        \draw[line width=0.8pt,color=brown,->] (k04) -- (k14);
        \draw[line width=0.8pt,color=brown,->] (k04) -- (k15);
        
        \draw[line width=0.8pt,color=brown,->] (k05) -- (k11);
        \draw[line width=0.8pt,color=brown,->] (k05) -- (k12);
        \draw[line width=0.8pt,color=brown,->] (k05) -- (k14);
        \draw[line width=0.8pt,color=brown,->] (k05) -- (k15);
        
        \node[circle, fill=blue!50,inner sep=2.25pt] (k21) at (5.5,11.5){\large$\sigma$};
        \node[circle, fill=blue!50,inner sep=2.25pt] (k22) at (5.5,12.5){\large$\sigma$};
        \node[circle, fill=blue!0] (k23) at (5.5,13.5){$\cdots$};
        \node[circle, fill=blue!50,inner sep=2.25pt] (k24) at (5.5,14.5){\large$\sigma$};
        \node[circle, fill=blue!50,inner sep=2.25pt] (k25) at (5.5,15.5){\large$\sigma$};
        \node[circle, fill=green!70,inner sep=5.5pt] (uk) at (7.0,13.5){};
        
        \draw[line width=0.8pt,color=brown,->] (k11) -- (k21);
        \draw[line width=0.8pt,color=brown,->] (k11) -- (k22);
        \draw[line width=0.8pt,color=brown,->] (k11) -- (k24);
        \draw[line width=0.8pt,color=brown,->] (k11) -- (k25);
        
        \draw[line width=0.8pt,color=brown,->] (k12) -- (k21);
        \draw[line width=0.8pt,color=brown,->] (k12) -- (k22);
        \draw[line width=0.8pt,color=brown,->] (k12) -- (k24);
        \draw[line width=0.8pt,color=brown,->] (k12) -- (k25);
        
        \draw[line width=0.8pt,color=brown,->] (k14) -- (k21);
        \draw[line width=0.8pt,color=brown,->] (k14) -- (k22);
        \draw[line width=0.8pt,color=brown,->] (k14) -- (k24);
        \draw[line width=0.8pt,color=brown,->] (k14) -- (k25);
        
        \draw[line width=0.8pt,color=brown,->] (k15) -- (k21);
        \draw[line width=0.8pt,color=brown,->] (k15) -- (k22);
        \draw[line width=0.8pt,color=brown,->] (k15) -- (k24);
        \draw[line width=0.8pt,color=brown,->] (k15) -- (k25);
        
        \draw[line width=0.8pt,color=brown,->] (k21) -- (uk);
        \draw[line width=0.8pt,color=brown,->] (k22) -- (uk);
        \draw[line width=0.8pt,color=brown,->] (k24) -- (uk);
        \draw[line width=0.8pt,color=brown,->] (k25) -- (uk);
        
        \node[rectangle, fill=green!70, minimum width=1.25cm, minimum height=12cm, inner sep=0.5pt] (lin) at (9.5,7.25){\rotatebox{-90}{\huge{linear~~~~}}};
        
        \draw[line width=0.8pt,color=brown,->] (uh) -- (8.75,1.0);
        \draw[line width=0.8pt,color=brown,->] (ui) -- (8.75,4.5);
        \draw[line width=0.8pt,color=brown,->] (uj) -- (8.75,8);
        \draw[line width=0.8pt,color=brown,->] (uk) -- (8.75,13.5);

        \node[circle, fill=green!70,inner sep=3.5pt] (y) at (12.5,8.0){\huge{$\hat{\bm{u}}$}};

        \draw[line width=0.8pt,color=brown,->] (10.25,8.0) -- (y);
    \end{tikzpicture}
    \end{center}
    \caption{A schematic diagram of MscaleDNN with $N$ subnetworks.}
    \label{fig: Sub-Fourier PINN Structure}
\end{figure}

From the viewpoint of function approximation, the first layer with activation functions for the DNN framework can be regarded as a series of basis functions and its output is the (nonlinear) combination of those basis functions. Recent works found the phenomenon of spectral bias or frequency preference for DNN, that is, DNN will first capture the low-frequency components of input data~\citep{Xu_2020} \cite{rahaman2018spectral}. After that, some corresponding mechanisms are made using \textcolor{black}{Neural Tangent Kernel (NTK)}~\cite{jacot2018neural, wang2020eigenvector}. Under these mechanisms, many efforts are made to improve the performance of DNN, such as the structures and the activation functions. By introducing some scale factors $\Lambda = (\bm{k}_1, \bm{k}_2, \bm{k}_3\cdots,\bm{k}_{Q-1},\bm{k}_Q)^T$( $\bm{k}_i$ is a vector or matrix), a variety of multi-scale DNN (MscaleDNN) frameworks are proposed which will use the radial scale factors $\Lambda$ to shift the high-frequency component into the low ones, then accelerate the convergence and improve the accuracy of DNN~\cite{wang2020eigenvector, liu2020multi, li2020elliptic}. Figure~\ref{fig: Sub-Fourier PINN Structure} presents the schematic diagram of the MscaleDNN with $N$ subnetworks.

Recent works have shown \textcolor{black}{that the adaptive activation function is able to improve the performance of DNN for solving nonlinear discontinuous problems and image processing~\citep{jagtap2020adaptive, jagtap2020locally, jagtap2022deep2,jagtap2023important}. In addition,} using Fourier feature mapping as an activation function for the first hidden layer of each subnetwork can remarkably improve the capacity of MscaleDNN, it can mitigate the pathology of spectral bias for DNN, and enable networks to learn well the target function~\cite{rahaman2018spectral, wang2020eigenvector, Matthew2020Fourier, Xu_2020, LI2023114963}. It is expressed as follows:
\begin{equation}
\bm{\zeta}(\bm{x}) = 
\left[\begin{array}{c}
\cos(\bm{\kappa} \bm{x})\\
\sin(\bm{\kappa} \bm{x}) 
\end{array}
\right],
\label{fourier}
\end{equation}
where $\bm{\kappa}$ is a user-specified vector or matrix (trainable or untrainable) which is consistent with the number of neural units in the first hidden layer for DNN. By performing the Fourier feature mapping for the input data, the input points in $\mathbb{R}^{d}$ can be mapped to the range $[-1,1]$.
After that, the subsequent layers of the neural network can process the feature information efficiently. For convenience, we denote the PINN model with a MscaleDNN performed by Fourier feature mapping being its solver as Sub-Fourier PINN(called SFPINN).

According to the above description, we denote the proper Fourier feature information of the $n_{th}$ subnetwork by $\bm{\zeta}_{n}(\tilde{\bm{x}})$ with $\tilde{\bm{x}}=(\bm{x},t)$ and obtain its output by performing this information through the remainder block of SFPINN model with general activation functions, such as sigmoid, tanh, and ReLU, etc. Finally, the overall output of the SFPINN model is the linear combination of all subnetwork outputs, denoted by $\boldsymbol{NN}(\tilde{\bm{x}})$. In sum, the detailed procedure is concluded as follows:
\begin{equation}
\begin{aligned}
 \hat{\bm{x}}&= \bm{k}_n\tilde{\bm{x}}, \quad n=1,2, \ldots, N,\\
 \bm{\zeta}_{n}(\tilde{\bm{x}}) & =\left[\cos \left( \boldsymbol{W}^{[n]}_{1} \hat{\bm{x}}\right),\sin \left(\boldsymbol{W}^{[n]}_{1} \hat{\bm{x}}\right)\right]^{\mathrm{T}}, \quad n=1,2, \ldots, N, \\
 \boldsymbol{F}_{n}(\tilde{\bm{x}}) & =\widetilde{\mathcal{F C N}}_{n}\left(\bm{\zeta}_{n}(\tilde{\bm{x}})\right), \quad n=1,2, \ldots, N, \\
 \boldsymbol{NN}(\tilde{\bm{x}}) & =\boldsymbol{W}_O \cdot\left[\boldsymbol{F}_{1}, \boldsymbol{F}_{2}, \cdots, \boldsymbol{F}_{N}\right]+\boldsymbol{b}_O,
\end{aligned}
\end{equation}
where $W^{[n]}_1$ represents the weight matrix of the first hidden layer for the $n_{th}$ subnetwork in the SFPINN model 
and $\widetilde{\mathcal{F C N}_{n}}$ stands for the remaining blocks of the $n_{th}$ subnetwork. $\boldsymbol{W}_O$ and $\boldsymbol{b}_O$ represent the weights and bias of the last linear layer, respectively (see Figure \ref{fig: Sub-Fourier PINN Structure}). Notably, all subnetworks in SFPINN are standalone and their sizes can be adjusted independently.

\section{The process of SFHCPINN algorithm}\label{subsec:algr2SFHCPINN}
Our proposed SFHCPINN is the combination of HCPINN and SFPINN that imposes ``hard'' constraints on the I/BCs and employs a subnetwork structure shown in Figure~\ref{fig: Sub-Fourier PINN Structure}. The solution for ADE \eqref{eq:ADE_IC_BC} is expressed as
\begin{equation}\label{eq:ansatz2SFHCPINN}
    u_{NN}(\bm{x},t)=G(\bm{x},t)+D(\bm{x},t) \boldsymbol{NN}(\bm{x},t ; \bm{\theta})
\end{equation}
where the parameter definitions are identical to Section~\ref{sec:SFHCPINN}. To start with, the smooth extension function $G(\bm{x},t)$ satisfying the I/BCs and smooth distance function $D(\bm{x},t)$ in the distance between interior points to the I/BCs $\partial\Omega \times \{t_0\}$ are constructed (see Section~\ref{sec:SFHCPINN}). Thus before the training procedure of the neural network, our proposed solution has already satisfied the I/BCs. For the SFPINN consisting of $N$ subnetworks (see Section~\ref{subsec:SFPINN}), the input data for each subnetwork will be transformed by the following operation 
\begin{equation*}
\textcolor{black}{\widehat{\bm{x}}= a_n*(\bm{x},t), \quad n=1,2, \cdots, N,}
\end{equation*}
where $a_n>0$ is a scalar factor. Denoting the output of each subnetwork as $\boldsymbol{F}_{n}(n=1,2,\ldots,N)$, then the overall output of the SFPINN model is obtained by 
\begin{equation*}\label{eq:initial solution}
    \boldsymbol{NN}(\bm{x},t) =\frac{1}{N}\sum_{n=1}^{N}\frac{\boldsymbol{F}_{n}}{a_n}.
\end{equation*}
 
In the SFHCPINN algorithm, we let  $u_{NN}^{0}(\bm{x},t)=G(\bm{x},t)+D(\bm{x},t) NN^{0}(\bm{x},t; \bm{\theta})$ be its initial stage. In this stage,  our proposed solution $u_{NN}^{0}(\bm{x},t)$ satisfied the I/BCs in \eqref{eq:ADE_IC_BC} automatically and we can focus on the loss of the interior points. Then in the $k_{th}$ iteration step, a set of randomly sampled collocation points $\mathcal{S}^k$ is provided, and then the $k_{th}$ loss can be obtained by \eqref{softlossHCPINNDir} or \eqref{softlossHCPINNNeu}. The loss function for the Dirichlet boundary is expressed as follows
\begin{equation}\label{loss_SFHCPINNDir}
  Loss^{Dir}_{SFHCPINN}(\mathcal{S}^k;\bm{\theta}) = Loss_{in}(\mathcal{S}_R^k;\bm{\theta})
\end{equation}
with
\begin{equation*}
  Loss_{in}=\frac{1}{|S_{R}^k|}\sum_{i=1}^{|S_{R}^k|}\left| \frac{\partial u_{NN}(\bm{x}^i,t^i)}{\partial t} - \mathbf{div}\big{(}\bm{p}\cdot \nabla u_{NN}(\bm{x}^i,t^i)\big{)}  + \bm{q}\cdot \nabla u_{NN}(\bm{x}^i,t^i) - f(\bm{x}^i,t^i) \right|^2.
\end{equation*}
for $(\bm{x}^i,t^i)\in \mathcal{S}_{R}^k$, \textcolor{black}{where $|S_{R}^k|$ stands for the number of set.} The loss function for the Neumann boundary is formulated as
\begin{equation}\label{loss_SFHCPINNNeu}
    {Loss}_{SFHCPINN}^{Neu}(\mathcal{S}^k;\bm{\theta}) = Loss_{in}(\mathcal{S}_R^k;\bm{\theta}) + \frac{\gamma}{|S_{B}^k|}\sum^{S_{B}^k}_{j=1}\bigg{|}\mathcal{B}u_{NN}(\bm{x}_B^i,t_B^i)-g(\bm{x}_B^i,t_B^i)\bigg{|}^2.
\end{equation}

To sum up, the SFHCPINN method for solving ADE with Dirichlet and/or Neumman boundaries is briefly described in Algorithm~\ref{alg1}.
\begin{algorithm} 
\caption{SFHCPINN algorithm for ADEs with Dirichlet and/or Neumman boundaries} 
\label{alg1}
\begin{algorithmic}[1]
\STATE Construct the extension function $G(\bm{x},t)$ and distance function D$(\bm{x},t)$ according to I/BCs;
\STATE Generate a train set $\mathcal{S}^k\subset \Omega \times T$;
\STATE Calculate the fitting part of loss function $\mathscr{L}(\cdot;\bm{\theta}^{k})$ for training set $\mathcal{S}^k$:
\begin{equation*}
  \mathscr{L}(\mathcal{S}^k;\bm{\theta}^{k})= Loss^{Dir}_{SFHCPINN}(\mathcal{S}_{R}^k;\bm{\theta}^k) ~\textup{or}~\mathscr{L}(\mathcal{S}^k;\bm{\theta}^{k})= Loss^{Neu}_{SFHCPINN}(\mathcal{S}^k;\bm{\theta}^k)
\end{equation*}
\STATE \textcolor{black}{Take a suitable optimization method to update the internal parameters of DNN, such as SGD method expressed as follows}:
\begin{equation*}
\bm{\theta}^{k+1}=\bm{\theta}^{k}-\alpha_k\nabla_{\bm{\theta}^k}\mathscr{L}(\bm{x},t;\bm{\theta}^{k})~~\text{with}~~(\bm{x}, t)\in \mathcal{S}^k
\end{equation*}
where the ``learning rate'' $\alpha_k$ decreases as $k$ increases.
\STATE Repeat steps \textcolor{black}{2-4} until the convergence criterion is satisfied or the loss function tends to be stable.
\end{algorithmic}
\end{algorithm}

\section{Numerical experiments}\label{sec:experiment}
In this section, we test the performance of the SFHCPINN method for advection-diffusion equations \textcolor{black}{with different boundaries in one to three-dimensional spaces}. For comparison, a standard PINN and a sub-Fourier PINN with soft constraints are introduced to validate the effectiveness and feasibility of our SFHCPINN model.

\subsection{Model and training setup}\label{sec:model}
\subsubsection{Model setup}\label{modelsetup}
The details of all the models in the numerical experiments are elaborated in the following and summarized in Table~\ref{tabNN}.
\begin{itemize}
    \item 
	\emph{SFHCPINN}: A solution approach for solving the ADE with different boundary constraints is presented in this article. The approach employs a composite \textcolor{black}{PINN} model composed of a distance function $D(x,t)$, a smooth extension function $G(x,t)$, and a subnetwork deep neural network (DNN). The distance function $D(x,t)$ and smooth extension function $G(x,t)$ are obtained through small-scale DNN training or defined analytically based on boundary conditions. The SFHCPINN model is composed of 20 subnetworks according to the manually defined frequencies $\Lambda=(1,2,3,\ldots, 20)$, and each subnetwork has 5 hidden layers with sizes (10, 25, 20, 20, 15). The activation function of the first hidden layer for each subnetwork is set as Fourier feature mapping $\bm{\zeta}(\hat{\bm{x}})$, and the other activation functions (except for the output layer) are set as $\sin$. The final output of the composite PINN model is a weighted sum of the outputs of all subnetworks. The overall structure of the sub-Fourier PINN model is depicted in Figure~\ref{fig: Sub-Fourier PINN Structure}.

	\item 
	\emph{SFPINN}: The PINN model we consider here uses a subnetwork DNN as the solver, with the activation function of the first hidden layer in each subnetwork set as a Fourier feature mapping and the other activation functions (except for the output layer) set as $\sin$. The I/B constraints in this model are applied in a soft manner, which is the classical approach.
 
	\item 
	\emph{PINN}: The solver for the vanilla PINN model is a normal DNN, where all activation functions except for the output layer are set to $\tanh$. The type of I/B constraints used in this model is a soft manner. 
\end{itemize}

\begin{table}[H]
	\centering
	\caption{ Comparisons for the above models}
	\label{tabNN}
	\begin{tabular}{cccccc}
		\toprule
		  Model &Subnetwork &Numbers of subnetwork & Activation &Constraint& Size of the network\\
            \midrule
		SFHCPINN &DNN   &20   &$\sin$    &hard&(10,25,20,20,10)    \\
		SFPINN   &DNN   &20   &$\sin$    &soft&(10,25,20,20,10)    \\ 
           PINN  &-     &-    &$\tanh$   &soft &(100,150,80,80,50) \\
           \bottomrule
	\end{tabular}
\end{table}
\subsubsection{Training setup}\label{trainingsetup}
We use the following mean square error and \textcolor{black}{ square relative $L^2$ error} to evaluate the accuracy of different models:
\begin{equation*}
MSE = \frac{1}{N'}\sum_{i=1}^{N'}\left(\tilde{u}(x^i,t^i)-u^*(x^i,t^i)\right)^2~~\text{and}~~
REL = \frac{\sum_{i=1}^{N'}\left(\tilde{u}(x^i,t^i)-u^*(x^i,t^i)\right)^2}{\sum_{i=1}^{N'}\left(u^*(x^i,t^i)\right)^2}
\end{equation*}
where $\tilde{u}(x^i,t^i)$ is the approximate DNN solution, $u^*(x^i,t^i)$ is the exact/reference solution, $\{(x^i,t^i)\}_{i=1}^{N'}$ is the set of testing points, and $N'$ is the number of testing points. 

In our numerical experiments, we uniformly sample all training and testing data within $\Omega$ (or $\partial\Omega$) and use the Adam optimizer~\cite{Adam} to train all networks. A step learning rate with an initial learning rate of 0.01 and a decay rate of $2.5\%$ every 100 training epochs is utilized. For visualization purposes, we evaluate our models every 1000 epochs during training and record the final results. The penalty parameter $\gamma$ for the boundary constraint in \eqref{softloss2ADE}, \eqref{softlossHCPINNNeu}, and \eqref{loss_SFHCPINNNeu} is specified as:
\begin{equation}
\gamma=\left\{
\begin{aligned}
\gamma_0, \quad &\textup{if}~~i_{\textup{epoch}}<0.1T_{\max}\\
10\gamma_0,\quad &\textup{if}~~0.1T_{\max}<=i_{\textup{epoch}}<0.2T_{\max}\\
50\gamma_0, \quad&\textup{if}~~ 0.2T_{\max}<=i_{\textup{epoch}}<0.25T_{\max}\\
100\gamma_0, \quad&\textup{if}~~ 0.25T_{\max}<=i_{\textup{epoch}}<0.5T_{\max}\\
200\gamma_0, \quad&\textup{if}~~ 0.5T_{\max}<=i_{\textup{epoch}}<0.75T_{\max}\\
500\gamma_0, \quad&\textup{otherwise}
\end{aligned}
\right.
\end{equation}
where $\gamma_0=20$ in all our tests and $T_{\max}$ represents the total epoch number. We implement our code in Pytorch (version 1.12.1) on a workstation (256 GB RAM, single NVIDIA GeForce GTX 2080Ti 12-GB).

\subsection{Performance of SFHCPINN for solving ADEs}

This section demonstrates the feasibility of SFHCPINN in solving the ADE with Dirichlet and/or Neumann BCs in one-dimensional to three-dimensional Euclidean space. These examples are common in engineering and reality, and four of them involve multi-frequency scenarios to illustrate the ability of SFHCPINN to handle high-frequency problems.

\subsubsection{One-dimensional ADE}

We first consider the one-dimension ADE problem expressed as follows:
\begin{equation}\label{1d_ADE}
\frac{\partial u(x, t)}{\partial t}=\tilde{p} \frac{\partial^{2} u(x, t)}{\partial x^{2}}-\tilde{q}\frac{\partial u(x, t)}{\partial x}+f(x, t) \text {, for } x \in[a, b] \text { and } t \in\left[t_{0}, T\right] \text {,}
\end{equation}
where $f(x,t)$ is the source term, $\tilde{p}$ is coefficient of diffusivity and $\tilde{q}$ is the coefficient of advection rate.
The IC can be formulated as:
\begin{equation}\label{initialC}
u(x,t_{0}) = A_{0}(x).
\end{equation}
The Dirichlet BCs are denoted as:
\begin{equation}\label{DirchC}
\begin{aligned}
&u(a,t)=A_{1}(t),\\
&u(b,t)=A_{2}(t),
\end{aligned}
\end{equation}
and the Neumann BCs could be formulated as follows:
\begin{align}
&\textcolor{black}{\frac{\partial u(a,t)}{\partial x}}=N_{1}(t),\\
&\textcolor{black}{\frac{\partial u(b,t)}{\partial x}}=N_{2}(t).
\end{align}

\begin{example}\label{EX1}
To emphasize the capability of the proposed SFHCPINN framework,
we first consider the \textcolor{black}{one-dimensional ADE \eqref{1d_ADE} with Dirichlet boundary and a high-frequency component. In this case, the \eqref{1d_ADE} with $\tilde{p}=0.02$ and $\tilde{q}=0.01$ is solved by aforementioned models in the spatio-temporal range $\Omega \times T=[0,2]\times[0,5]$}. The source term $f(x,t)$, I/BCs are given by the following known solution:
\begin{equation}\label{eq:0402}
    {u}(x, t)=e^{- \alpha t} [\sin (\pi x) + 0.1\sin (\beta \pi x)]
\end{equation}
with $\alpha = 0.1$ and $\beta=30$.
\end{example}

\textcolor{black}{When employing the proposed SFHCPINN model to solve one-dimensional ADE \eqref{1d_ADE} with Dirichlet boundary, we determine the} distance function $D(x,t) = \frac{x(x-2)t}{5\times 2^2}$ and the smooth function $G(x,t)=\sin (\pi x) + 0.1\sin (\beta \pi x)$ according to the boundary and initial conditions.

To compare with the SFHCPINN model, we use two simple neural networks with only one hidden layer of 20 neurons each to fit the distance function $D(x,t)$ and the extension function $G(x,t)$ before starting the training process. This model is called SFHCPINN$_{NN}$. Both SFHCPINN models are the same except for the distance and extension functions.

\begin{figure}[H]
	\centering
	\subfigure[Exact Solution of Example \ref{EX1}]{
		\label{Exact:EX1}
		\includegraphics[scale=0.35]{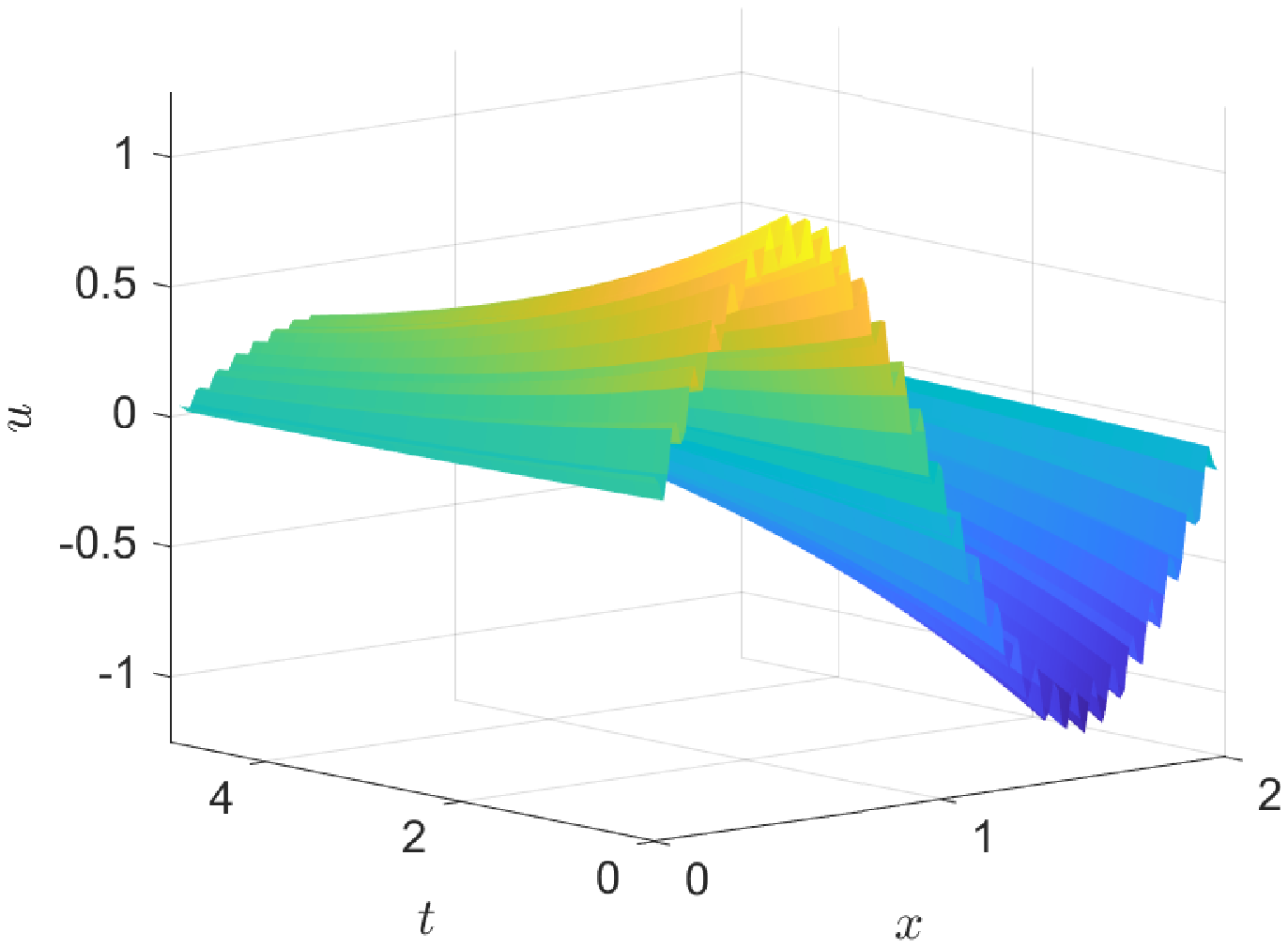}
    }        
	\subfigure[Point-wise error for PINN]{
		\label{ex1:pinnpwe}
		\includegraphics[scale=0.35]{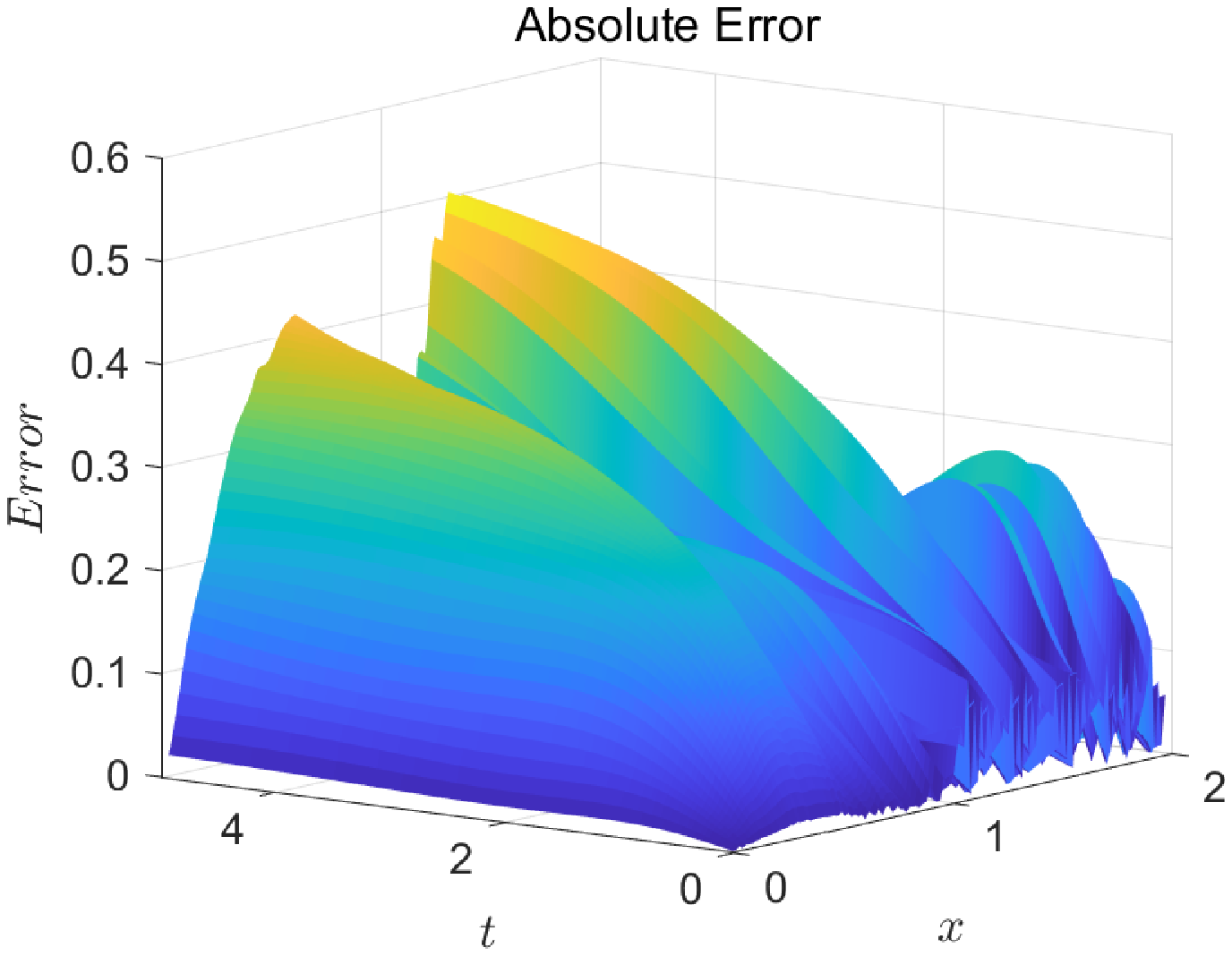}
     }
    \subfigure[Point-wise error for SFPINN]{
		\label{ex1:sfpinnpwe}
		\includegraphics[scale=0.35]{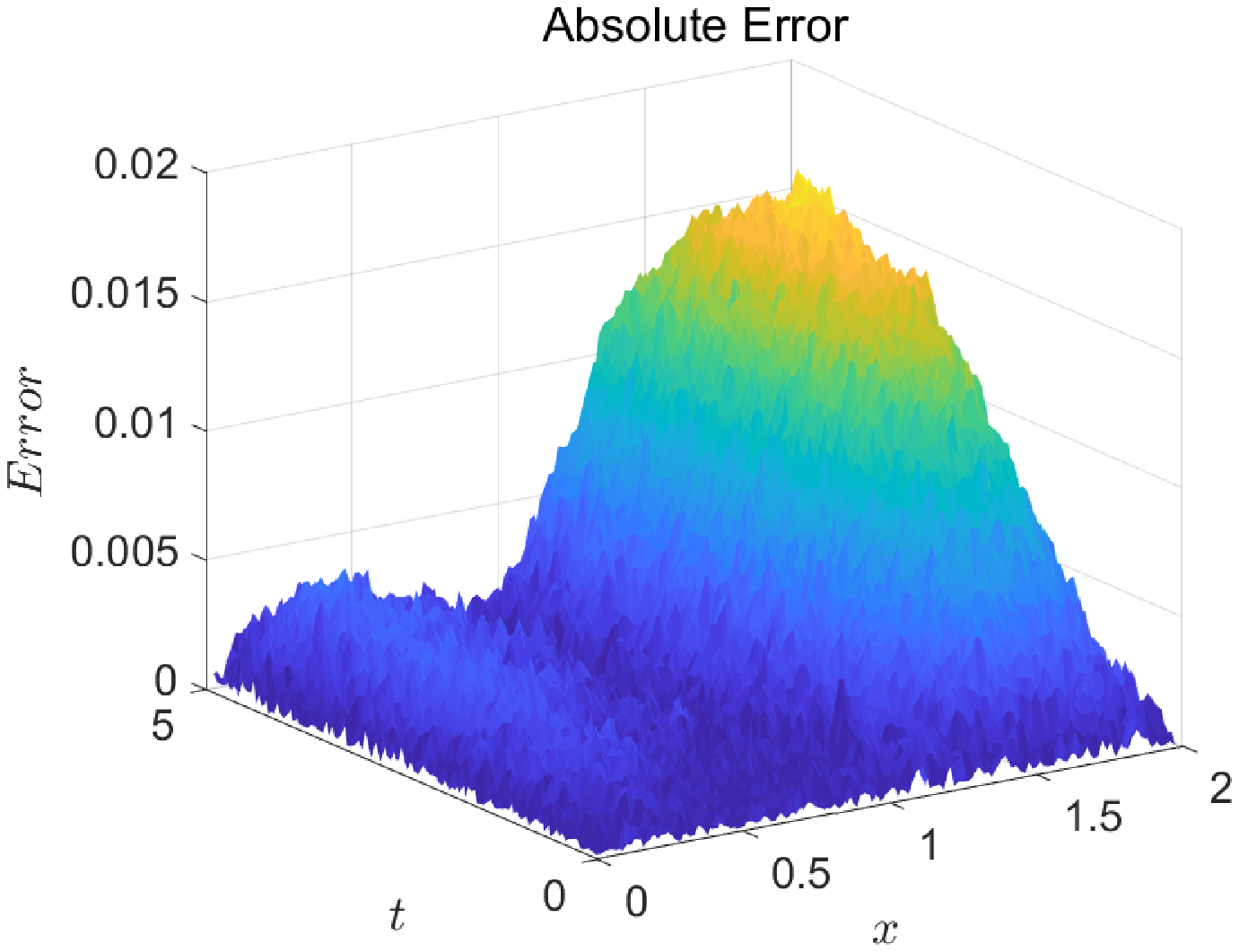}
    }
    \subfigure[Point-wise error for SFHCPINN]{
		\label{ex1:sfhcpinnpwe}
		\includegraphics[scale=0.35]{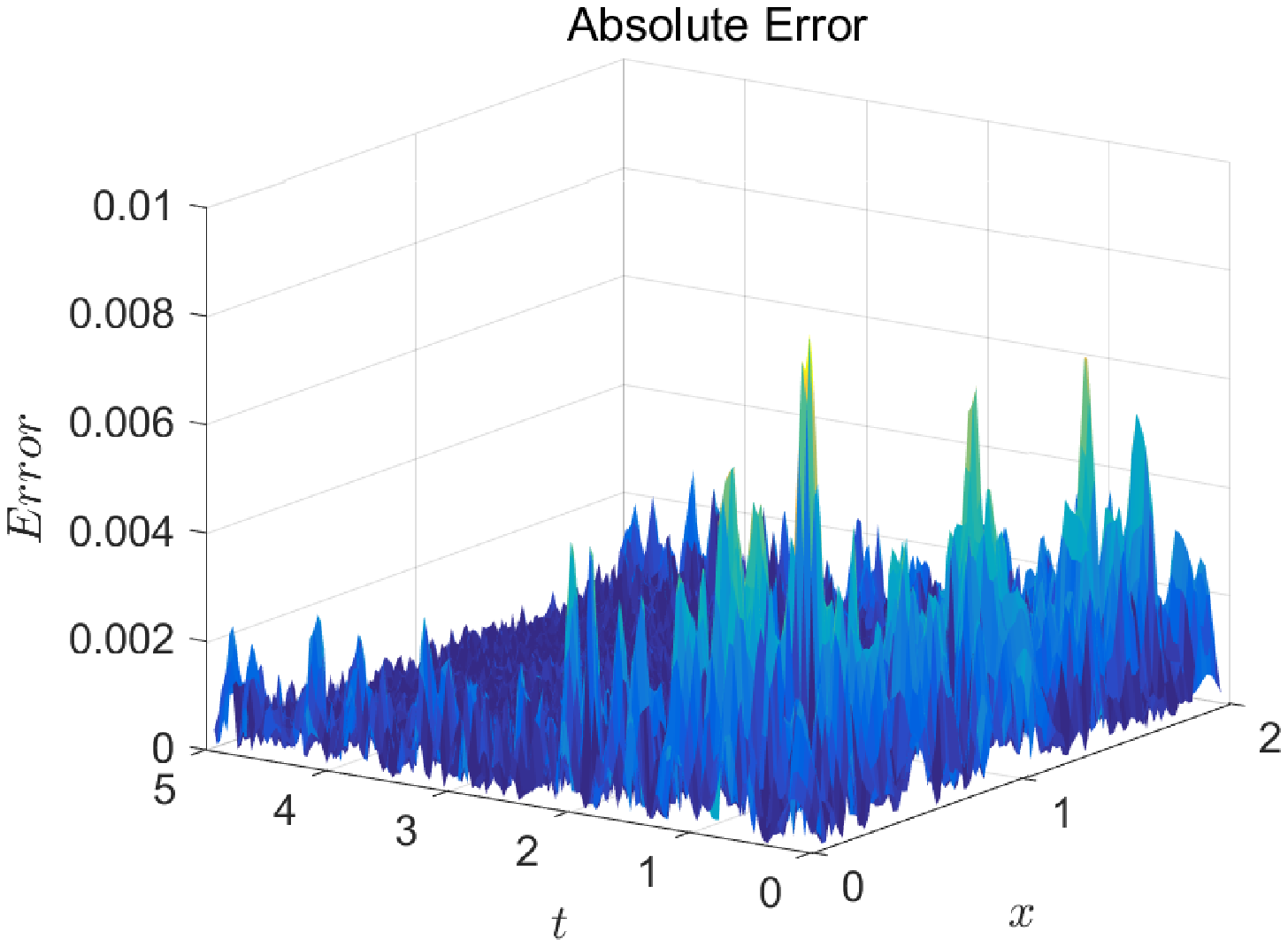}
    }
    \subfigure[Point-wise error for SFHCPINN$_{\textrm{NN}}$]{
		\label{PWE2SFHCPINN_NN}
		\includegraphics[scale=0.35]{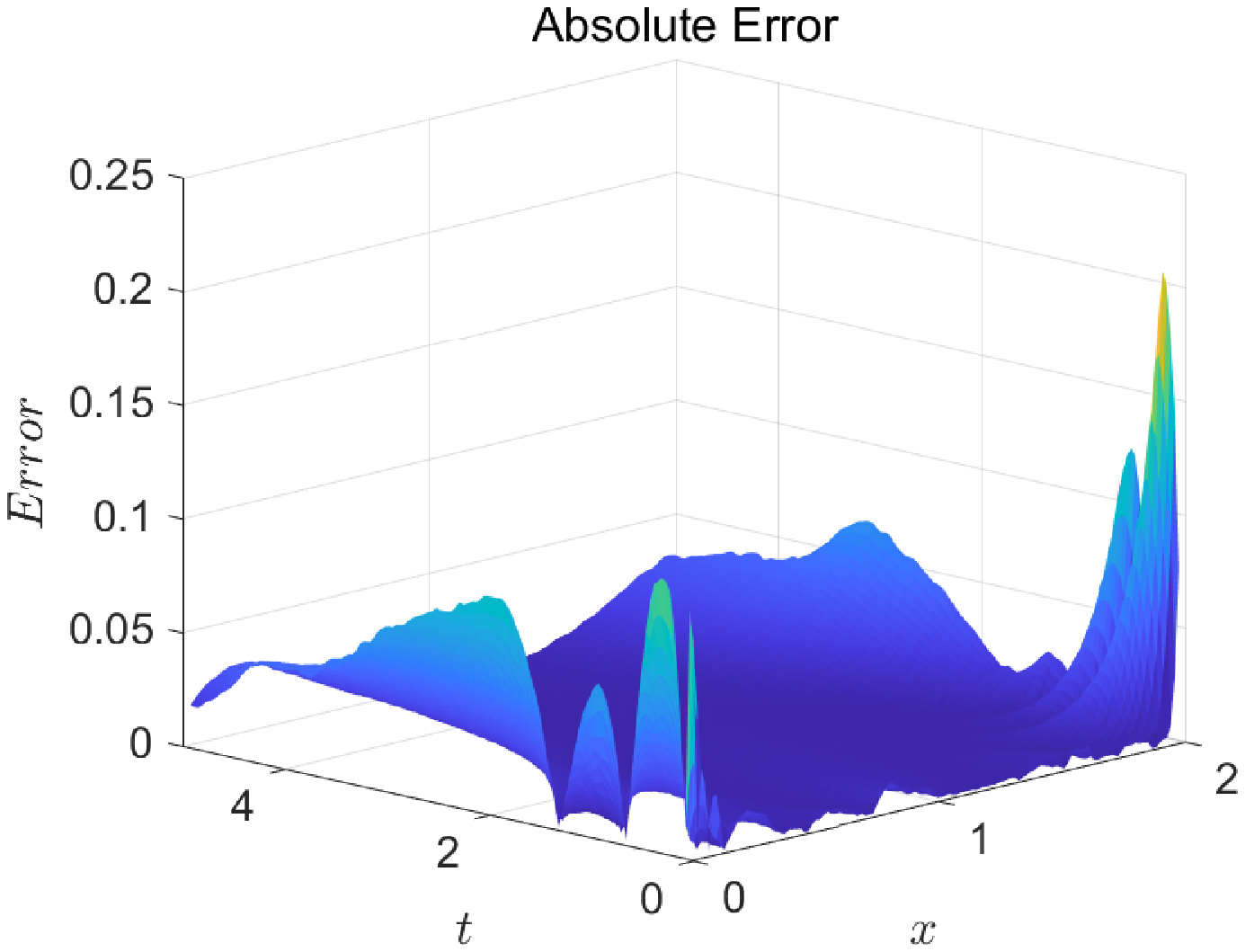}
	}
	\subfigure[REL of PINN, SFPINN and SFHCPINNs]{
		\label{EX1_testrel}
		\includegraphics[scale=0.35]{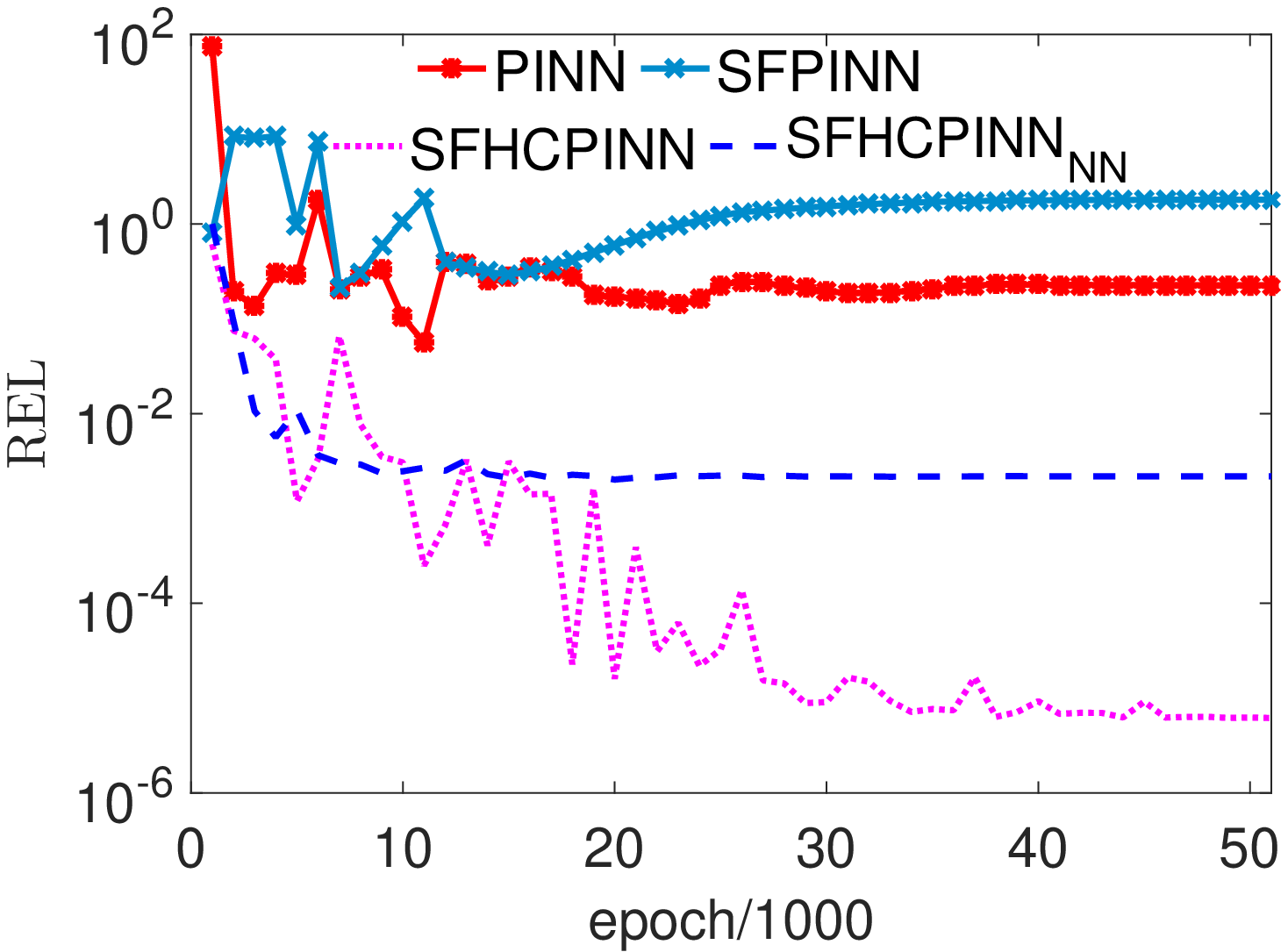}
	}
    \caption{Testing results for Example \ref{EX1}.}
    \label{fig2E1}
\end{figure}

The models have comparable network sizes and parameter numbers, which are listed in Table~\ref{tabNN}. All models are trained for 50000 epochs, and in each epoch, PINN and SFPINN are trained with $N_R=8,000$ collocation points, $N_B=4,000$ boundary points, and $N_I=3,000$ initial points, while SFHCPINN is trained with $N_R=8,000$ collocation points. To assess the accuracy of the neural network approximations, we uniformly sample 10,000 test points from $\Omega \times T$ \textcolor{black}{and make them be the testing set}. We present the results of the four models in Fig. \ref{fig2E1} and Table \ref{tab1}.

\begin{table}[H] 
\caption{MSE and REL of SFHCPINNs, SFPINN, and PINN for Example \ref{EX1}}
\label{tab1}
\begin{center}
\setlength{\tabcolsep}{3pt}
    \begin{tabular}{llll} 
        \toprule 
        & constraint & MSE & REL \\
        \midrule 
          PINN  & soft   & \textcolor{black}{0.0421}  & \textcolor{black}{0.226}\\
          SFPINN  & soft & \textcolor{black}{$4.37\times 10^{-5}$} & \textcolor{black}{$2.34\times 10^{-4}$}\\
          SFHCPINN$_{NN}$& hard &\textcolor{black}{$4.06\times 10^{-4}$} & \textcolor{black}{$2.18\times 10^{-3}$}\\
          SFHCPINN & hard & \textcolor{black}{$1.14\times 10^{-6}$} & \textcolor{black}{$6.16\times 10^{-6}$}\\
        \bottomrule 
    \end{tabular}
\end{center}
\end{table}

First, the heatmaps in Figs.~\ref{ex1:pinnpwe} -- \ref{ex1:sfhcpinnpwe} and the error curves in Fig.~\ref{EX1_testrel} show that SFHCPINN has a higher level of accuracy than PINN and SFPINN, with testing REL decreasing at a faster rate. This suggests that SFHCPINN is effective in addressing the issue of gradient oscillation in DNN parameters, thanks to its use of Fourier expansion and subnetwork framework.

The second observation is that the SFHCPINN exhibits a smaller initial error and a faster convergence rate compared to both the standard PINN and the SFPINN, as shown in Fig.~\ref{EX1_testrel}. This suggests that the hard constraint included in the SFHCPINN allows for better adherence to the boundary conditions, leading to a significant improvement in the performance of SFHCPINN.
 
Lastly, we can observe from the experimental results that SFHCPINN, using numerically determined distance and extension functions, outperforms SFHCPINN$_{NN}$ in terms of both accuracy and training speed. 
This is because the numerically determined distance and extension functions provide a more precise expression of the I/BCs and simple NNs cannot capture functions that vary frequently on the boundaries.
Therefore, we use the numerically determined distance and extension function in all following experiments. In summary, SFHCPINN proves to be superior to PINN and SFPINN in one-dimensional problems with Dirichlet BCs.

\textcolor{black}{\emph{Influence of the choice for activation function:} We study the influence of activation function for PINN and SFHCPINN models. In the test, we set the activation function as sin, tanh, enhance tangent function $tanh(0.5\pi x)$ , sigmoid, elu, gelu, relu, and leaky$\_$relu, all other setups are identical to the above experiments. Based on the results in Table \ref{Table2activation}, we observe that the performance of gelu is slightly superior to that of sin and tanh for our SFHCPINN model, and sin outperforms tanh. Then, the sin activation function is an ideal candidate for our SFHCPINN model.}

\begin{table}[H]
	\centering
	\caption{REL of different activation functions for PINN and SFHCPINN models used to solve Example \ref{EX1}.}
	\label{Table2activation}
	\begin{tabular}{|l|c|c|c|c|c|c|c|c|}
	\hline
         &sin                                    &tanh                                   &enhance tanh              &sigmoid                   &elu                  & gelu      & relu      & leaky relu\\  \hline
	{PINN}     &\textcolor{black}{$0.215$}       &\textcolor{black}{$0.226$}             &\textcolor{black}{$0.224$}&\textcolor{black}{$0.912$}&\textcolor{black}{0.302}&\textcolor{black}{0.0686} &\textcolor{black}{0.362}&\textcolor{black}{0.713} \\  \hline
	SFHCPINN   &\textcolor{black}{$6.16\times 10^{-6}$}&\textcolor{black}{$1.19\times 10^{-4}$}&\textcolor{black}{$0.916$}&\textcolor{black}{$0.015$}&\textcolor{black}{$0.319$}&\textcolor{black}{$9.74\times 10^{-7}$}    &\textcolor{black}{$186.13$}&\textcolor{black}{$20.46$}\\
	\hline
	\end{tabular}
\end{table}

\textcolor{black}{\emph{Influence of the hidden units:}  We study the influence of hidden units including width and depth for our SFHCPINN model. In the test, the hidden units are set as $(10,25,10)$, $(20,25,10)$,$(10,25,20)$, $(20,25,20)$ and $(10,25,20,20)$, respectively, besides the original configuration of hidden units. All other setups are identical to the aforementioned experiments. Based on the results in Table \ref{Table2hiddens}, the performance of SFHCPINN will be improved with the depth and width of hidden units increasing.}

\begin{table}[H]
	\centering
	\caption{REL of different depth and width for hidden units when SFHCPINN is used to solve Example \ref{EX1}.}
	\label{Table2hiddens}
	\begin{tabular}{|l|c|c|c|c|c|}
		\hline
	    $(10,25,10)$                          &$(20,25,10)$                          &$(10,25,20)$                          &$(20,25,20)$                          &$(10,25,20,20)$                        &$(10,25,20,20,10)$  \\ \hline
		\textcolor{black}{$3.24\times10^{-4}$}&\textcolor{black}{$2.03\times10^{-3}$}&\textcolor{black}{$4.05\times10^{-5}$}&\textcolor{black}{$3.91\times10^{-5}$}&\textcolor{black}{$2.93\times 10^{-5}$}&\textcolor{black}{$6.16\times 10^{-6}$}   \\
		\hline
	\end{tabular}
\end{table}

\textcolor{black}{\emph{Influence of the learning rate:} We study the influence of learning rate for our SFHCPINN model. In the test, the SFHCPINN models with fixed learning rates 0.01, 0.005, 0.001, and 0.0005 are introduced to compare with the one with adjustable learning rate. Their other setups are identical to the above experiments. Based on the results in Table \ref{Table2learning}, the adjustable learning rate  is the optimal strategy for the SFHCPINN model.}

\begin{table}[H]
\centering
\caption{REL of the different learning rates for SFHCPINN used to solve Example \ref{EX1}.}
\label{Table2learning}
\begin{tabular}{|l|c|c|c|c|c|}
    \hline
    $0.01$                               &$0.005$                              &$0.001$                              &$0.005$                              &$0.001$                        &Adjustable\\ \hline
    \textcolor{black}{$2.12\times10^{-3}$}&\textcolor{black}{$9.59\times10^{-4}$}&\textcolor{black}{$2.56\times10^{-4}$}&\textcolor{black}{$3.04\times10^{-4}$}&\textcolor{black}{$0.0299$}&\textcolor{black}{$6.16\times 10^{-6}$}\\
    \hline
\end{tabular}
\end{table}

\begin{example}\label{EX2}
\textcolor{black}{We consider the multiscale case of the one-dimensional ADE \eqref{1d_ADE} with Neumann boundaries in spatial domain $\Omega =[0,1]$ and time range $T=[0,1]$},  the source term $f(x,t)$ and the I/BCs are specified by the following prescribed solution: \textcolor{black}{$u(x,t)=e^{-\alpha t}[\sin(\pi x)+0.1\sin(10\pi x)]$ with $\alpha=0.25$}. \textcolor{black}{In this example and next one, $\tilde{p}=0.002$ and $\tilde{q}=0.001$ for equation \eqref{1d_ADE}}. The details of the SFHCPINN with the Neumann BCs are outlined in Algorithm~\ref{alg1}.
At first, \textcolor{black}{the distance function $D(x,t)=x(1-x)t$ and extension function $G(x,t)=\sin(x)$} are well-defined according to the given boundary to specify the distance between the interior and the boundary and the evaluation on the boundary, respectively.
All settings for SFHCPINN, SFPINN, and PINN are the same as in Example \ref{EX1}. 
In each epoch, we randomly collect 8000 points from the interior of the defined domain and 3000 points from the Neumann boundary. The sampling procedures of PINN and SFPINN are identical to those in Example \ref{EX1}. We train all models for 50,000 epochs and use Adam with default parameters as the optimizer. The testing data are uniformly generated from $\Omega \times T$.
\end{example}

\begin{figure}[!ht]
	\centering
	\subfigure[Exact Solution of Example \ref{EX2}]{
		\label{Exact:EX2}
		\includegraphics[scale=0.35]{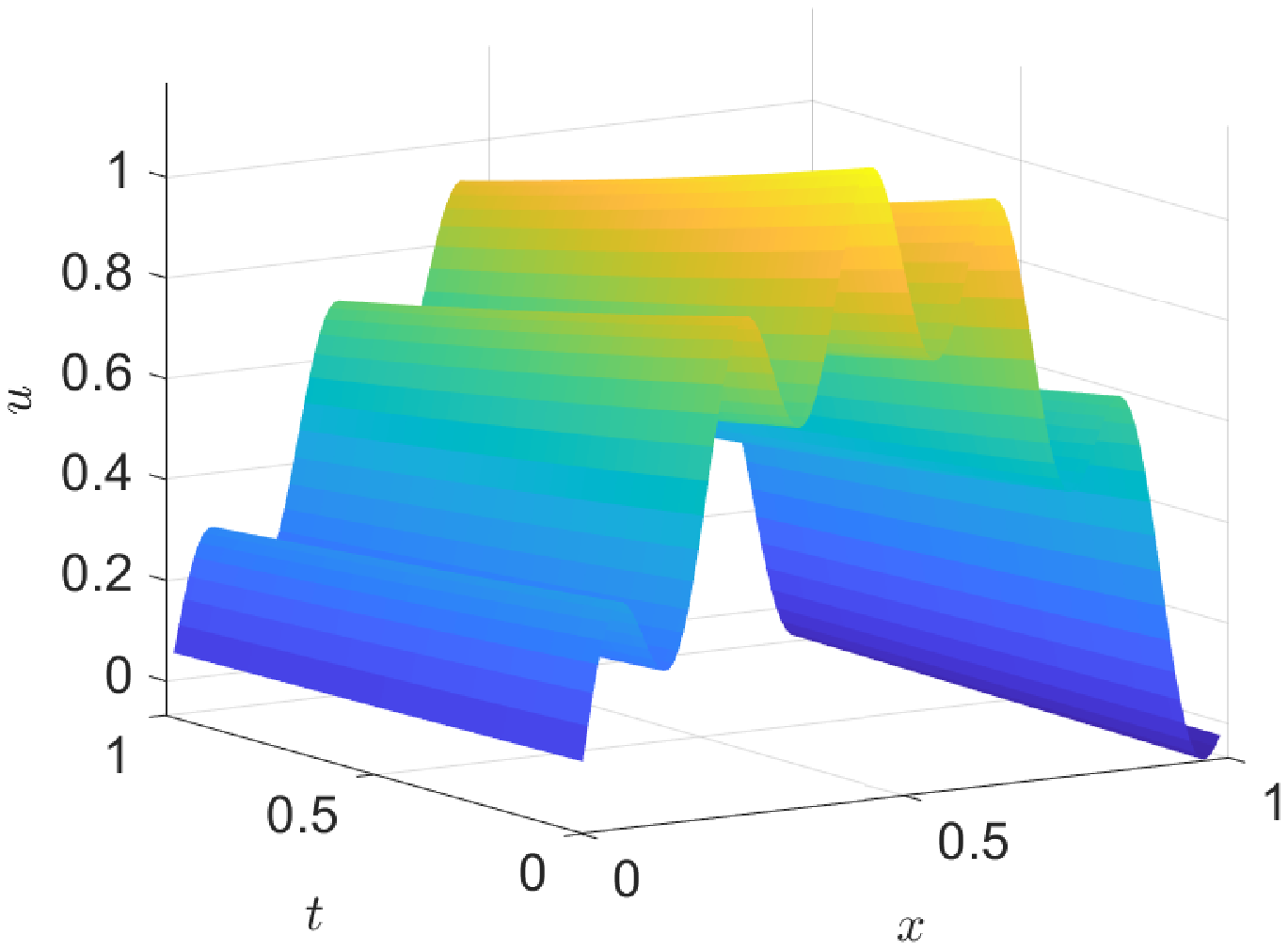}
    }        
	\subfigure[Point-wise error for PINN]{
		\label{ex2:pinnpwe}
		\includegraphics[scale=0.35]{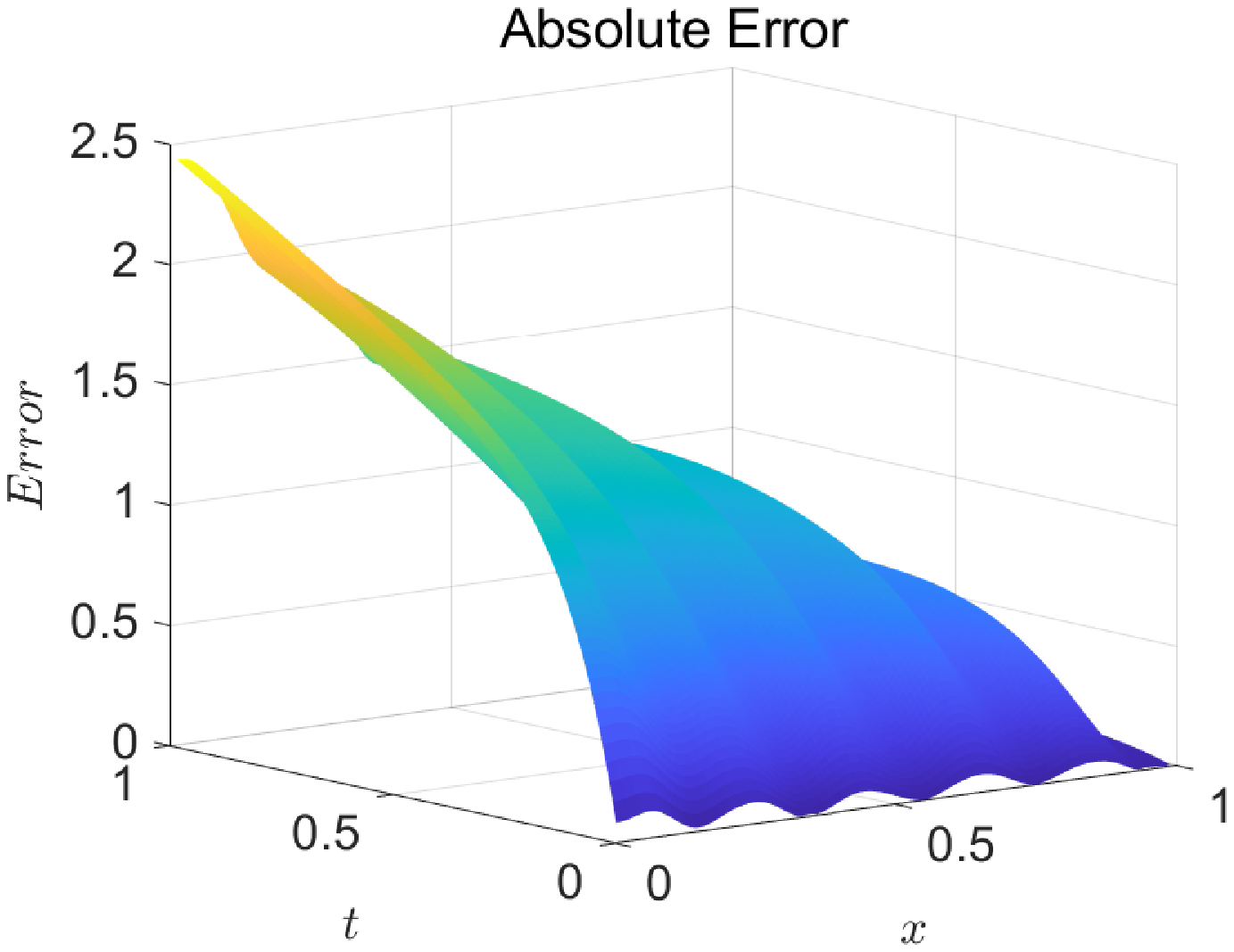}
    }
    \subfigure[Point-wise error for SFPINN]{
	   \label{ex2:sfpinnpwe}
	   \includegraphics[scale=0.35]{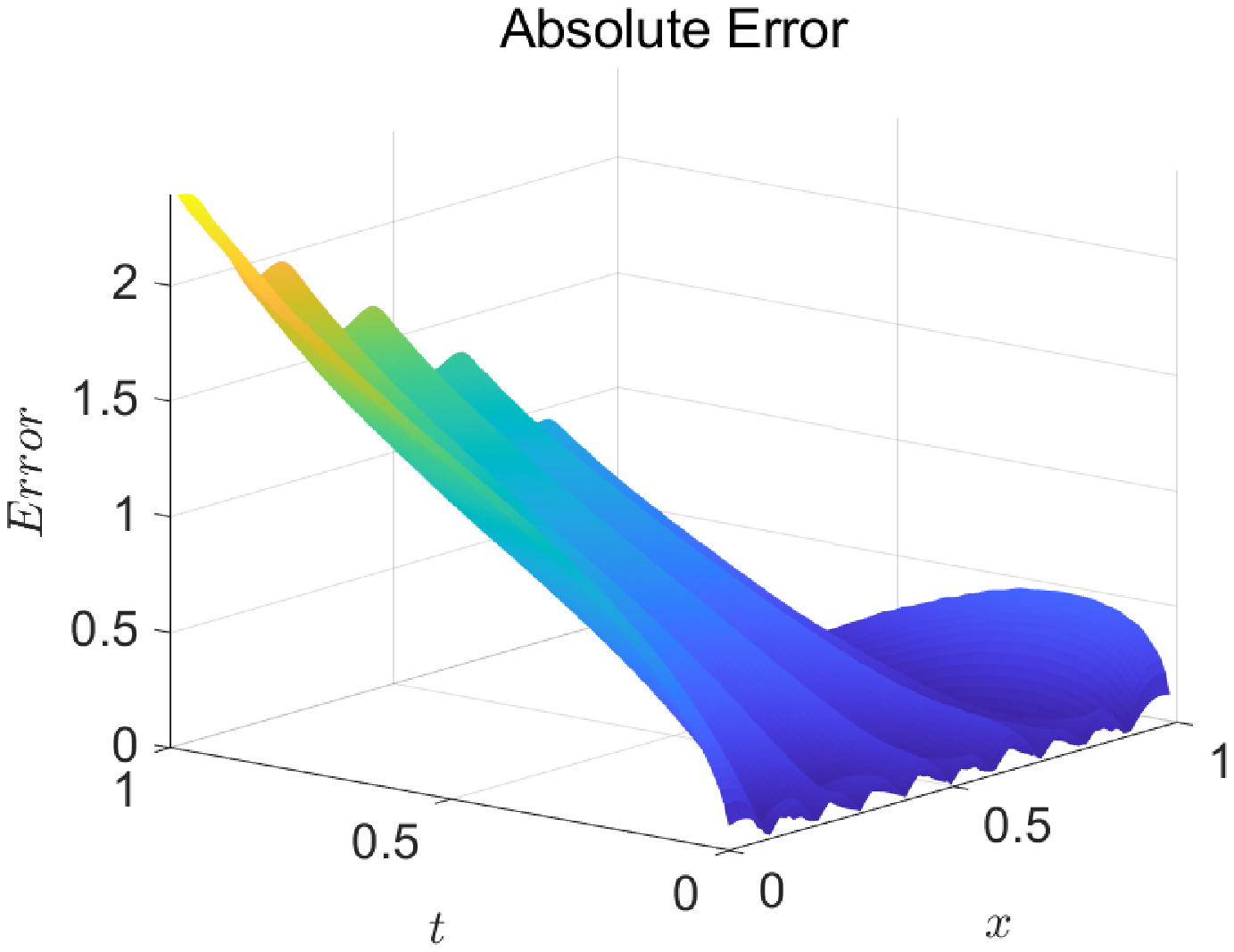}
    }
    \subfigure[Point-wise error for SFHCPINN]{
	   \label{ex2:sfhcpinnpwe}
	   \includegraphics[scale=0.35]{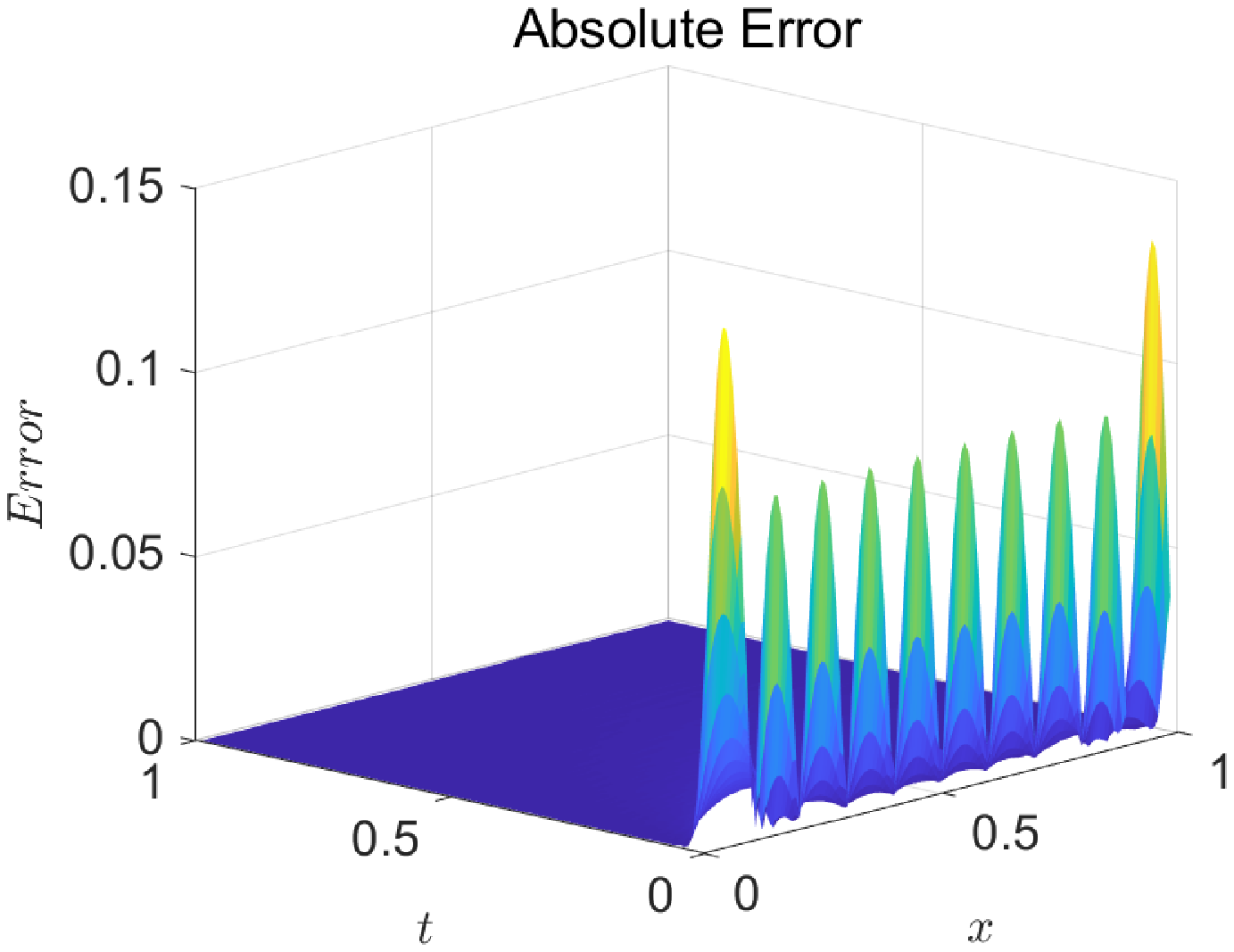}
    }
    \subfigure[ MSE of PINN, SFPINN and SFHCPINN]{
	   \label{EX2_test_mse}
	   \includegraphics[scale=0.35]{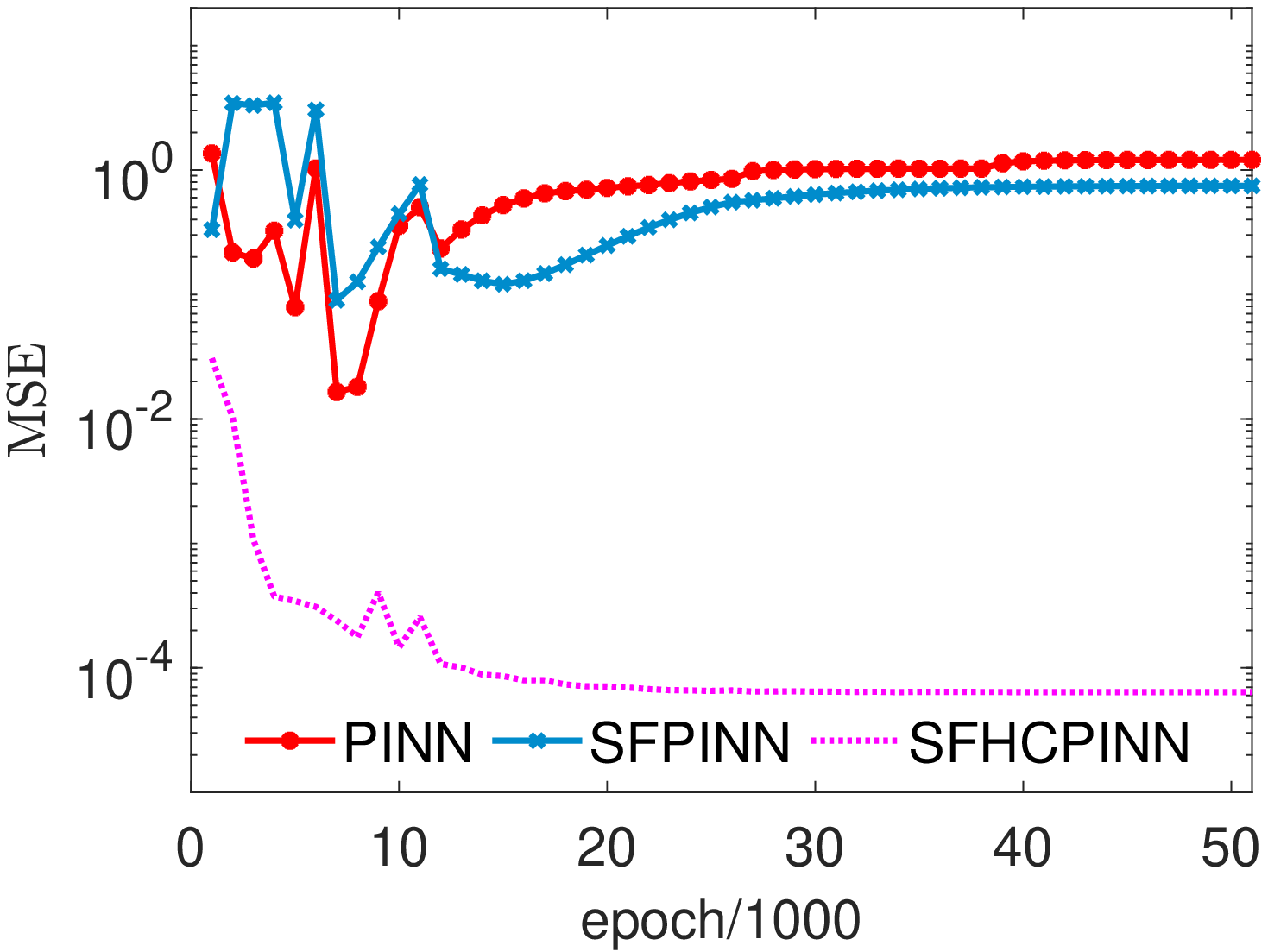}
	}
	\subfigure[ REL of PINN, SFPINN and SFHCPINN]{
		\label{EX2_testrel}
		\includegraphics[scale=0.365]{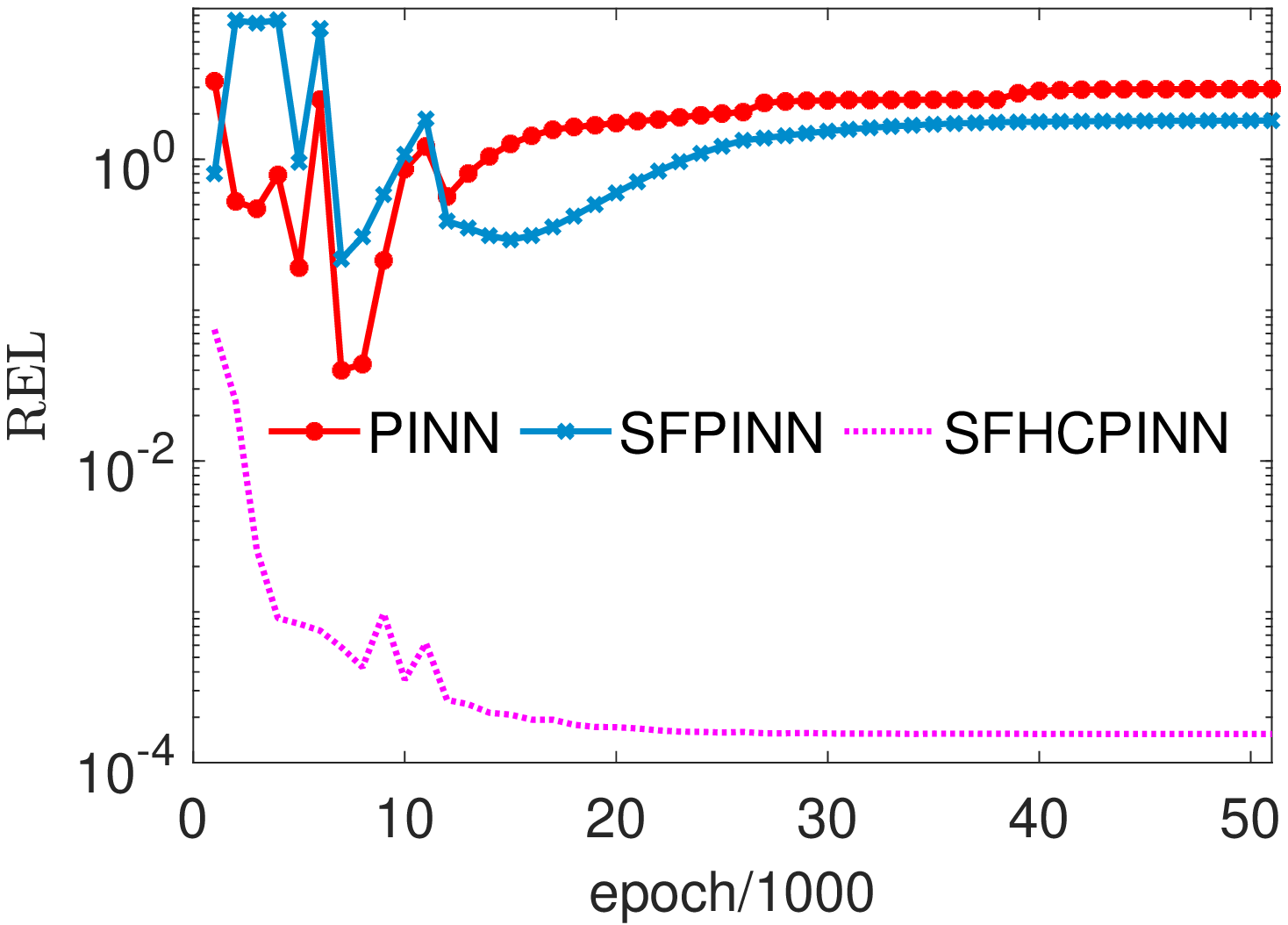}
	}
        \caption{Testing results for Example \ref{EX2}.}
\end{figure}

\begin{table}[!ht] 
	\caption{MSE and REL of SFHCPINN, SFPINN, and PINN for Example \ref{EX2}}
	\label{tab2}
	\begin{center}
		\setlength{\tabcolsep}{3pt}
		\begin{tabular}{llll} 
			\toprule 
			& constraint & MSE & REL \\
			\midrule 
			PINN      & soft & \textcolor{black}{1.204} & \textcolor{black}{2.921}\\
			SFPINN    & soft & \textcolor{black}{0.745} & \textcolor{black}{1.806}\\
			SFHCPINN & hard  &\textcolor{black}{$6.399\times 10^{-5}$} &\textcolor{black}{$1.551\times 10^{-4}$}\\
			\bottomrule 
		\end{tabular}
	\end{center}
\end{table}

The following conclusions can be drawn as follows: First, the diminishing color depth of the thermal maps in Figs.~\ref{ex2:pinnpwe} -- \ref{ex2:sfhcpinnpwe} indicates that \textcolor{black}{the accuracy of the SFHCPINN model improves steadily, but the PINN and SFPINN fail to be converged}. In addition, by comparing the point-wise absolute error, MSE and REL tracks of PINN, SFPINN and SFHCPINN in Figs.~\ref{EX2_test_mse} and \ref{EX2_testrel}, it is possible to conclude that by adopting a subnetwork architecture performed by Fourier activation function and hard-constraint technique, the performance of the DNN has been enhanced with a faster training rate and higher precision under Neumann BCs.
In addition, the accuracy of the SFHCPINN is significantly higher than PINN and SFPINN almost all the time, especially at the initial stage. This is because the ansatz of the hard-constraint PNNN always satisfies the BCs throughout training, preventing the approximations from breaching physical restrictions at the \textcolor{black}{boundaries}. Table~\ref{tab2} further reveals that SFHCPINN with hard constraints and subnet topology outperforms PINN and SFPINN by multiple orders of magnitude when Neumann BCs are present under one-dimensional settings.

\begin{figure}[H]
	\centering
	\subfigure[Exact Solution of Example \ref{Smooth2mixBd_1D}]{
		\label{Exact2MixBD_1D}
		\includegraphics[scale=0.35]{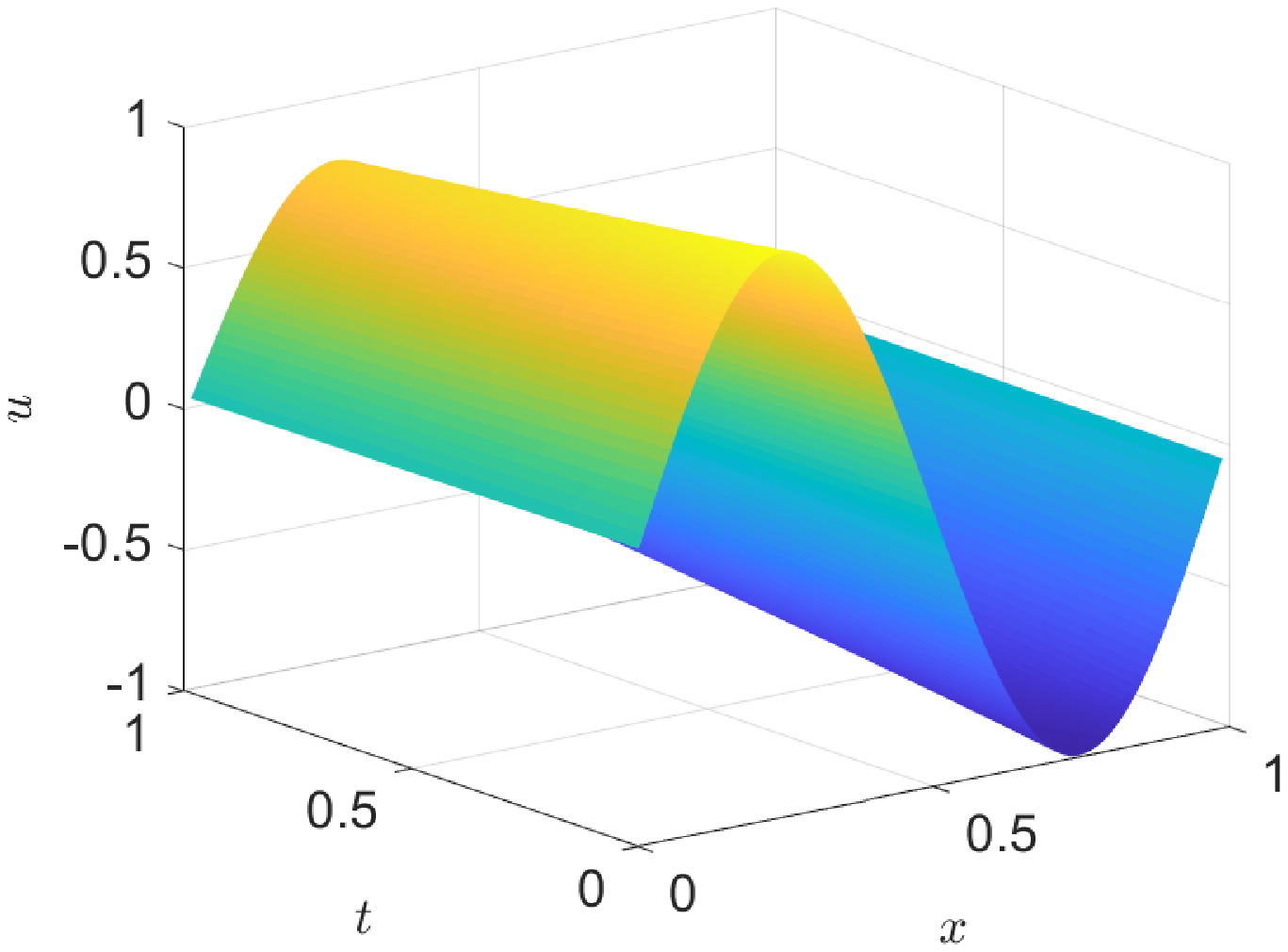}
	}        
	\subfigure[Point-wise error for PINN]{
		\label{PINN_ERR2MixBD_1D}
		\includegraphics[scale=0.35]{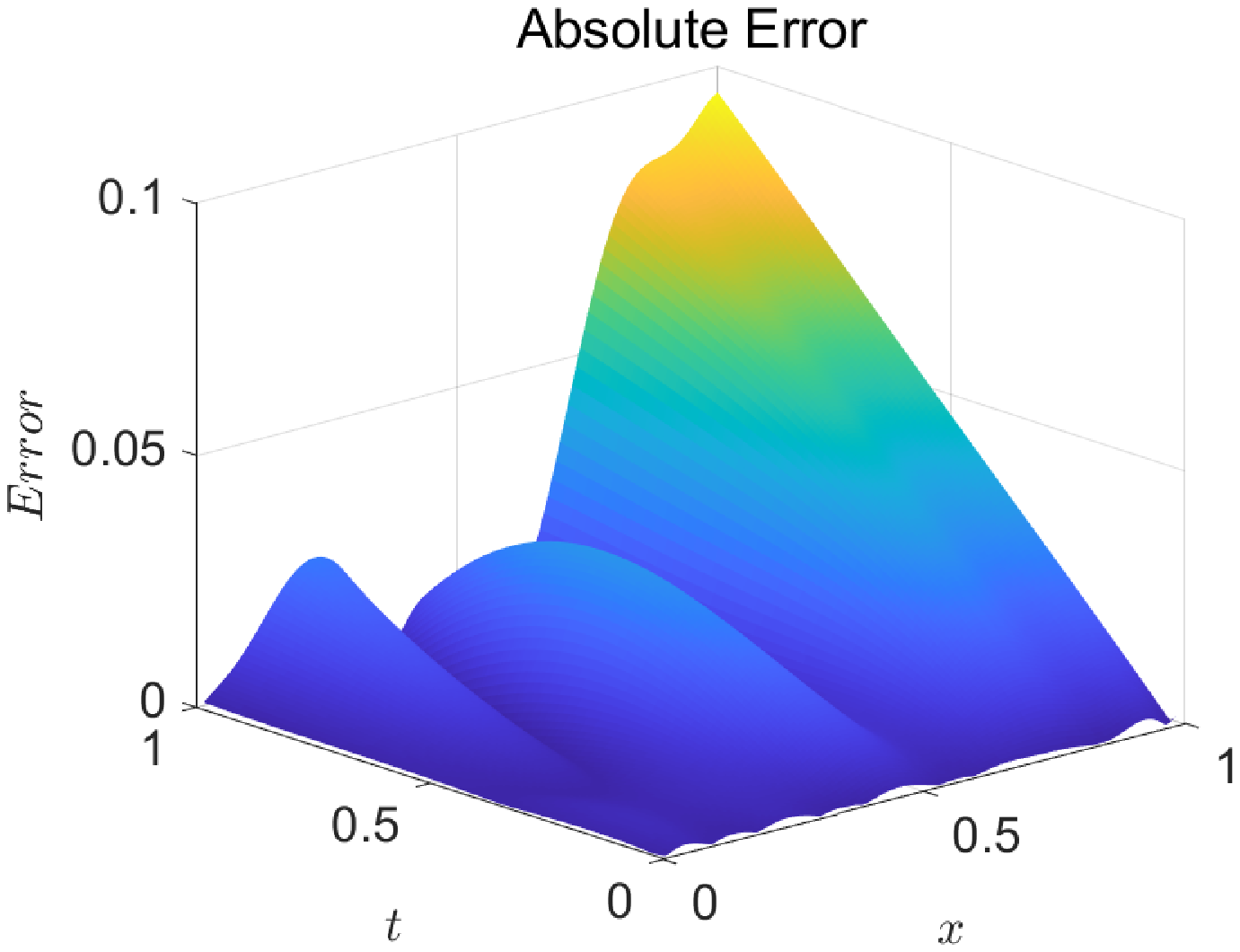}
	}
	\subfigure[Point-wise error for SFPINN]{
		\label{SFPINN_ERR2MixBD_1D}
		\includegraphics[scale=0.35]{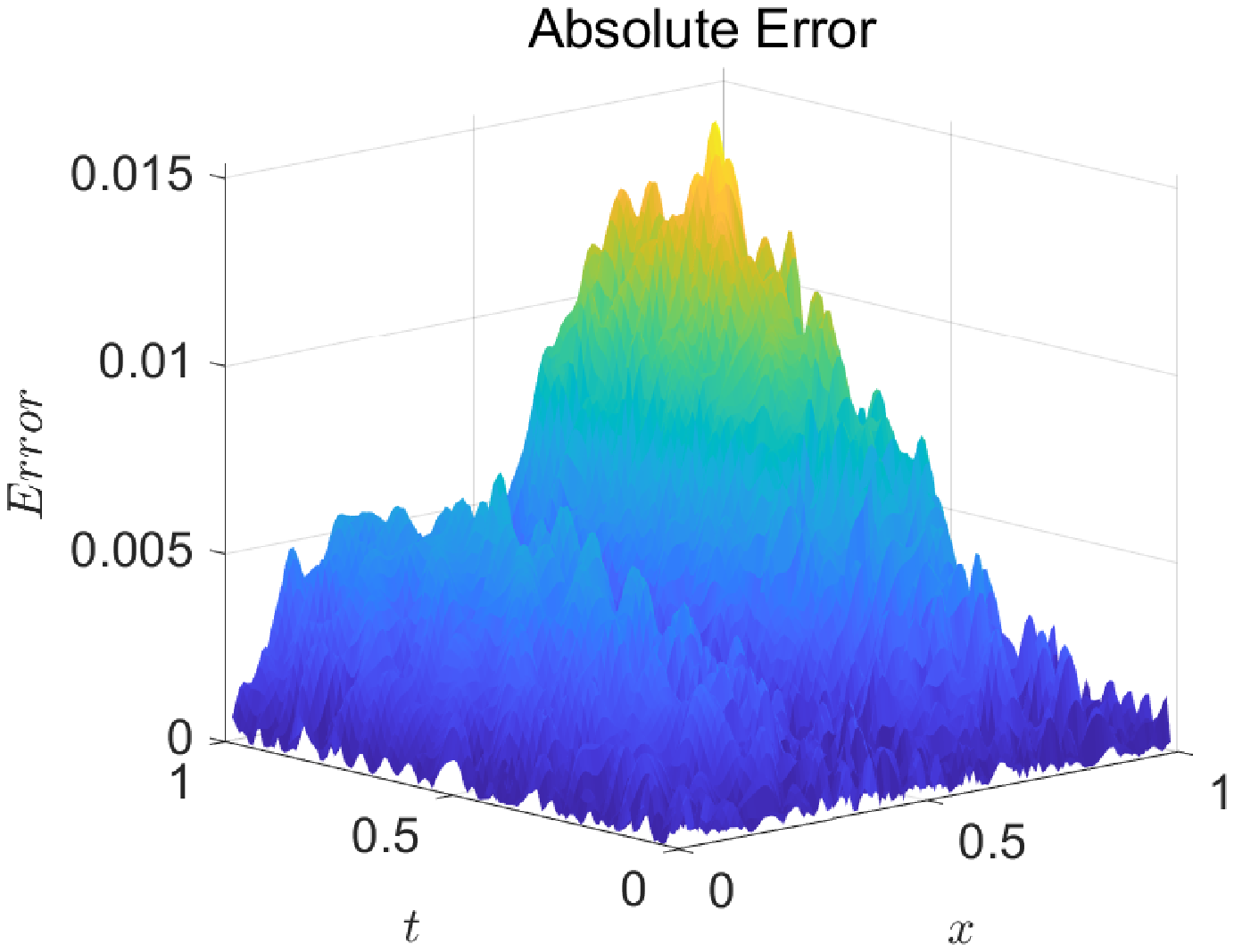}
	}
	\subfigure[Point-wise error for SFHCPINN]{
		\label{SFHCPINN_ERR2MixBD_1D}
		\includegraphics[scale=0.35]{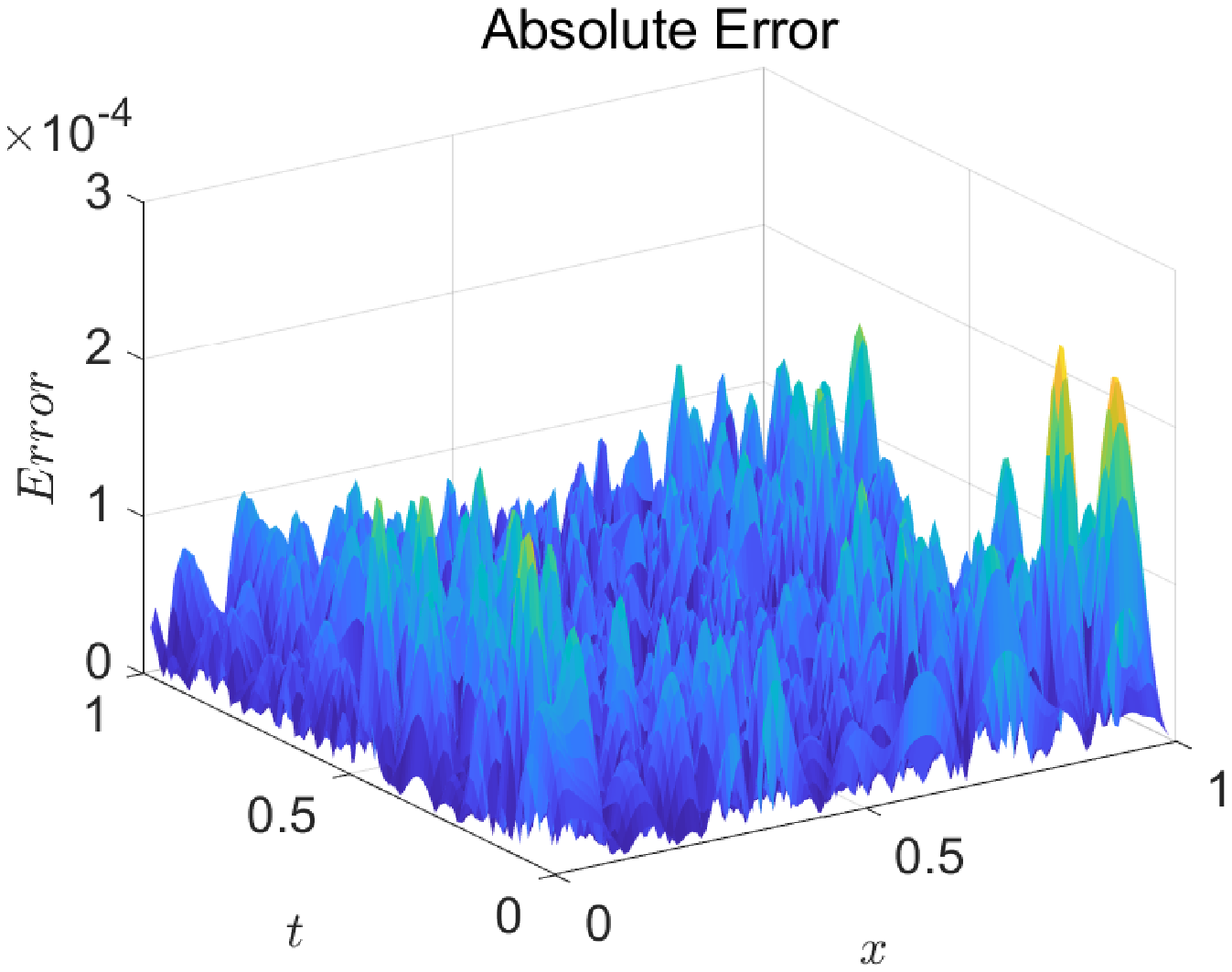}
	}
	\subfigure[ MSE of PINN, SFPINN and SFHCPINN]{
		\label{Test_MSE2MixBD_1D}
		\includegraphics[scale=0.35]{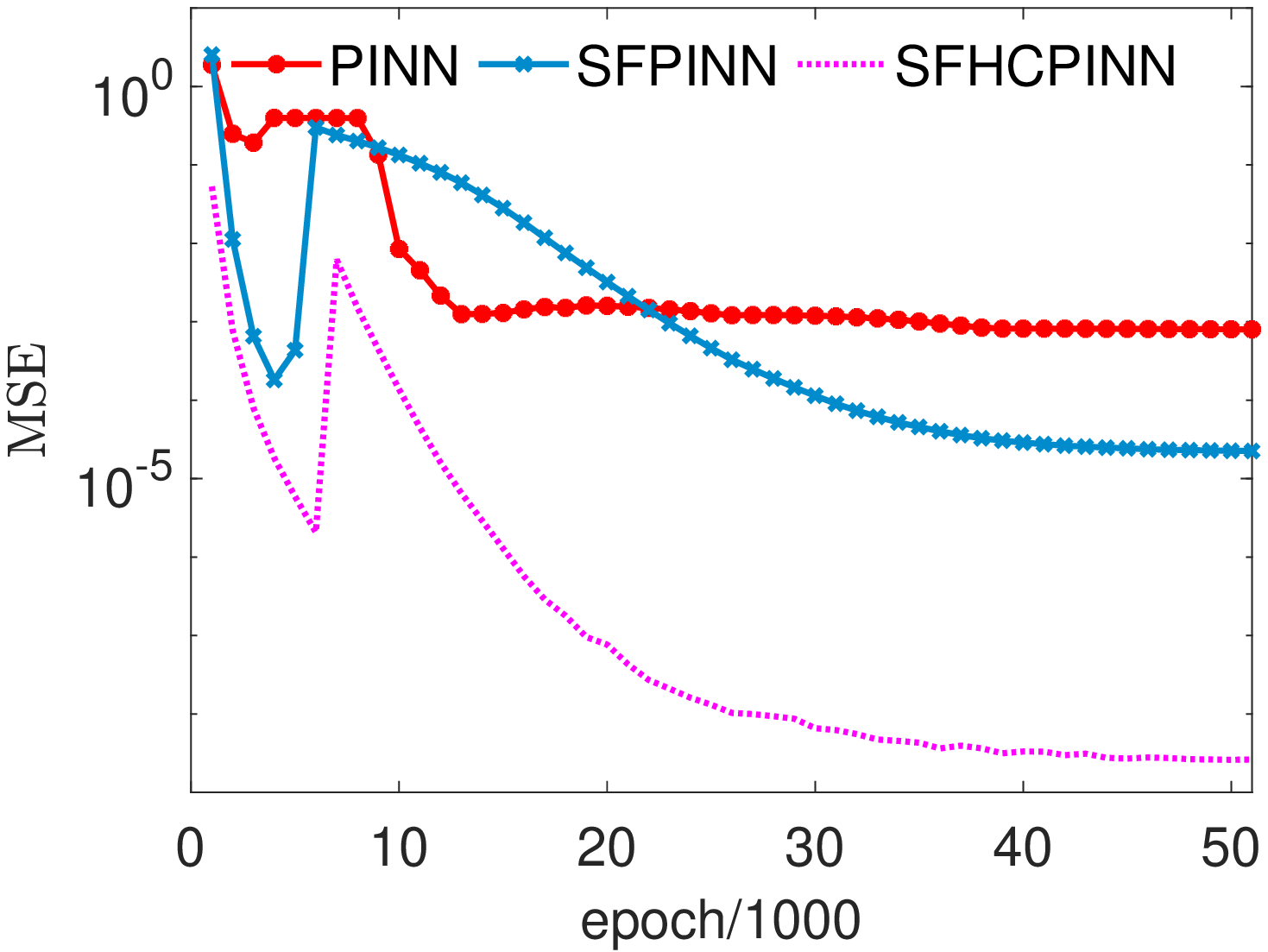}
	}
	\subfigure[ REL of PINN, SFPINN and SFHCPINN]{
		\label{Test_REL2MixBD_1D}
		\includegraphics[scale=0.35]{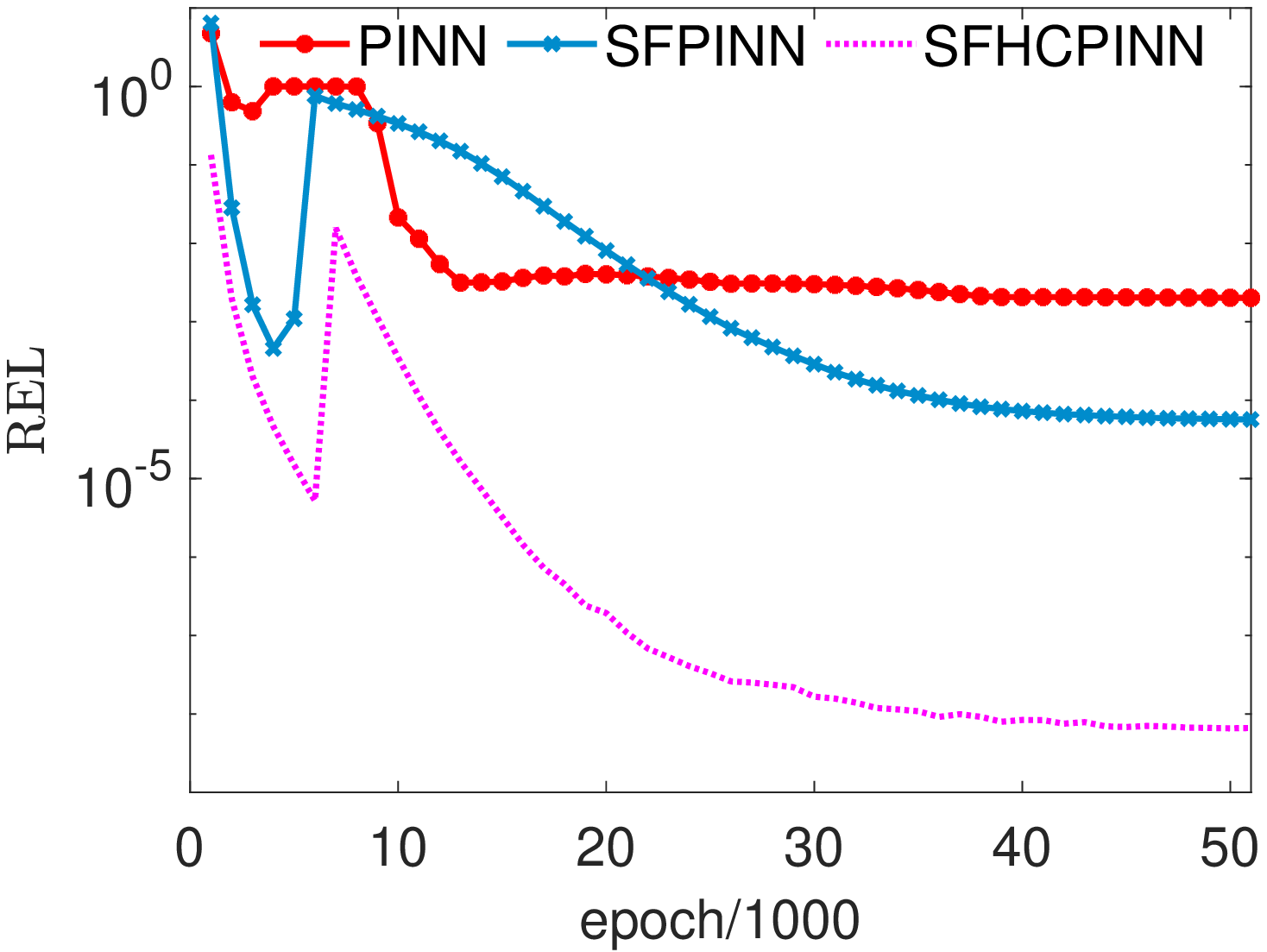}
	}
	\caption{Testing results for Example \ref{Smooth2mixBd_1D}.}
    \label{fig2Smooth2mixBd_1D}
\end{figure}

\begin{example}\label{Smooth2mixBd_1D}
\textcolor{black}{The one-dimensional low-frequency ADE \eqref{1d_ADE} with the left boundary being Dirichlet and the right boundary being Neumann is considered in spatial domain $\Omega =[0,1]$ and time range $T=[0,1]$, the source term $f(x,t)$ and the I/BCs are specified by the known solution: $u(x,t)=e^{-\alpha t}\sin(2\pi x)$ with $\alpha=0.25$. 
At first, \textcolor{black}{the distance function $D(x,t)=x(1-x)t$ and extension function $G(x,t)=\sin(x)$} are well-defined according to the prescribed boundary and initial conditions, respectively.
All settings for SFHCPINN, SFPINN, and PINN are the same as in Example \ref{EX1}. In each epoch, we randomly collect 8000 points from the interior of the defined domain and 3000 points from the Neumann boundary. The sampling procedures of PINN and SFPINN are identical to those in Example \ref{EX1}. We train all models for 50,000 epochs and test our models on a uniform mesh grid generated from $\Omega \times T$. The results are plotted in Fig. \ref{fig2Smooth2mixBd_1D} and listed in Table \ref{Table2Smooth2mixBd_1D}.}
\end{example}

\begin{table}[H] 
	\caption{MSE and REL of SFHCPINN, SFPINN, and PINN for Example \ref{Smooth2mixBd_1D}}
	\label{Table2Smooth2mixBd_1D}
	\begin{center}
		\setlength{\tabcolsep}{3pt}
		\begin{tabular}{llll} 
			\toprule 
			& constraint & MSE & REL \\
			\midrule 
			PINN     & soft &\textcolor{black}{$7.99\times 10^{-4}$} &\textcolor{black}{$2.02\times 10^{-3}$}\\
			SFPINN   & soft &\textcolor{black}{$2.25\times 10^{-5}$} &\textcolor{black}{$5.68\times 10^{-5}$}\\
			SFHCPINN & hard &\textcolor{black}{$2.62\times 10^{-9}$} &\textcolor{black}{$6.63\times 10^{-9}$}\\
			\bottomrule 
		\end{tabular}
	\end{center}
\end{table}

\textcolor{black}{Fig. \ref{fig2Smooth2mixBd_1D} shows that the SFHCPINN model still is well able to capture the solution for the mixed boundaries problem, and its performance outperforms PINN and SFPINN.  Figs. \ref{PINN_ERR2MixBD_1D} -- \ref{SFHCPINN_ERR2MixBD_1D} not only show the point-wise errors of SFHCPINN for major points that are close to zero but also reveal the point-wise error of SFHCPINN is very smaller than that of the PINN and the SFPINN models. Additionally, Figs. \ref{Test_MSE2MixBD_1D} and \ref{Test_REL2MixBD_1D} and Table \ref{Table2Smooth2mixBd_1D} illustrate that the errors of SFHCPINN are superior to that of PINN and SFPINN by more than four orders of magnitude.}

\textcolor{black}{From the above results, we conclude that the SFHCPINN model is remarkable in addressing the ADE \eqref{1d_ADE} with different boundaries in one-dimensional space, it generally outperforms the PINN and SFPINN models.}

\subsubsection{Two-dimensional ADE}
Considering the following two-dimensional ADE with Dirichlet and Neumann boundaries, it is
\begin{equation}\label{TwoDim2ADE}
    \frac{\partial u}{\partial t}+ \hat{p} \frac{\partial u}{\partial x}+ \hat{q} \frac{\partial u}{\partial y}-\left(\frac{\partial^{2} u}{\partial x^{2}}+\frac{\partial^{2} u}{\partial y^{2}}\right)=f(x, y, t) ~~(x, y, z, t) \in \Omega \times(t_0, T]
\end{equation}
\textcolor{black}{where $\hat{p}>0$ and $\hat{q}>0$ are the corresponding diffusion coefficients for different directions.  $\Omega$ is the interested domain, $t_0\geqslant0$ is the initial time and $T>0$ is the end time, and $u$ is the function to be solved. Please refer \citet{NAZIR20164586} to see more details.}

\begin{example}\label{Smooth_Porous}
\textcolor{black}{Let us approximate the solution of \eqref{TwoDim2ADE} with Dirichlet boundary in time interval $[0,5]$ and a porous domain $\Omega$ when $\hat{p}=4$ and $\hat{q}=4$, the $\Omega$ is inserted in a square domain $[0,4]\times[0,4]$.  To address this problem, we perform our model in this regular domain and then obtain their approximations on an interested domain.} An precise solution is given by $u(x,y,t)=e^{-0.25t}xy(4-x)(4-y)$ and it specifies the initial and Dirichlet conditions. 
\textcolor{black}{For SFHCPINN model, we set the distance function as $\displaystyle D(x,y,t)=\frac{t}{5}\sin(0.25x)\sin(0.25y)$ which will be vanish on the boundary of $\Omega$. The extension function can be chosen as $G(x,y,t)=(4x-x^2)(4y-y^2)$ according to the boundary and initial conditions.} The above three models for solving \eqref{TwoDim2ADE} with Dirichlet boundary mentioned in Section~\ref{modelsetup} are identical to the ones described in Examples~\ref{EX1} and~\ref{EX2} in that they share the same optimizer, learning rate, and decay rate. 
In the training stage, we randomly sample 8,000 interior points as training points for SFHCPINN, and we sample an additional 3,000 initial points and 3,000 boundary points for PINN and SFPINN in every epoch.
\textcolor{black}{In addition, 17706 random collocations sampled from porous domian $\Omega$ at $t=2.5$ are used as the testing set. The results are plotted in Fig. \ref{fig2smooth_porous_dirichlet2d} and listed in Table \ref{table2smooth_porous_dirichlet2d}.}

\begin{figure}[H]
   \centering
    \subfigure[Exact Solution of Example \ref{Smooth_Porous}]{
		\label{Exact:EX3}
		\includegraphics[scale=0.375]{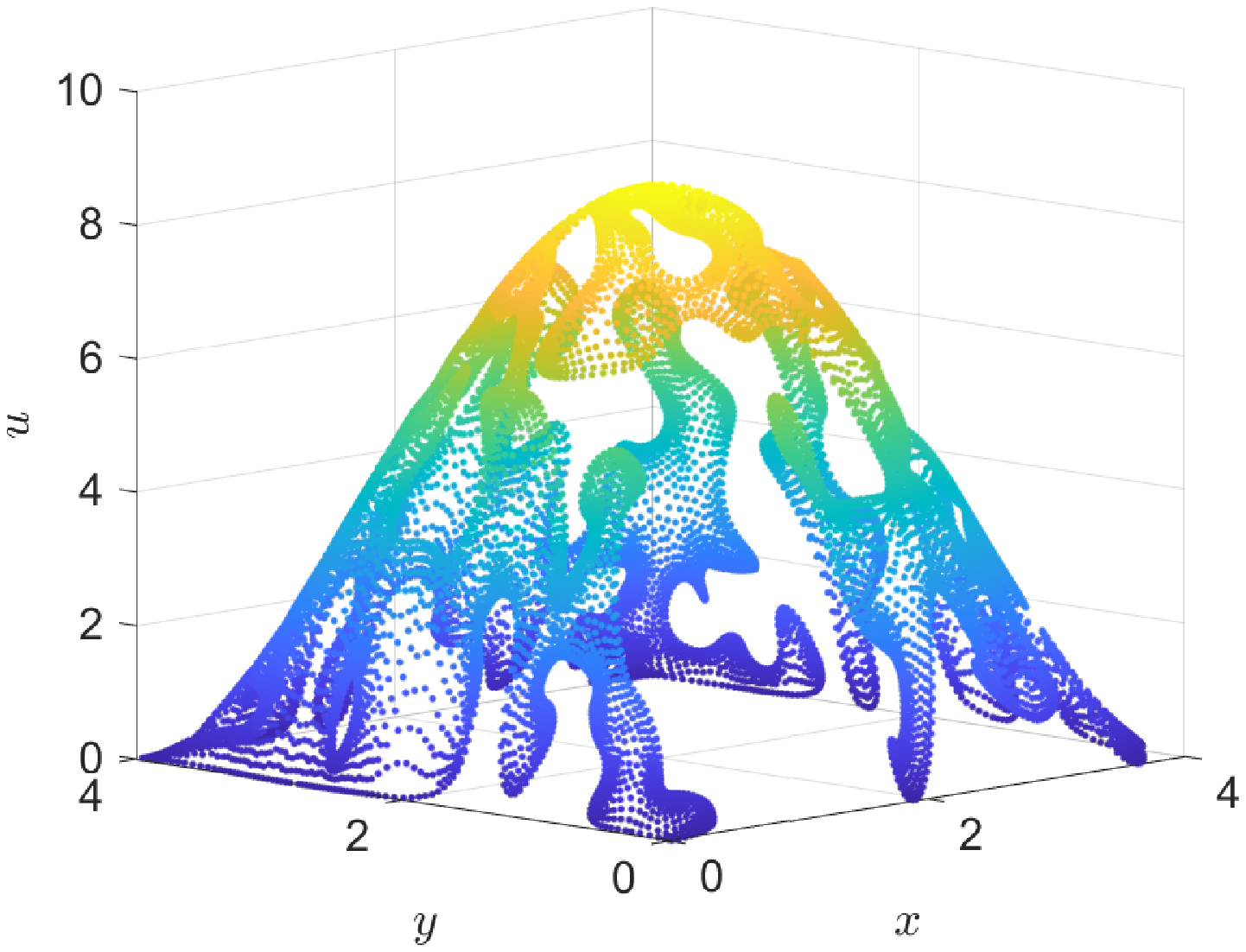}
     }        
    \subfigure[Point-wise error for PINN]{
		\label{ex3:pinnpwe}
		\includegraphics[scale=0.35]{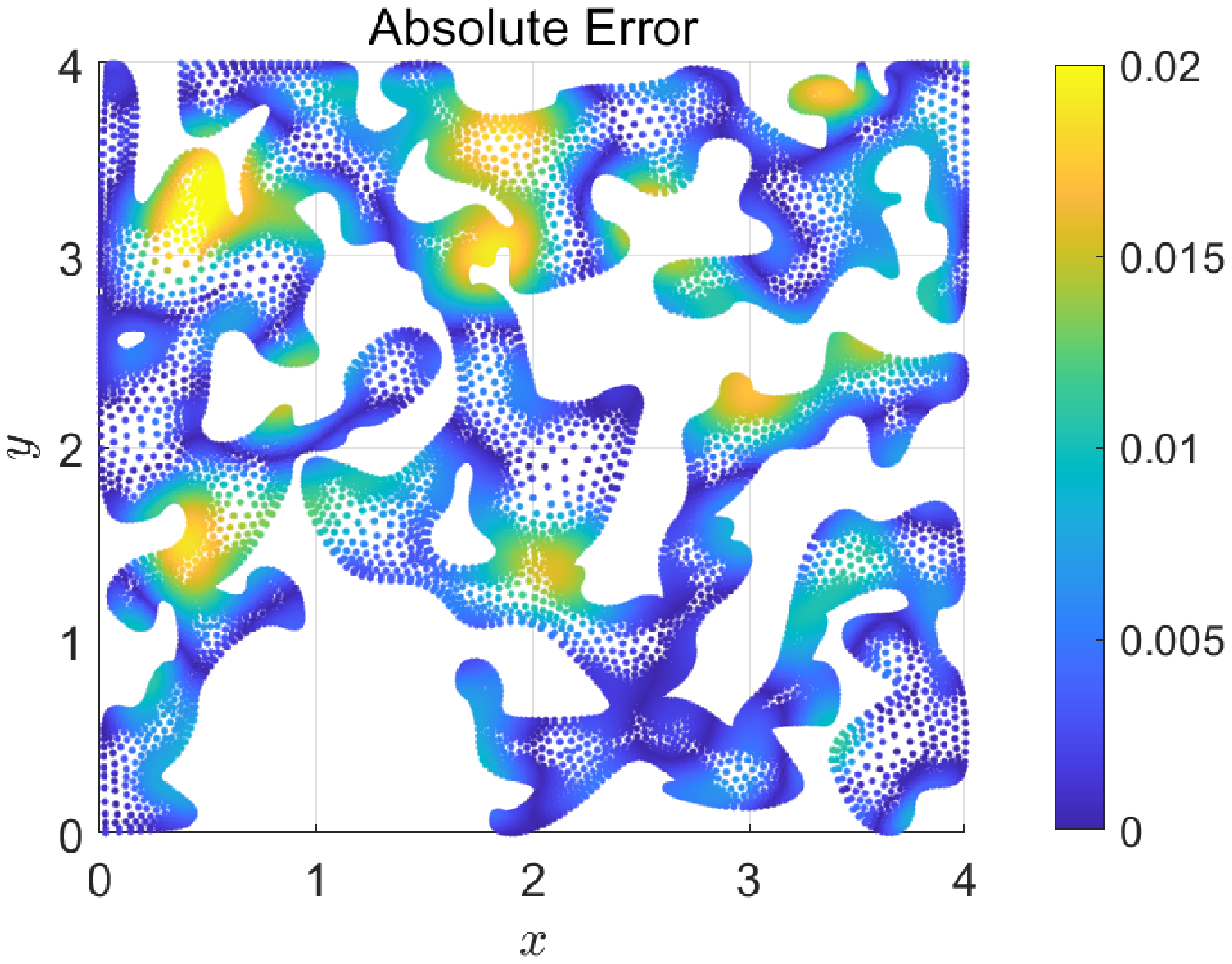}
    }
    \subfigure[Point-wise error for SFPINN]{
    	\label{ex3:sfpinnpwe}
    	\includegraphics[scale=0.35]{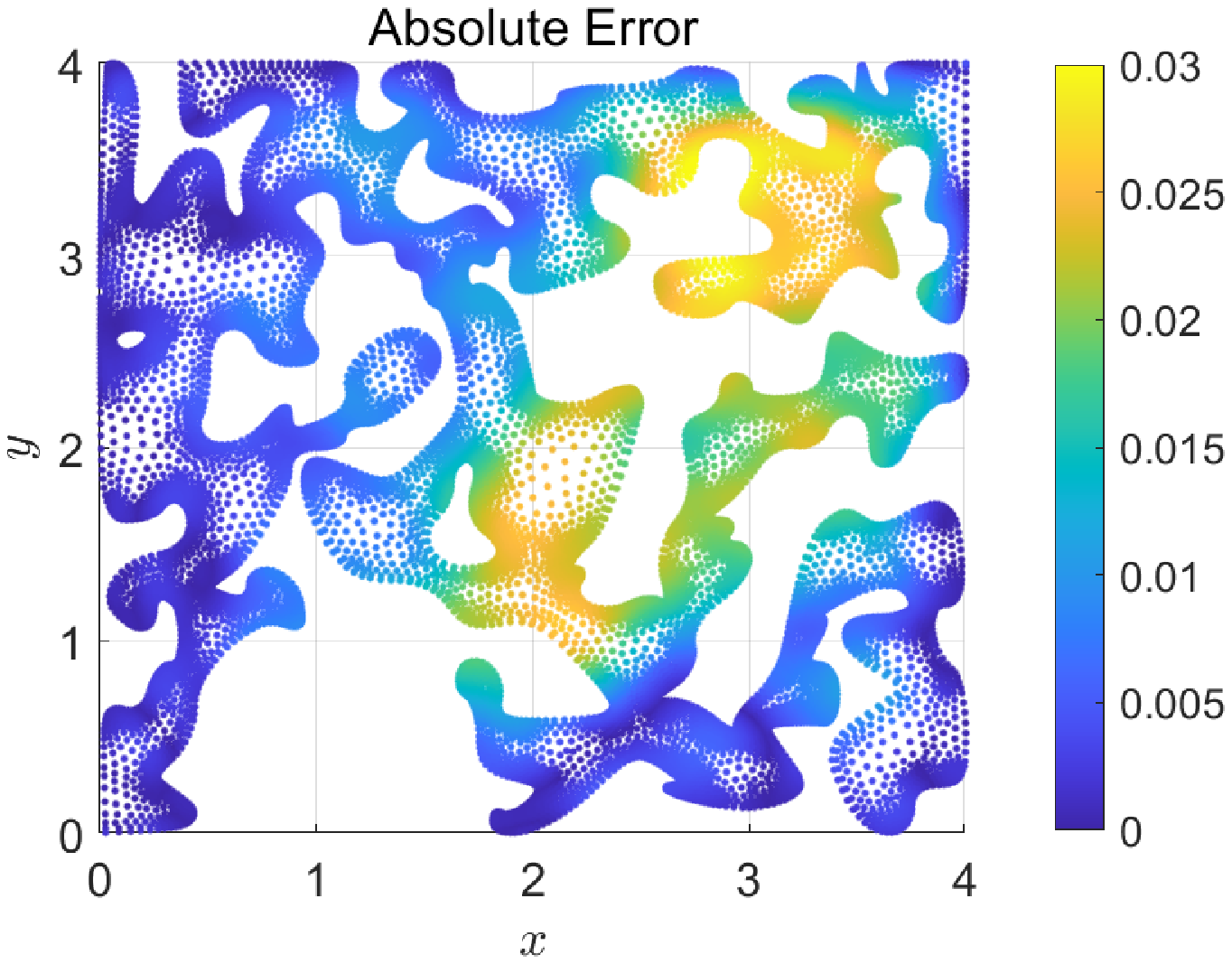}
    }
    \subfigure[Point-wise error for SFHCPINN]{
    	\label{ex3:sfhcpinnpwe}
    	\includegraphics[scale=0.35]{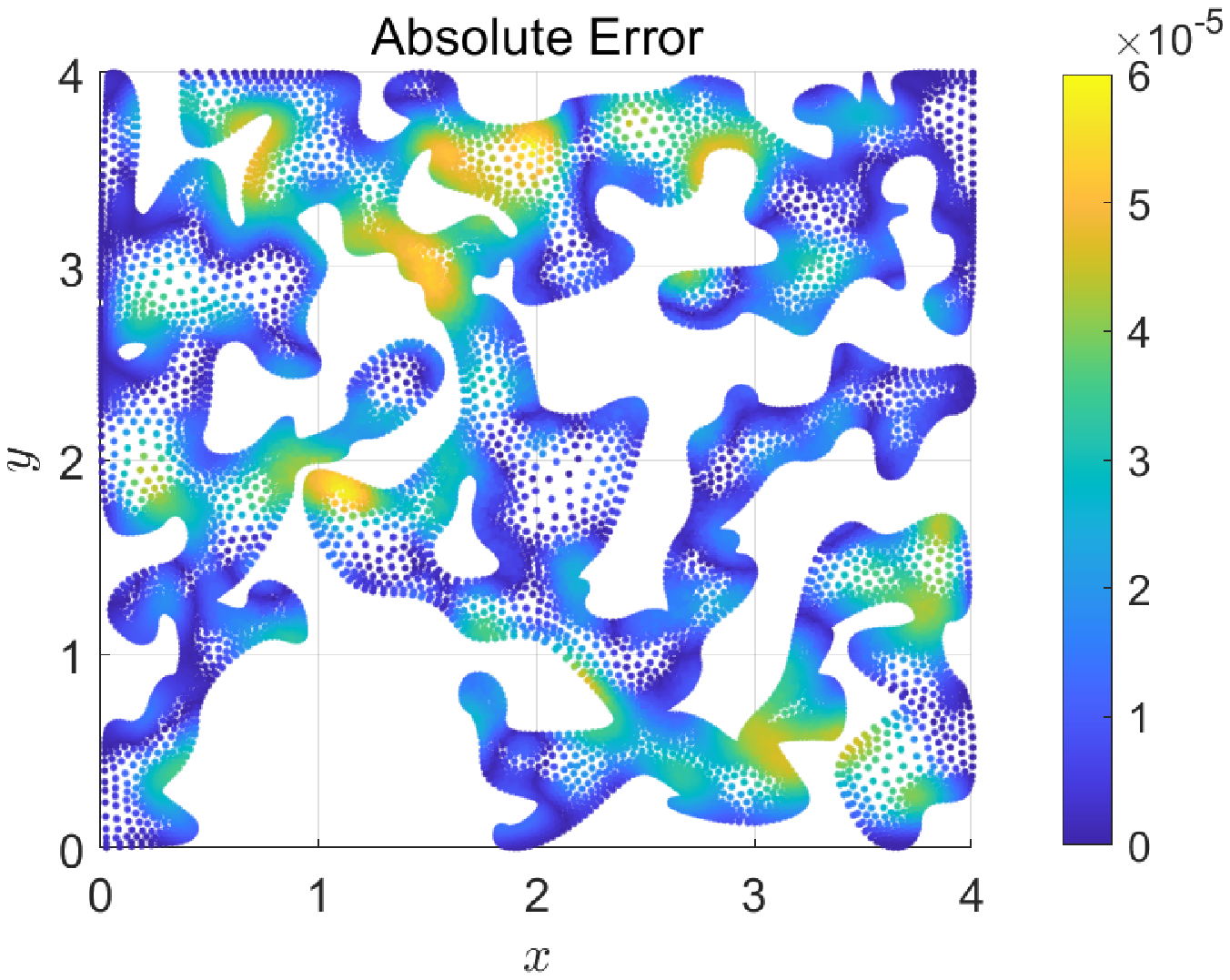}
    }
    \subfigure[MSE of PINN, SFPINN and SFHCPINN]{
    	\label{EX3_test_mse}
    	\includegraphics[scale=0.35]{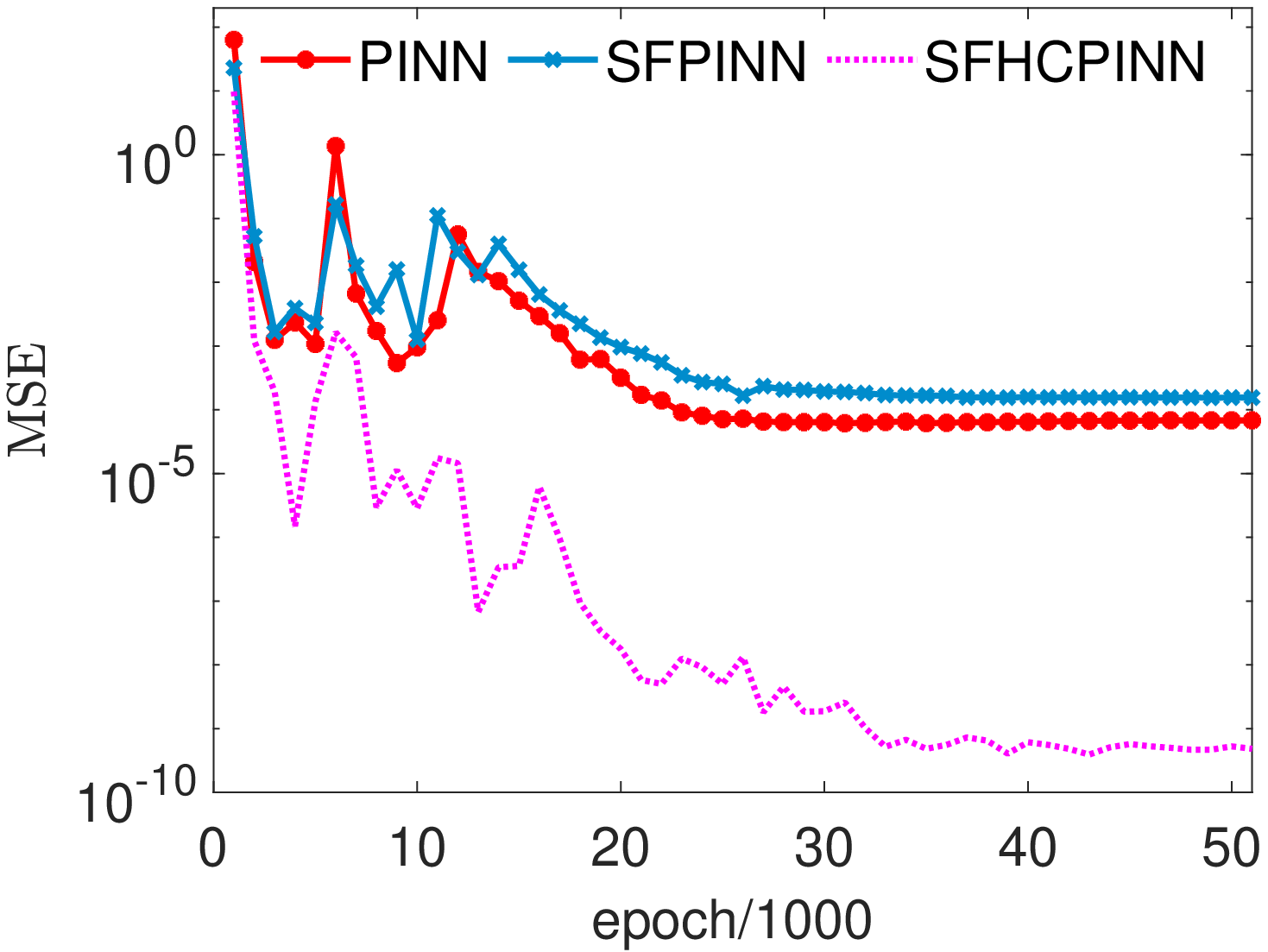}
    }
    \subfigure[REL of PINN, SFPINN and SFHCPINN]{
        \label{EX3_testrel}
        \includegraphics[scale=0.35]{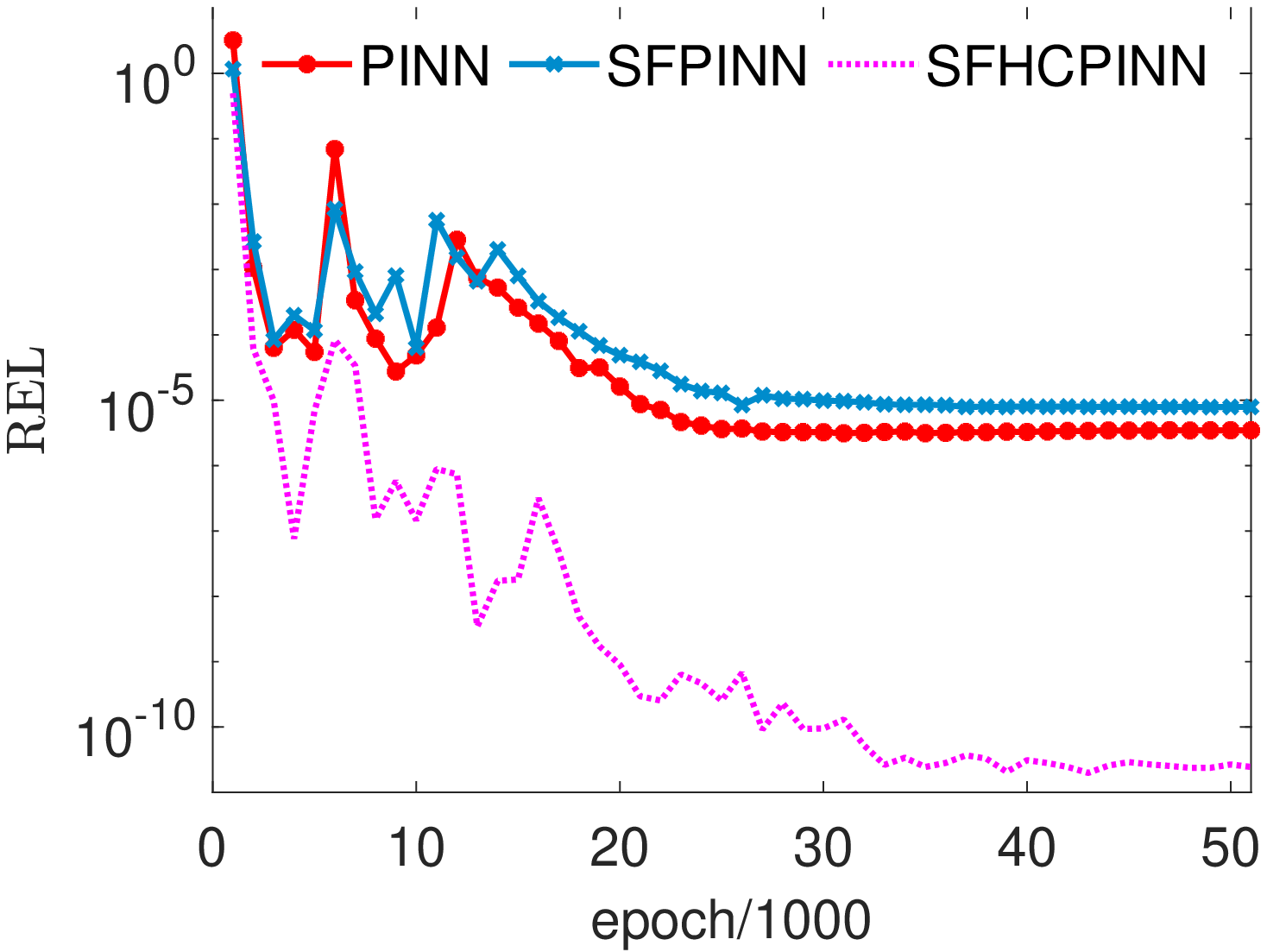}
    }
    \caption{Testing results for Example \ref{Smooth_Porous}.}
    \label{fig2smooth_porous_dirichlet2d}
\end{figure}

\begin{table}[H] 
\caption{MSE and REL of three models for Example \ref{Smooth_Porous} at $t=2.5$}
\label{table2smooth_porous_dirichlet2d}
\begin{center}
\setlength{\tabcolsep}{3pt}
    \begin{tabular}{llll} 
        \toprule 
        &constraint & MSE & REL \\
        \midrule 
          PINN     & soft &\textcolor{black}{$6.82\times10^{-5}$}  & \textcolor{black}{$3.46\times10^{-6}$}\\
          SFPINN   & soft &\textcolor{black}{$1.55\times10^{-4}$}  & \textcolor{black}{$7.88\times10^{-6}$}\\
          SFHCPINN & hard &\textcolor{black}{$4.83\times10^{-10}$} & \textcolor{black}{$2.45\times10^{-11}$}\\
        \bottomrule 
    \end{tabular}
\end{center}
\end{table}

The decreasing colors of the point-wise absolute error heatmap of three models in Figs.~\ref{ex3:pinnpwe} -- \ref{ex3:sfhcpinnpwe} and the data in Table \ref{table2smooth_porous_dirichlet2d} suggest that \textcolor{black}{the accuracy of SFHCPINN is higher than that of the PINN and SFPINN models}. After training using the sub-Fourier and hard boundary technique, the MSE of the model drops from $10^{-5}$ to $10^{-10}$, as compared to the normal PINN techniques. This is consistent with our study that utilizing a sub-Fourier structure during training may be more effective in addressing the spectrum bias induced by frequency items in the issue.
In addition, by incorporating a hard-constraint architecture, the model's accuracy jumps from $10^{-4}$ to $10^{-10}$ and has a faster convergence rate compared to SFPINN. This is owing to the hard-constraint model meeting the BCs before the training process.
In conclusion, SFHCPINN may outperform the baseline models under the two-dimensional ADE problem with Dirichlet BCs.
\end{example}

\begin{example}\label{Multiscale2D_Neumman}
In this example, we aim to approximate the solution of \eqref{TwoDim2ADE} with Neumann boundary for given 
regular domain $\Omega = [0, 1] \times [0, 1]$ and time range $[0,1]$. An exact solution is given by
\begin{equation}
u(x,t)=e^{-0.25t}[\sin(\pi x)\sin(\pi y)+0.1\sin(10\pi x)\sin(10\pi y)].
\end{equation}
then the boundary and initial functions can be easily obtained from the exact solution, we will omit them here. By performing careful calculations, one can determine the expression for the force side

The exact solution contains two frequency elements in which we are interested. We first define $D(x,y,t)=1-t$ and $G(x,y,t)=\sin(\pi x)\sin(\pi y)+0.1\sin(10\pi x)\sin(10\pi y)$ according to the BCs. All the model setups are identical to those in Example~\ref{Smooth_Porous}. 
 We train 50000 epochs for each model and in each epoch, we randomly generate 8000 initial points and 3000 boundary points from the interior of the defined domain $\Omega$ and Neumann boundary.
In addition,  we uniformly sample 16384 points in $[0,1] \times [0,1]$ at $t=0.5$ as the testing set to validate the feasibility of three models.

\begin{figure}[H]
    \centering
    \subfigure[Exact Solution of Example \ref{Multiscale2D_Neumman}]{
        \label{Exact:EX4}
        \includegraphics[scale=0.35]{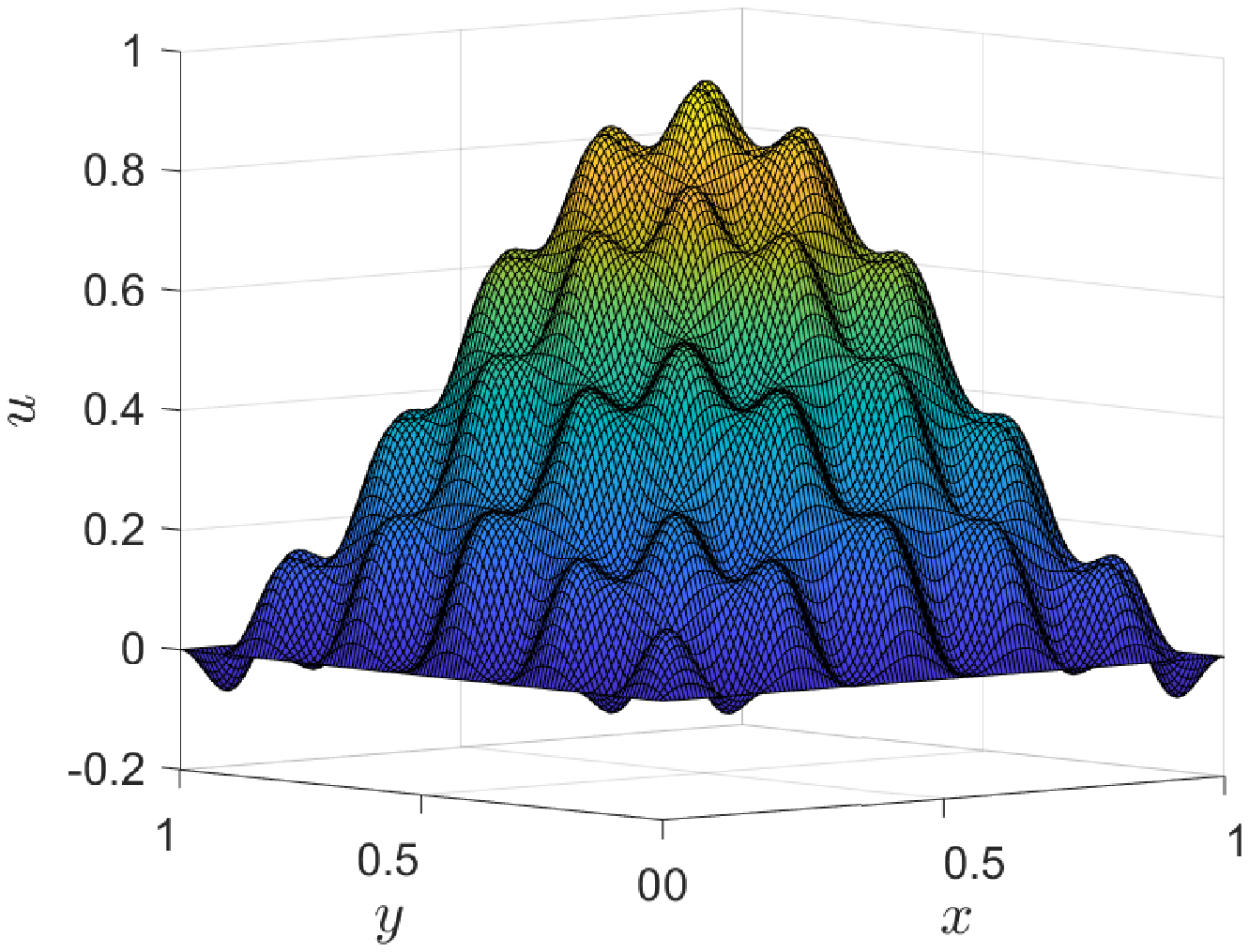}
    }        
    \subfigure[Point-wise error for PINN]{
        \label{ex4:pinnpwe}
        \includegraphics[scale=0.35]{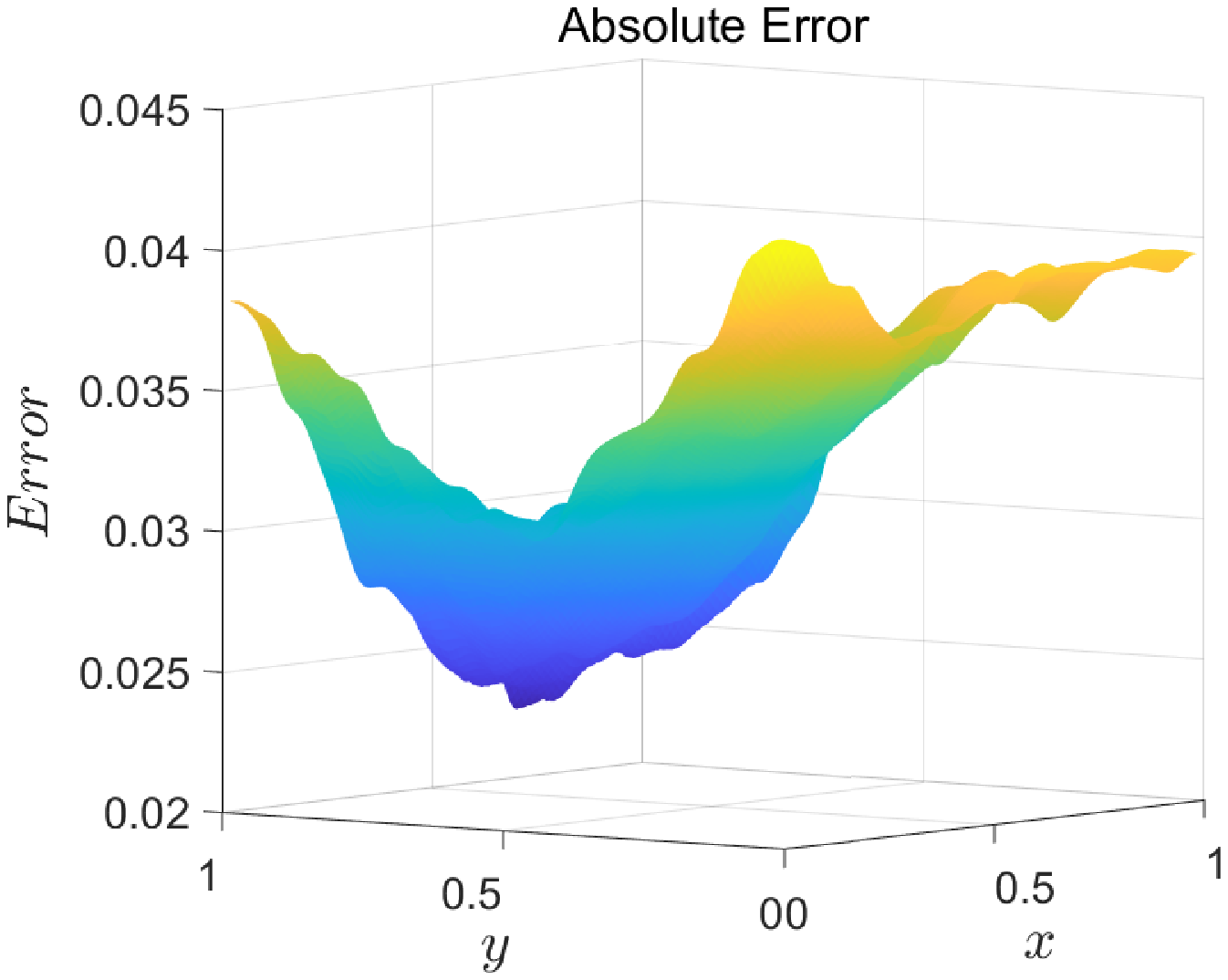}
    }
    \subfigure[Point-wise error for SFPINN]{
        \label{ex4:sfpinnpwe}
        \includegraphics[scale=0.35]{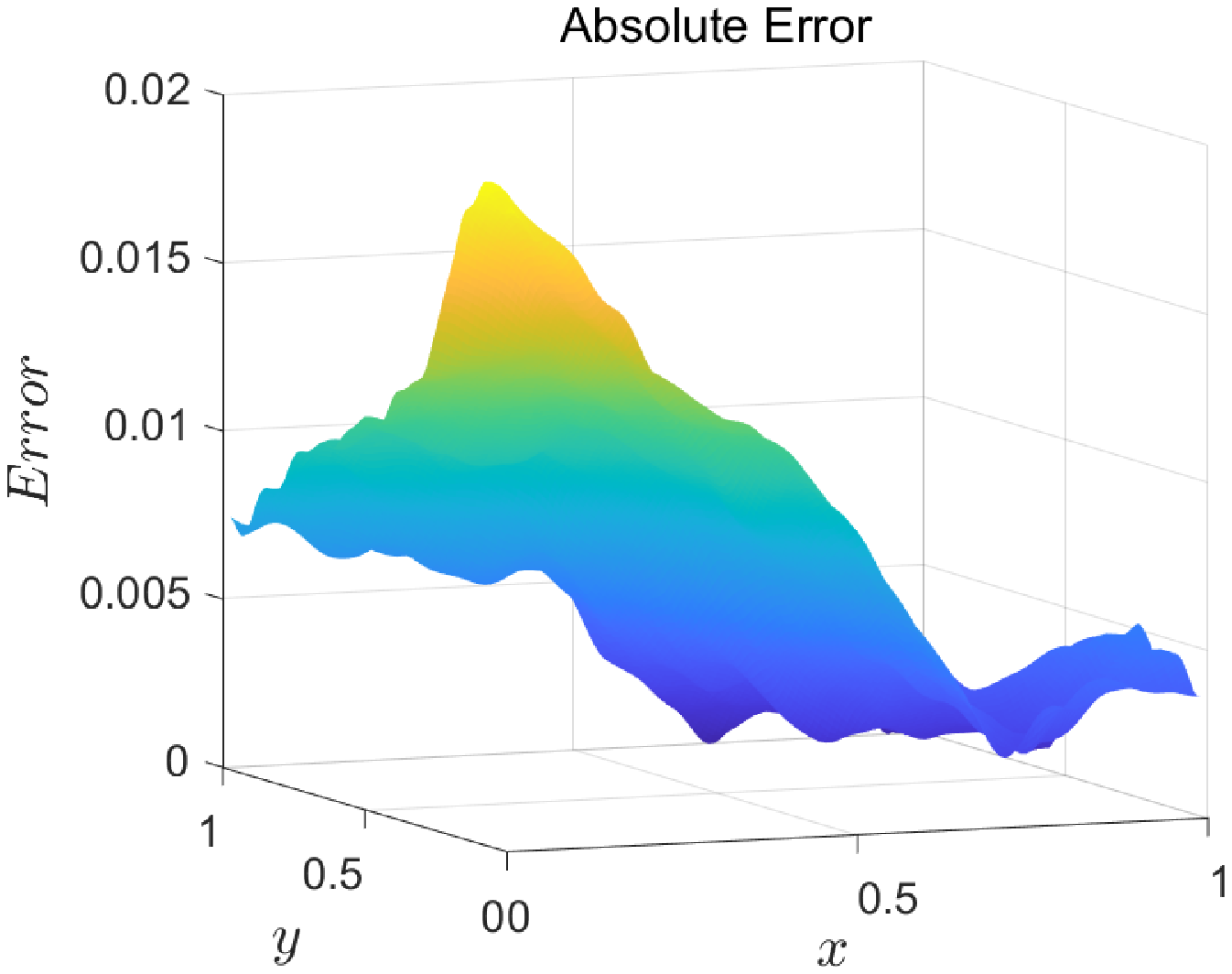}
    }
    \subfigure[Point-wise error for SFHCPINN]{
        \label{ex4:sfhcpinnpwe}
        \includegraphics[scale=0.35]{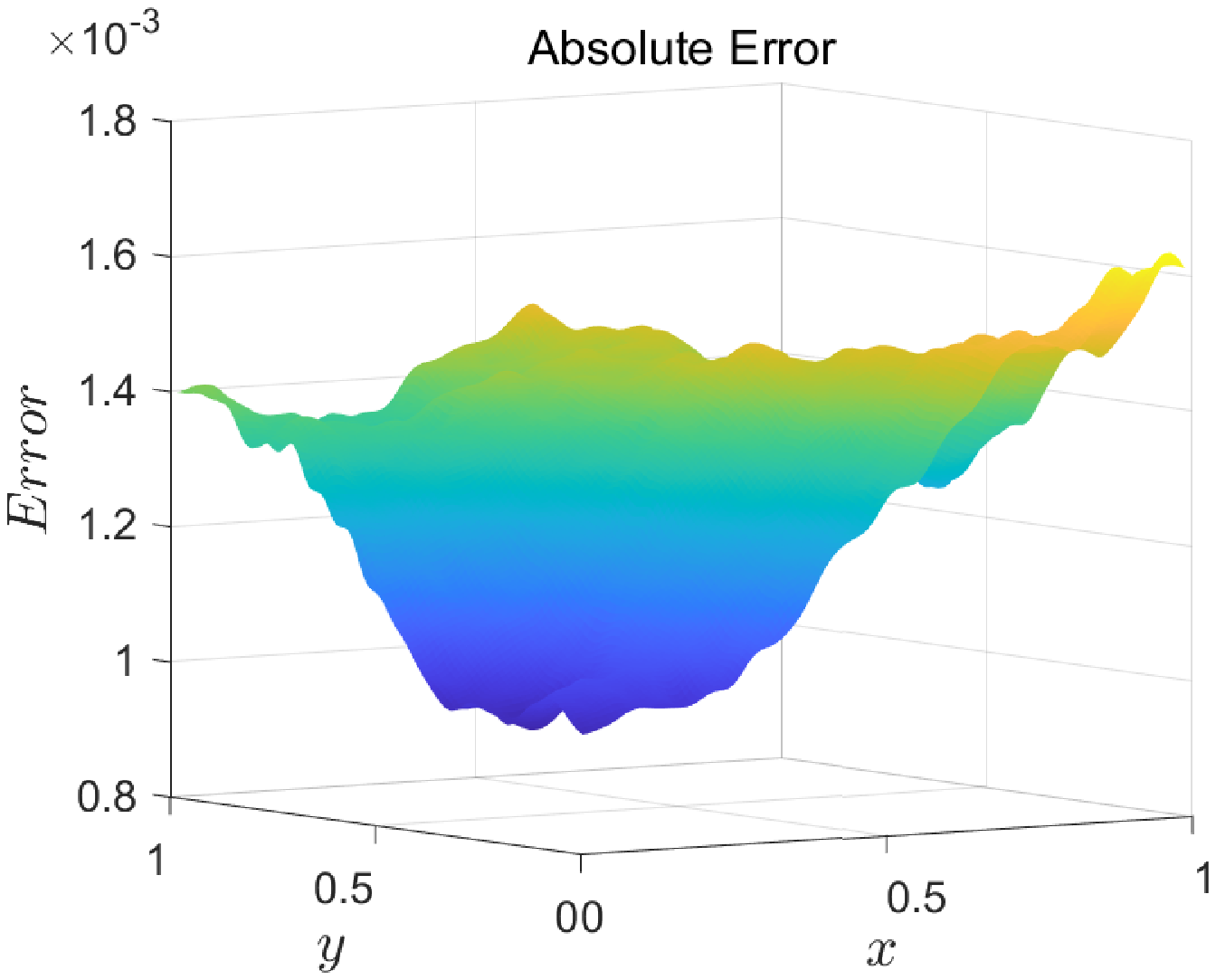}
    }
    \subfigure[MSE of PINN, SFPINN and SFHCPINN]{
        \label{EX4_test_mse}
        \includegraphics[scale=0.35]{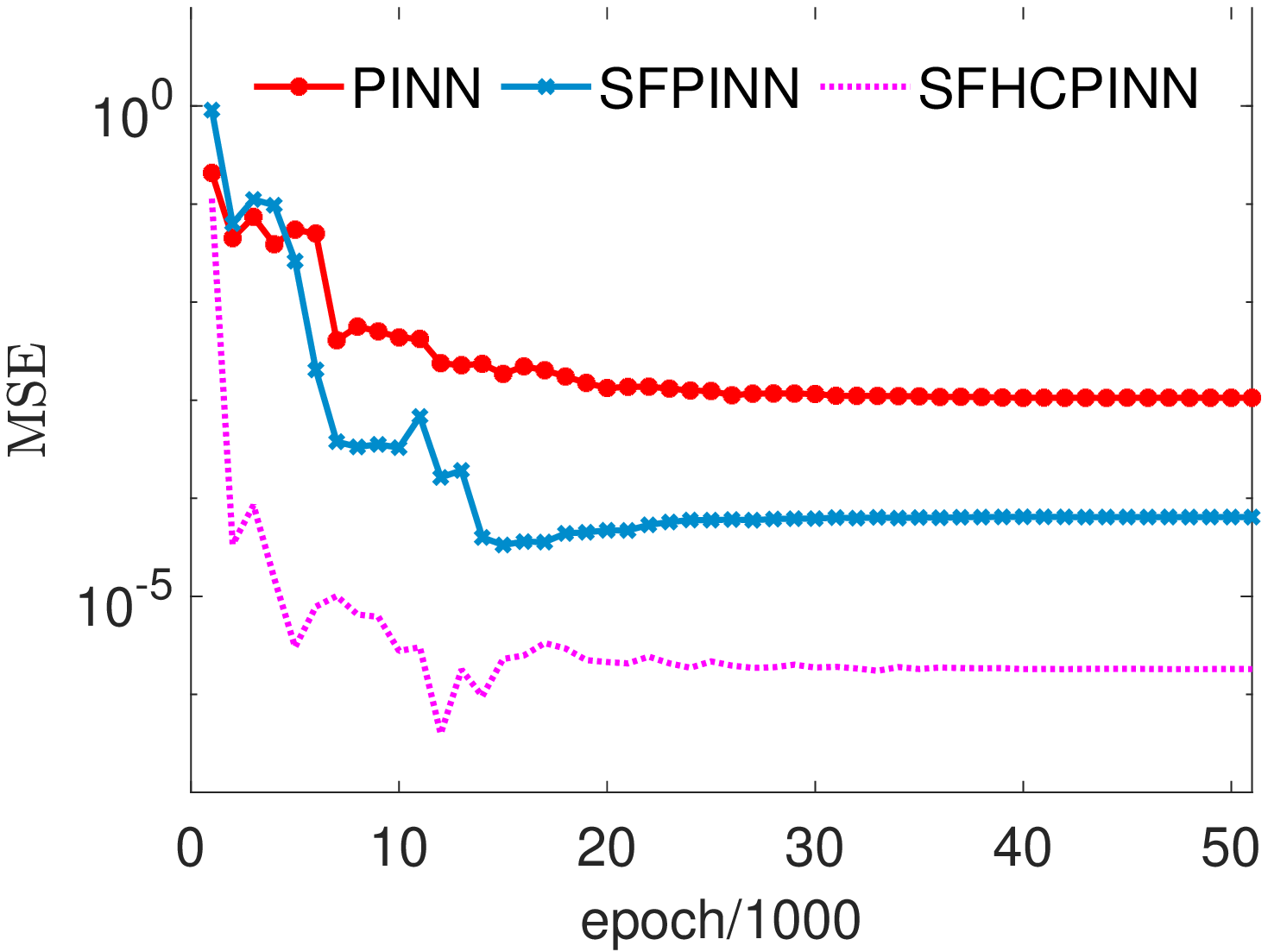}
	}
    \subfigure[REL of PINN, SFPINN and SFHCPINN]{
        \label{EX4_testrel}
        \includegraphics[scale=0.35]{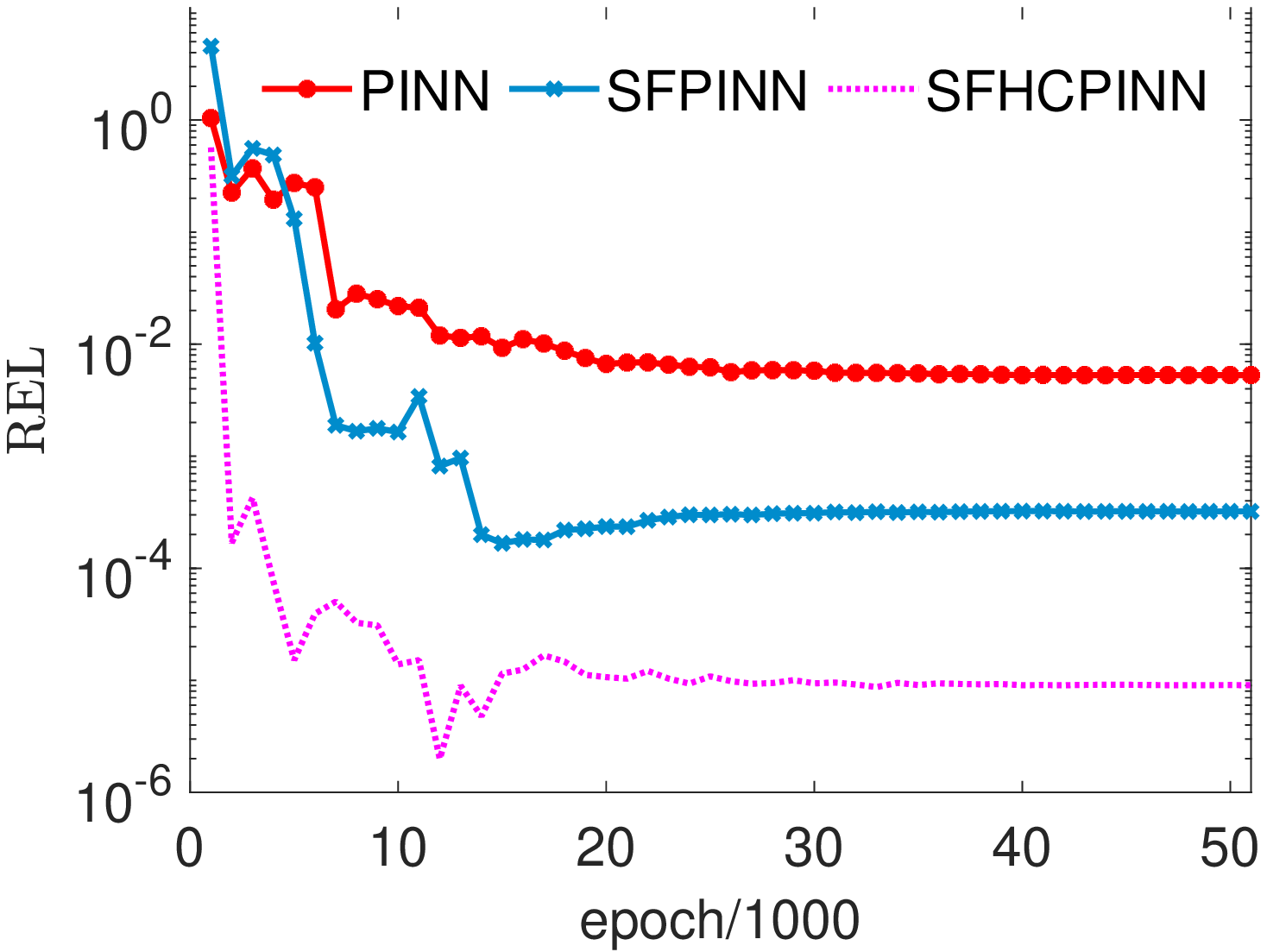}
	}
    \caption{Testing results for Example \ref{Multiscale2D_Neumman}.}
    \label{E4}
\end{figure}

\begin{table}[H]
\caption{MSE and REL of SFHCPINN, SFPINN, and PINN for Example \ref{Multiscale2D_Neumman} at $t=0.5$}
\label{tab4}
\begin{center}
\setlength{\tabcolsep}{3pt}
\resizebox{!}{!}{  
    \begin{tabular}{llll} 
        \toprule 
        &constraint & MSE & REL \\
        \midrule 
          PINN     & soft &$1.06\times 10^{-3}$ &$5.31\times 10^{-3}$\\
          SFPINN   & soft &$6.43\times 10^{-5}$ &$3.22\times 10^{-4}$\\
          SFHCPINN & hard &$1.81\times 10^{-6}$ &$9.05\times 10^{-6}$\\
        \bottomrule 
    \end{tabular}
}
\end{center}
\end{table}

 The results in Fig.~\ref{E4} show that employing the subnetwork structure and the Fourier expansion on the multi-scaled input the SFHCPINN has lower MSE and relative errors as well as a faster convergence rate than the standard PINN. \textcolor{black}{The accuracy of PINN, SFPINN, and SFHCPINN rises in that order.} In addition, by decomposing the solution, the solution automatically meets the IC, enabling the SFHCPINN to be tuned to achieve greater precision regarding the two-dimensional issues with Neumann boundaries. This example exhibits the feasibility and excellent accuracy of SFHCPINN in solving the two-dimensional ADE under the Neumann BCs, whereas the performance of PINN and SFPINN is only ordinary.
\end{example}

 \subsubsection{Three-dimensional ADE}
 We consider the following spatio-temporal ADE \eqref{eq0101} with Dirichlet boundary, then obtain their numerical solutions on given three-dimensional regular and irregular domain, respectively. It is
  \begin{equation}\label{ThreeDim2ADE}
 \frac{\partial u}{\partial t}-\Delta u+ \hat{p} \frac{\partial u}{\partial x}+ \hat{q} \frac{\partial u}{\partial y}+ \hat{r} \frac{\partial u}{\partial z}=f, ~~(x, y, z, t) \in \Omega \times(t_0, T]
 \end{equation}
 where $\Delta$ is the Laplace operator and $\hat{p}>0$, $\hat{q}>0$, $\hat{r}>0$ are the corresponding diffusion coefficients for different directions.  $\Omega$ is the interested domain, $t_0\geqslant0$ is the initial time and $T>0$ is the end time, and $u$ is the function to be solved.

 \begin{example}\label{SmoothCase_porous3D}
     \textcolor{black}{Let us now approximate the solution of \eqref{ThreeDim2ADE} with $\hat{p}=1$, $\hat{q}=1$ and $\hat{r}=1$ in time interval $[0,5]$ and a unit cubic domain $\Omega=[0,1]\times[0,1]\times[0,1]$ with 1 big hole (cyan) and 8 smaller holes (red and black), see Fig.~\ref{fig:ThreeDim_Holes}.}  
 \begin{figure}[H]
 	\centering
	 \subfigure[Unit cubic domain with some holes]{
	 	\label{Holes_3D}
	 	\includegraphics[scale=0.4]{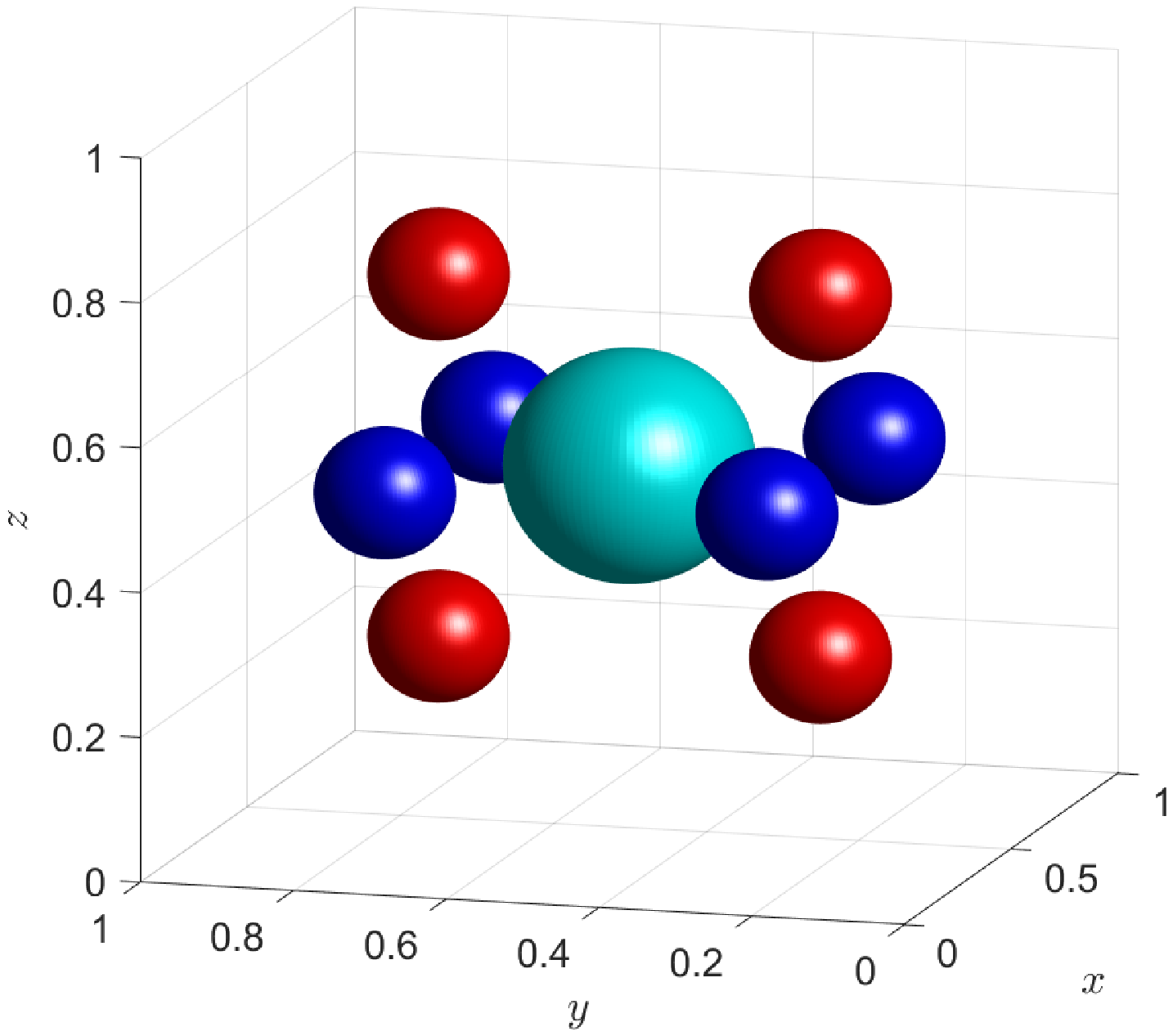}
	 }        
	 \subfigure[An $xoy$ cutting plane of domain at $z=0.5$.]{
	 	\label{Sile2_Holes_3DD}
	 	\includegraphics[scale=0.4]{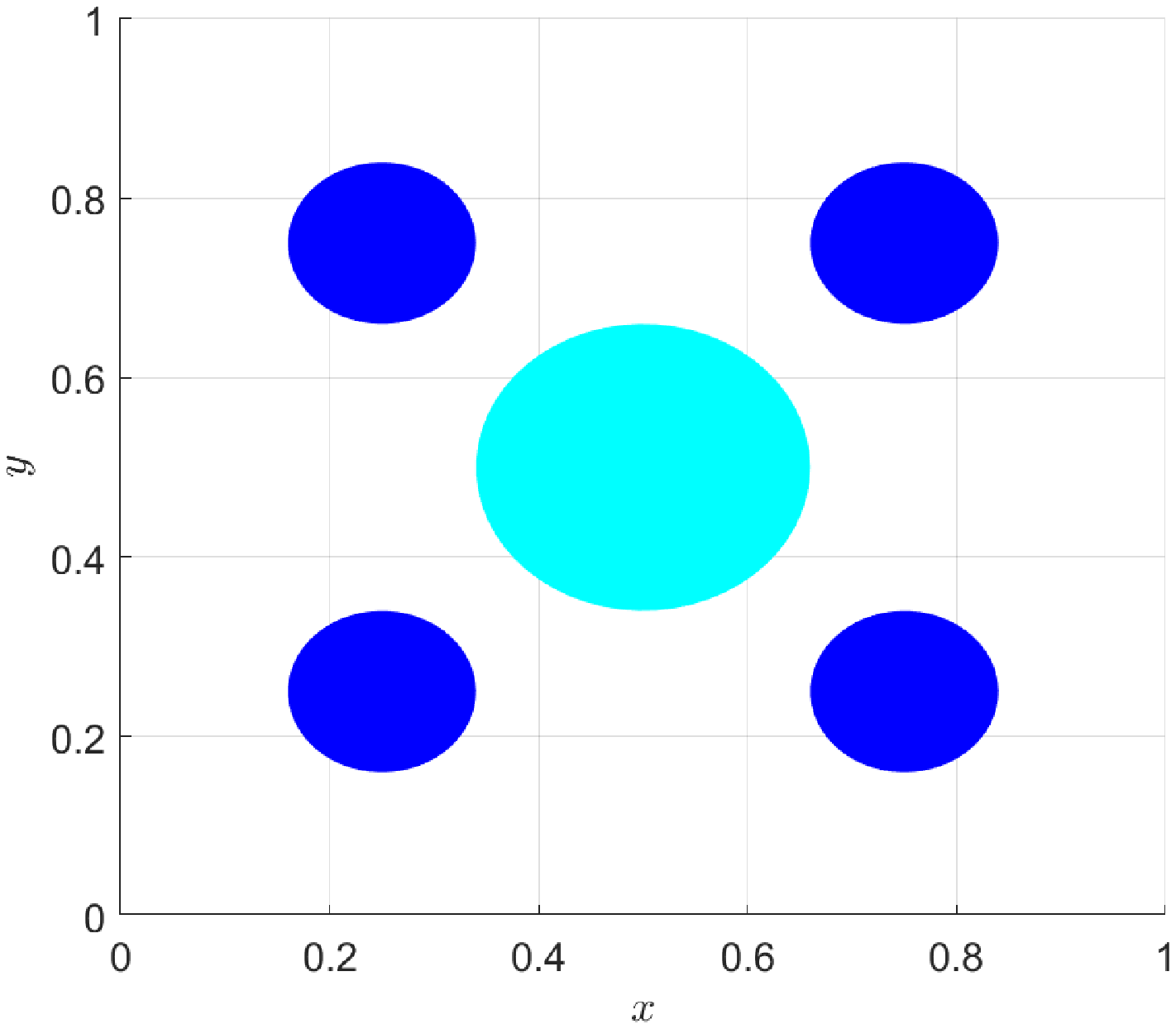}
	 }
 	\caption{The interested domain and its cutting plane for Example \ref{SmoothCase_porous3D}.}
 	\label{fig:ThreeDim_Holes}
 \end{figure}
\textcolor{black}{A homogeneous  solution is given by $u=e^{-0.25t}\sin(\pi x)\sin(\pi y) \sin(\pi z)$, which prescribes the corresponding boundary and initial conditions and right-side $f(x,y,z,t)$ of \eqref{ThreeDim2ADE}. 
 To solve this problem, we define the distance function $\displaystyle D(x,y,z,t)=\frac{1}{5}xyz(1-x)(1-y)(1-z)t$ on the interested domain and the smooth extension function $G(x,y,z,t)=\sin(\pi x )\sin(\pi y) \sin(\pi z)$ according to prescribed boundary and initial conditions.}
 The network setup, optimizer, the hyperparameters of three models are unified with that in Section~\ref{trainingsetup}. 
 In each training epoch, PINN and SFPINN are trained with $N_R=15,000$ collocation points, $N_B=4,000$ boundary points, and $N_I=3,000$ initial points while the SFHCPINN is trained with $N_R=15,000$ collocation points. We train the above models for 50,000 epochs and 
 test them on 80000 random points sampled from $[0,1] \times [0,1]$ at $z=0.5$ and $t=0.5$.

 \begin{figure}[!ht]
     \centering
     \subfigure[Exact Solution of Example \ref{SmoothCase_porous3D}]{
         \label{Exact:Smooth_porous3D}
         \includegraphics[scale=0.35]{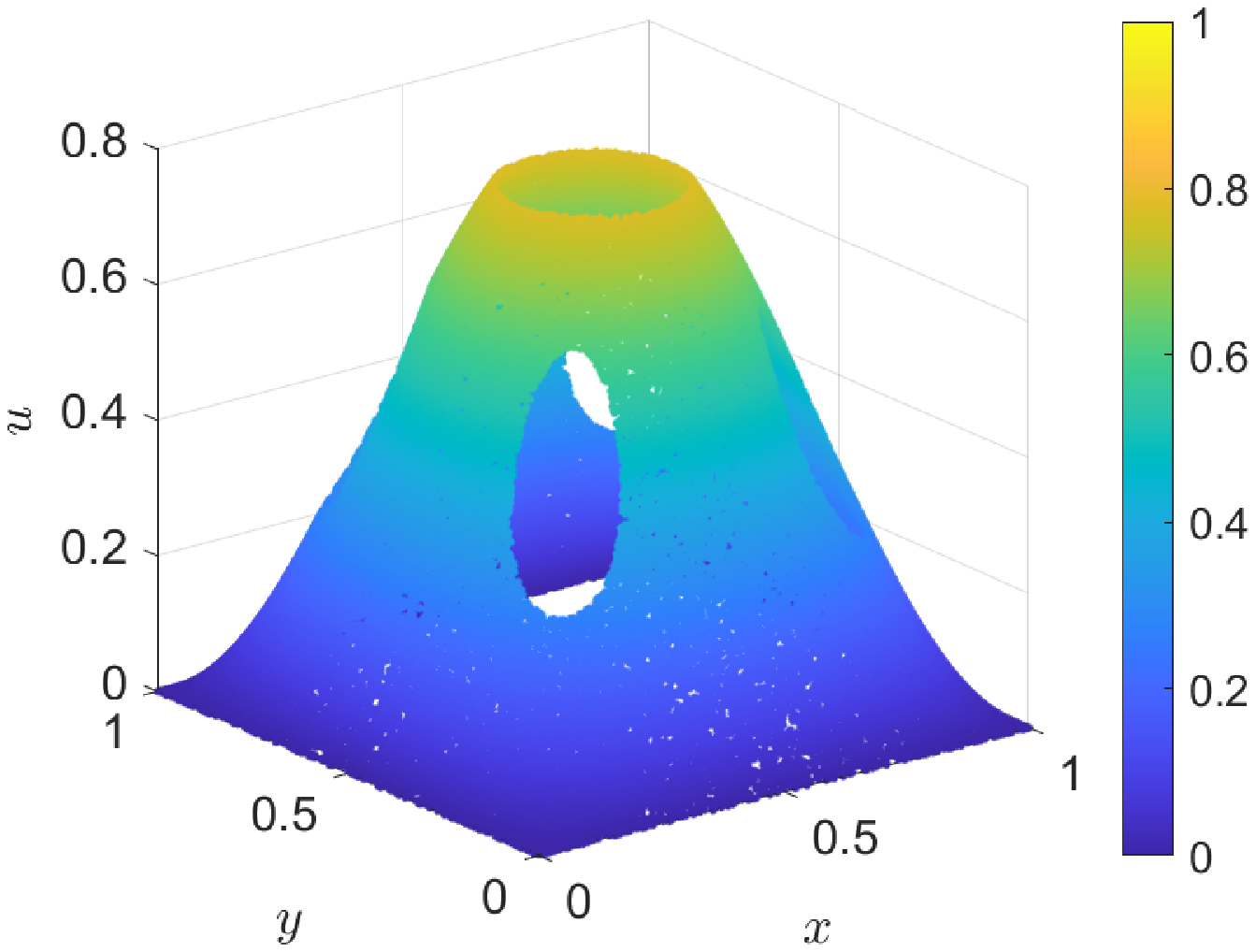}
     }        
     \subfigure[Point-wise error for PINN]{
         \label{PINNPWE:Smooth_porous3D}
         \includegraphics[scale=0.35]{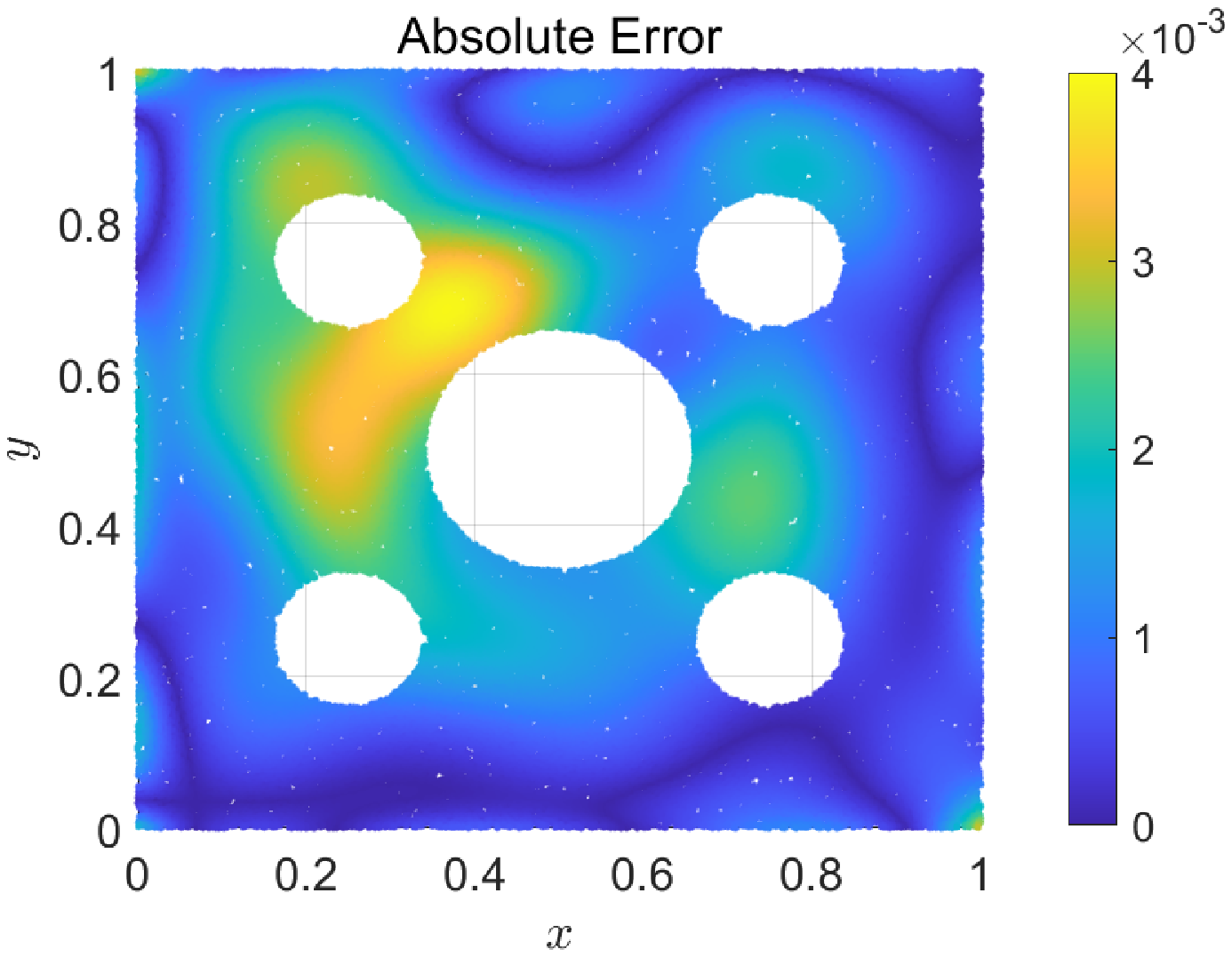}
     }
     \subfigure[Point-wise error for SFPINN]{
         \label{SFPINNPWE:Smooth_porous3D}
         \includegraphics[scale=0.35]{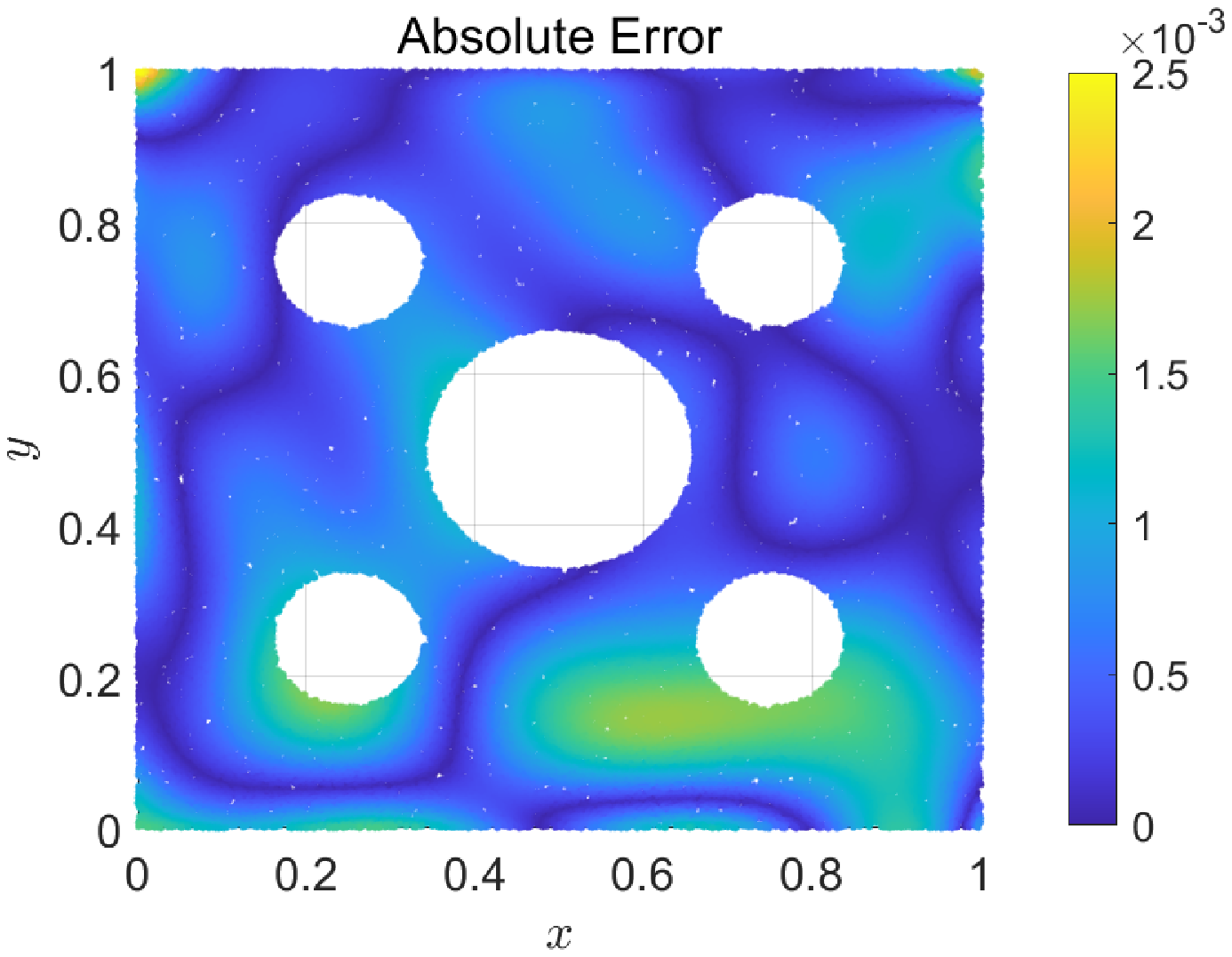}
     }
     \subfigure[Point-wise error for SFHCPINN]{
         \label{SFHCPINNPWE:Smooth_porous3D}
         \includegraphics[scale=0.35]{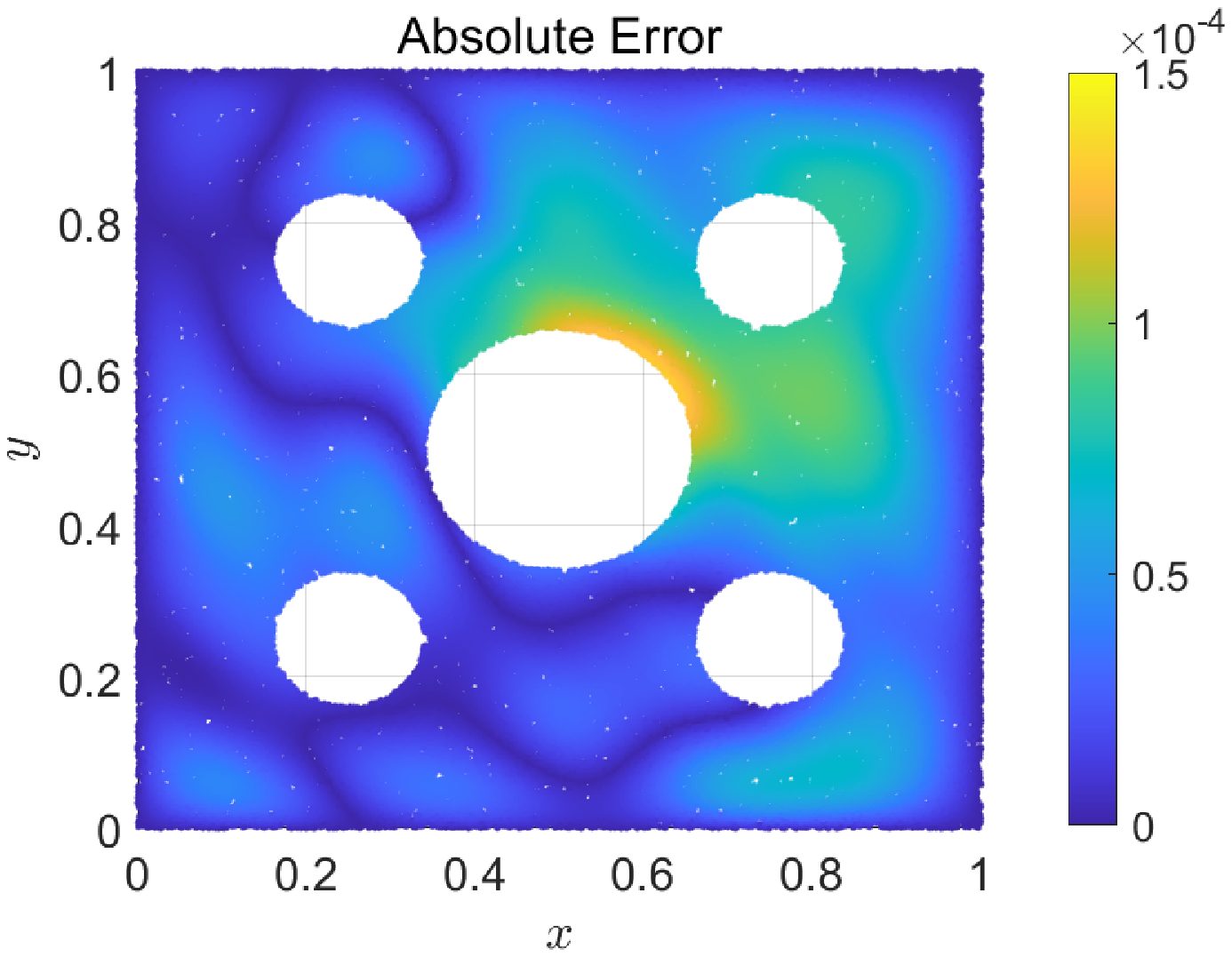}
     }
     \subfigure[MSE of PINN, SFPINN and SFHCPINN]{
         \label{Test_MSE:Smooth_porous3D}
         \includegraphics[scale=0.35]{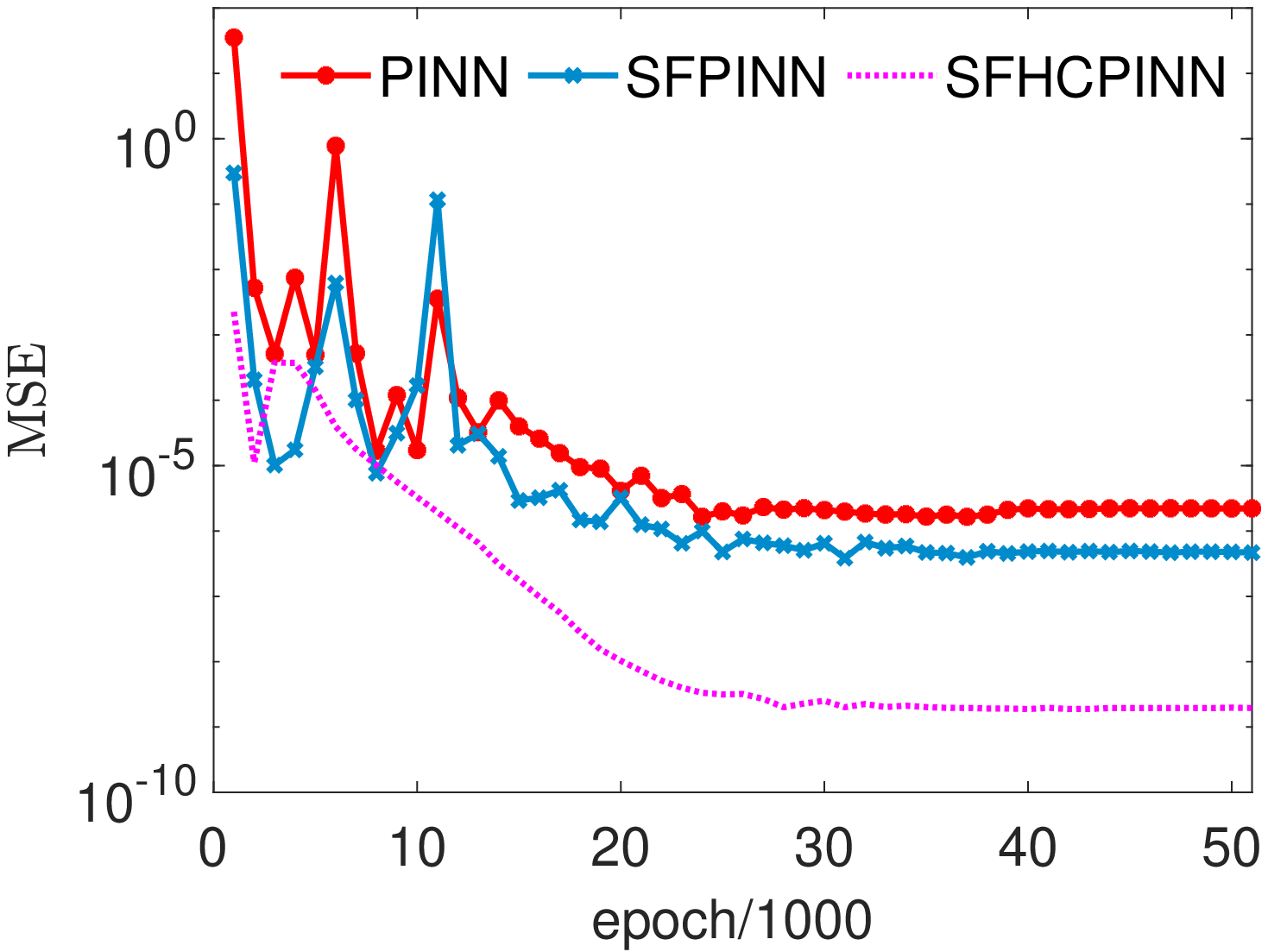}
 	}
     \subfigure[REL of PINN, SFPINN and SFHCPINN]{
         \label{Test_REL:Smooth_porous3D}
         \includegraphics[scale=0.35]{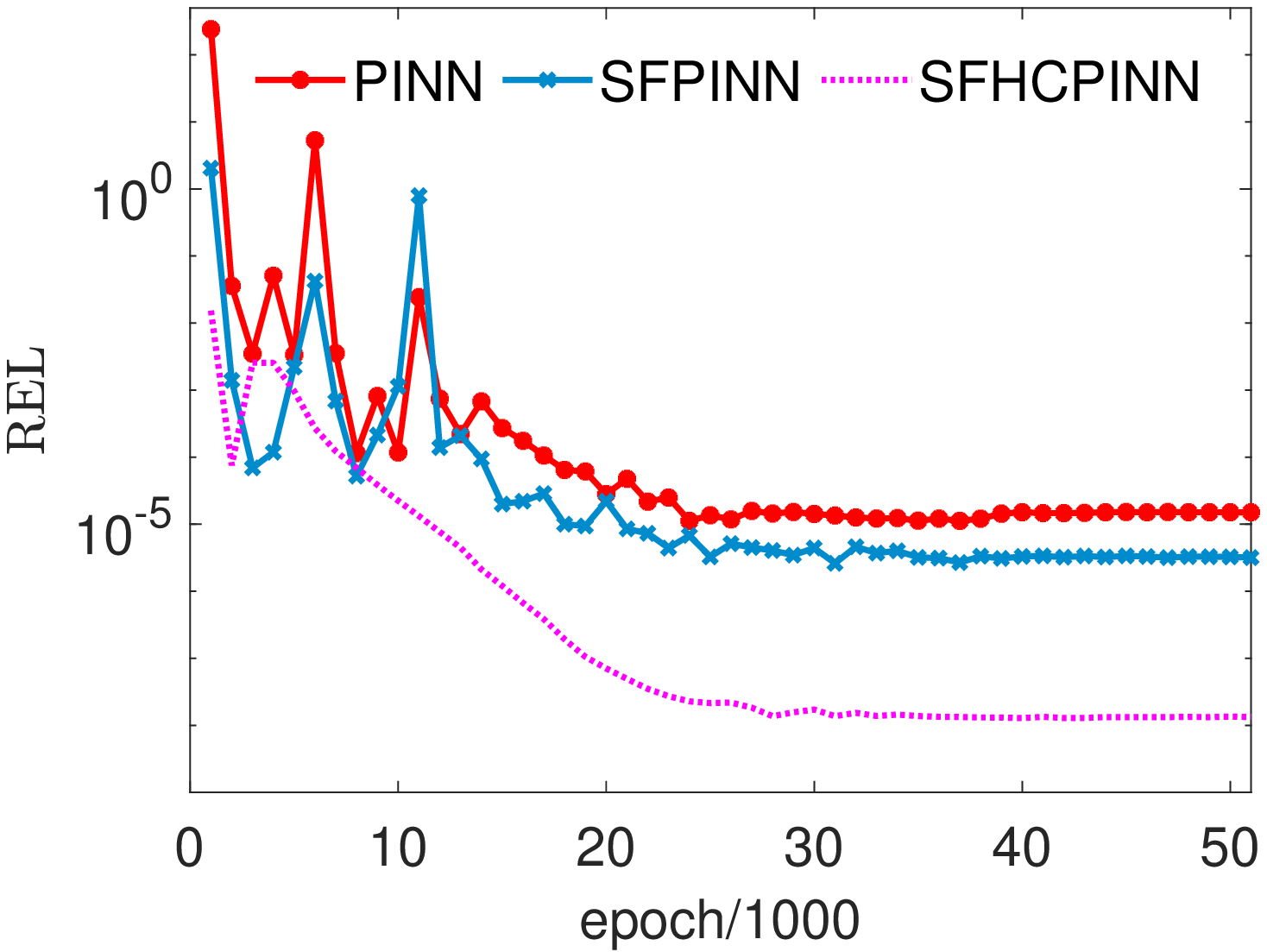}
 	}
     \caption{Testing results for Example \ref{SmoothCase_porous3D} at $z=0.5$ and $t=0.5$.}
 \end{figure}

\begin{table}[!ht] 
 \caption{MSE and REL of SFHCPINN, SFPINN, and PINN for Example \ref{SmoothCase_porous3D} at $z=0.5$ and $t=0.5$}
 \label{Table:Smooth_porous3D}
 \begin{center}
 \setlength{\tabcolsep}{3pt}
 \resizebox{!}{!}{  
     \begin{tabular}{llll} 
         \toprule 
         &constraint & MSE & REL \\
         \midrule 
           PINN      & soft &\textcolor{black}{$2.21\times 10^{-6}$}  &\textcolor{black}{$1.51\times 10^{-5}$}\\
           SFPINN    & soft &\textcolor{black}{$2.21\times 10^{-7}$}  &\textcolor{black}{$3.18\times 10^{-6}$}\\
           SFHCPINN  & hard &\textcolor{black}{$1.96\times 10^{-9}$}  &\textcolor{black}{$1.96\times 10^{-8}$}\\
         \bottomrule 
     \end{tabular}
 }
 \end{center}
 \end{table}

 Based on the data in Table~\ref{Table:Smooth_porous3D} and the heatmap of three models, it is evident that the normal PINN is less accurate than the other two models with hard constraints and/or sub-Fourier architecture when solving 3D Dirichlet problem in a porous domain.
 In addition, we can deduce from Figs.~\ref{Test_MSE:Smooth_porous3D} and \ref{Test_REL:Smooth_porous3D} that SFHCPINN with hard-constraint design has lower initial errors, faster convergence rate, and higher precision. This is also consistent with our analysis that models with hard constraints may enhance the performance of PINN because they naturally satisfy the BCs and transform the problem into a simpler optimization problem without additional physics constraints.
 \textcolor{black}{In conclusion, SFHCPINN retains its high precision, convergence rate, and robust stability in solving three-dimensional ADE problems for irregular domains.}
 \end{example}

\begin{example}\label{MultiscaleCase_Spere3D}
\textcolor{black}{We obtain the numerical approximation of \eqref{ThreeDim2ADE} with $\hat{p}=1$, $\hat{q}=1$ and $\hat{r}=1$ in time interval $[0,5]$ and a closed domain $\Omega$ restrained by two spherical surfaces with radius $r_1=0.1$ and $r_2=1.0$, respectively. Fig. \ref{fig:domain} depicts a cutting plane parallel to $xoz$ for our interested domain.}
\begin{figure}
    \centering
    \includegraphics[scale=0.45]{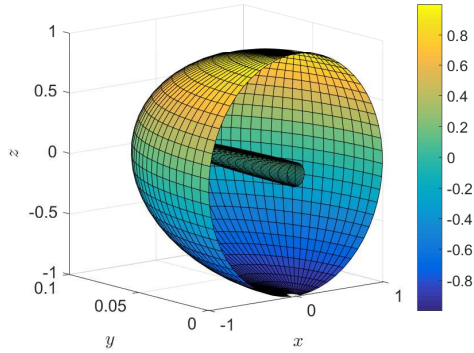}
    \caption{A $xoz$ cutting plane of interested domain.}
    \label{fig:domain}
\end{figure}
\textcolor{black}{An exact solution with two frequencies component is given by}
\begin{equation}\label{exact_solu2spere}
     \textcolor{black}{u(x, y,z, t) = e^{-0.25t}[\sin(\pi x)\sin(\pi y)\sin(\pi z)+0.1\sin(10\pi x)\sin(10\pi y)\sin(10\pi z)],~~~(x, y,z, t) \in \Omega \times [0, 1],}
\end{equation}
\textcolor{black}{it will naturally lead to the determination of the boundary condition and initial condition. By careful calculations, one can obtain the force side $f(x,y,z,t)$.}

\textcolor{black}{According to the BCs, we define the distance function $D(x,y,z,t)=\frac{(r-r_1)(r_2-r)}{(r_2-r_1)^2}t$ with $r=\sqrt{x^2+y^2+z^2}$ and the extension function $G(x,y,z,t)=\sin(\pi x)\sin(\pi y)\sin(\pi z)+0.1\sin(10\pi x)\sin(10\pi y)\sin(10\pi z)$.
The network setup, optimizer, the hyperparameters of three models are unified with that in Section~\ref{trainingsetup}. 
In each epoch, PINN and SFPINN are trained with $N_R=15,000$ collocation points, $N_B=4,000$ boundary points, and $N_I=3,000$ initial points while the SFHCPINN is trained with $N_R=15,000$ collocation points. We train the above models for 50,000 epochs and test them on 8000 mesh grid points on the spherical surface with radius $r=0.45$ and $t=0.5$. }
\begin{figure}[H]
    \centering
    \subfigure[Exact Solution of Example \ref{MultiscaleCase_Spere3D}]{
        \label{Exact:Multiscale_Spere3D}
        \includegraphics[scale=0.35]{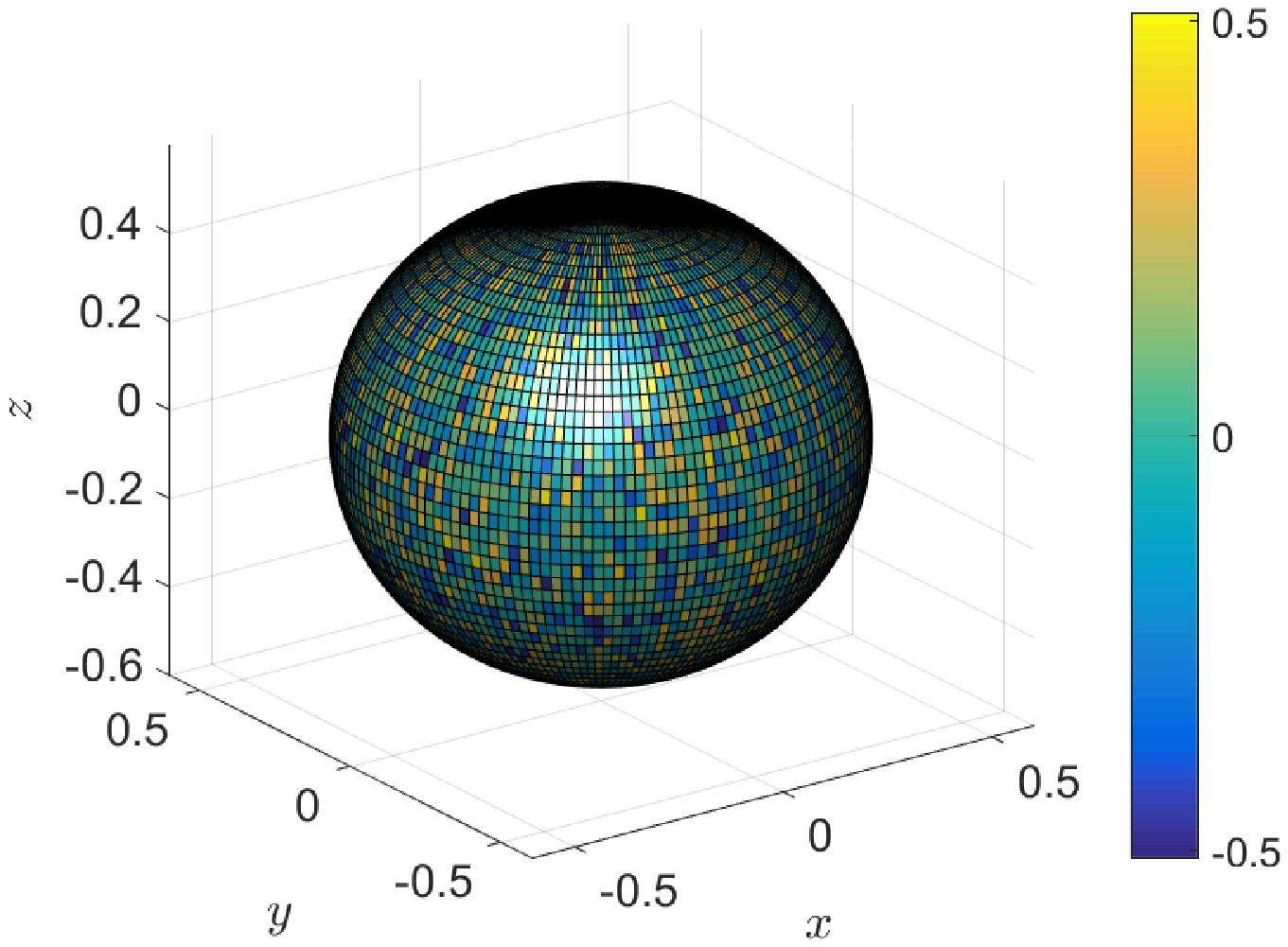}
    }        
    \subfigure[Point-wise error for PINN]{
        \label{PINNPWE:Multiscale_Spere3D}
        \includegraphics[scale=0.35]{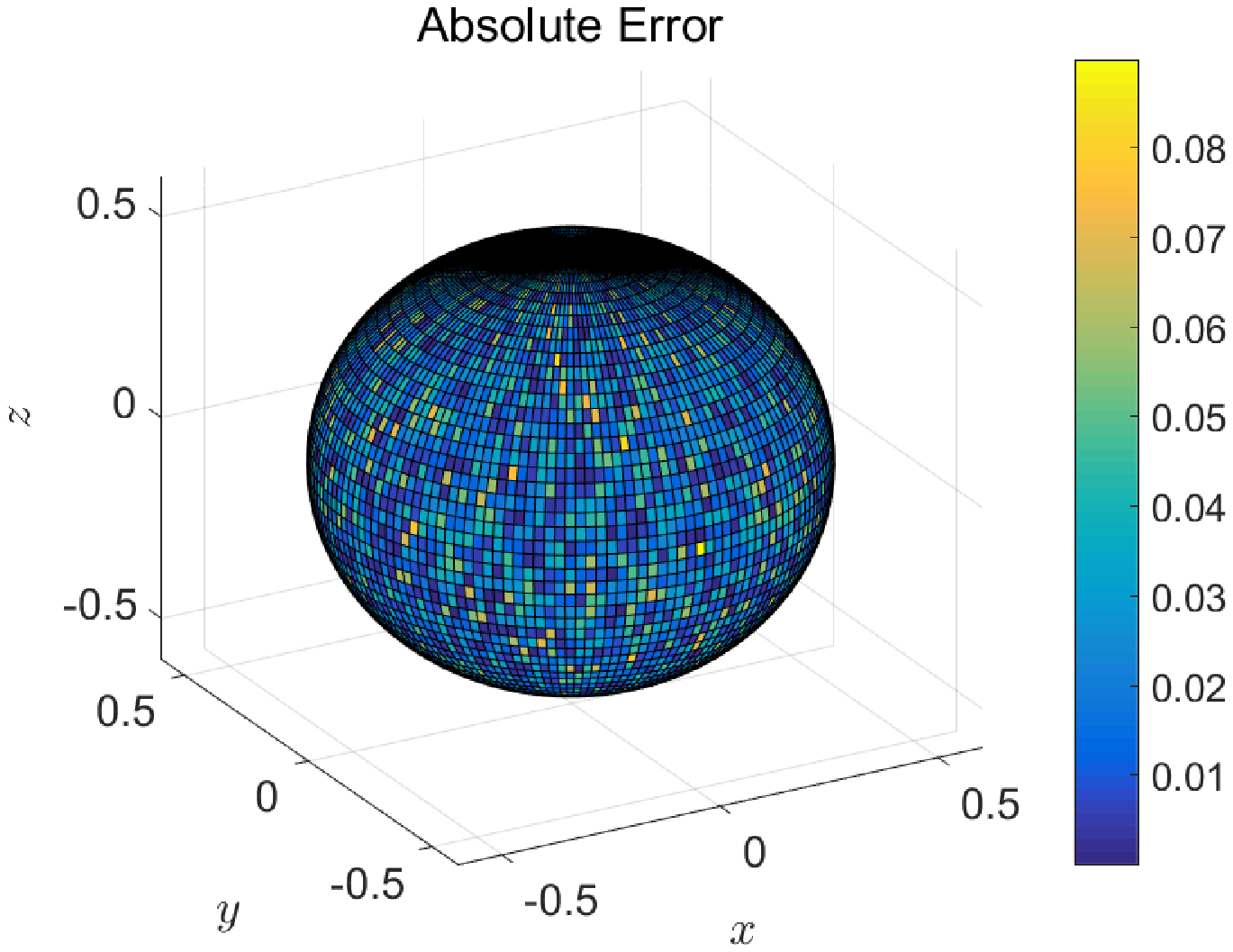}
    }
    \subfigure[Point-wise error for SFPINN]{
        \label{SFPINNPWE:Multiscale_Spere3D}
        \includegraphics[scale=0.35]{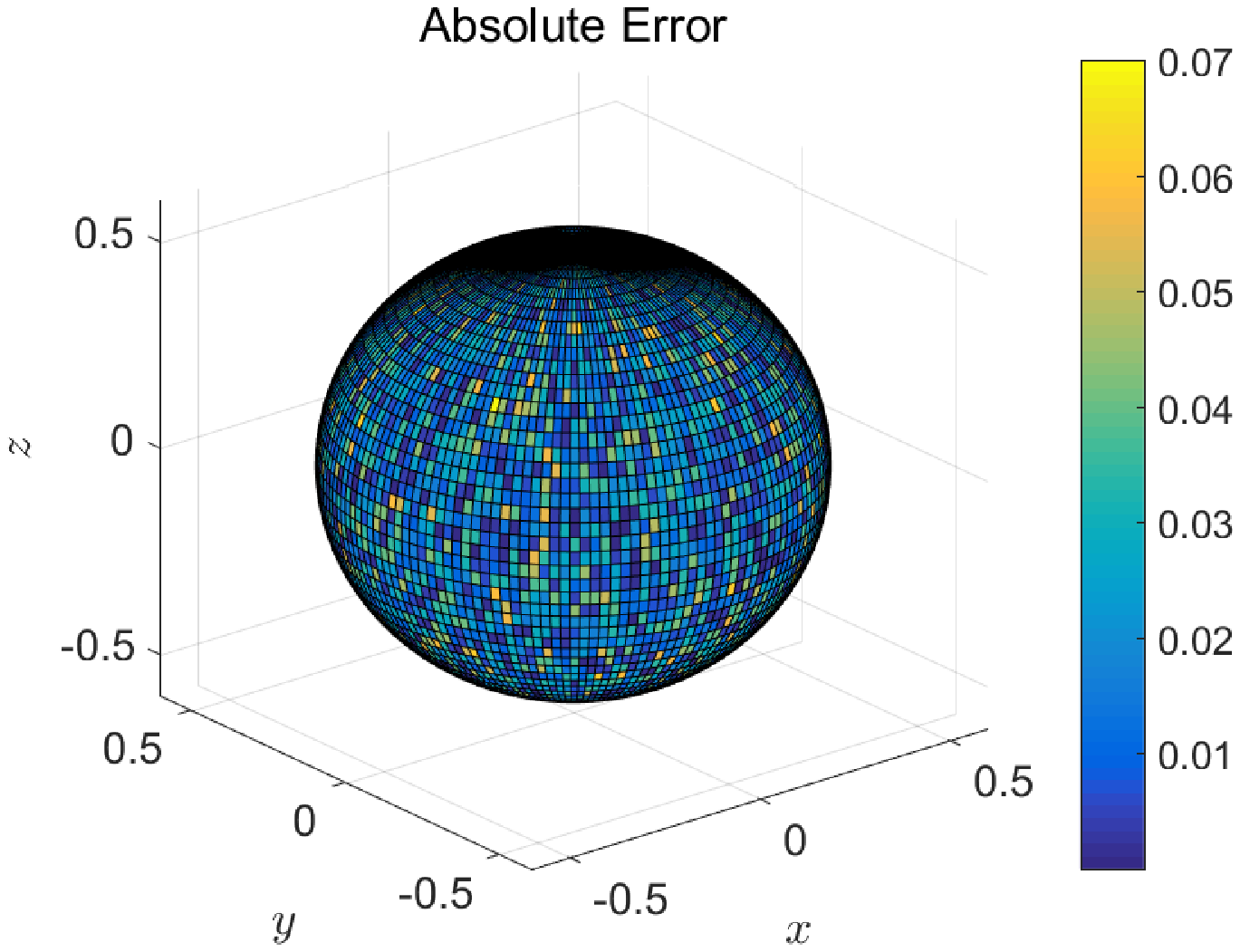}
    }
    \subfigure[Point-wise error for SFHCPINN]{
        \label{SFHCPINNPWE:Multiscale_Spere3D}
        \includegraphics[scale=0.35]{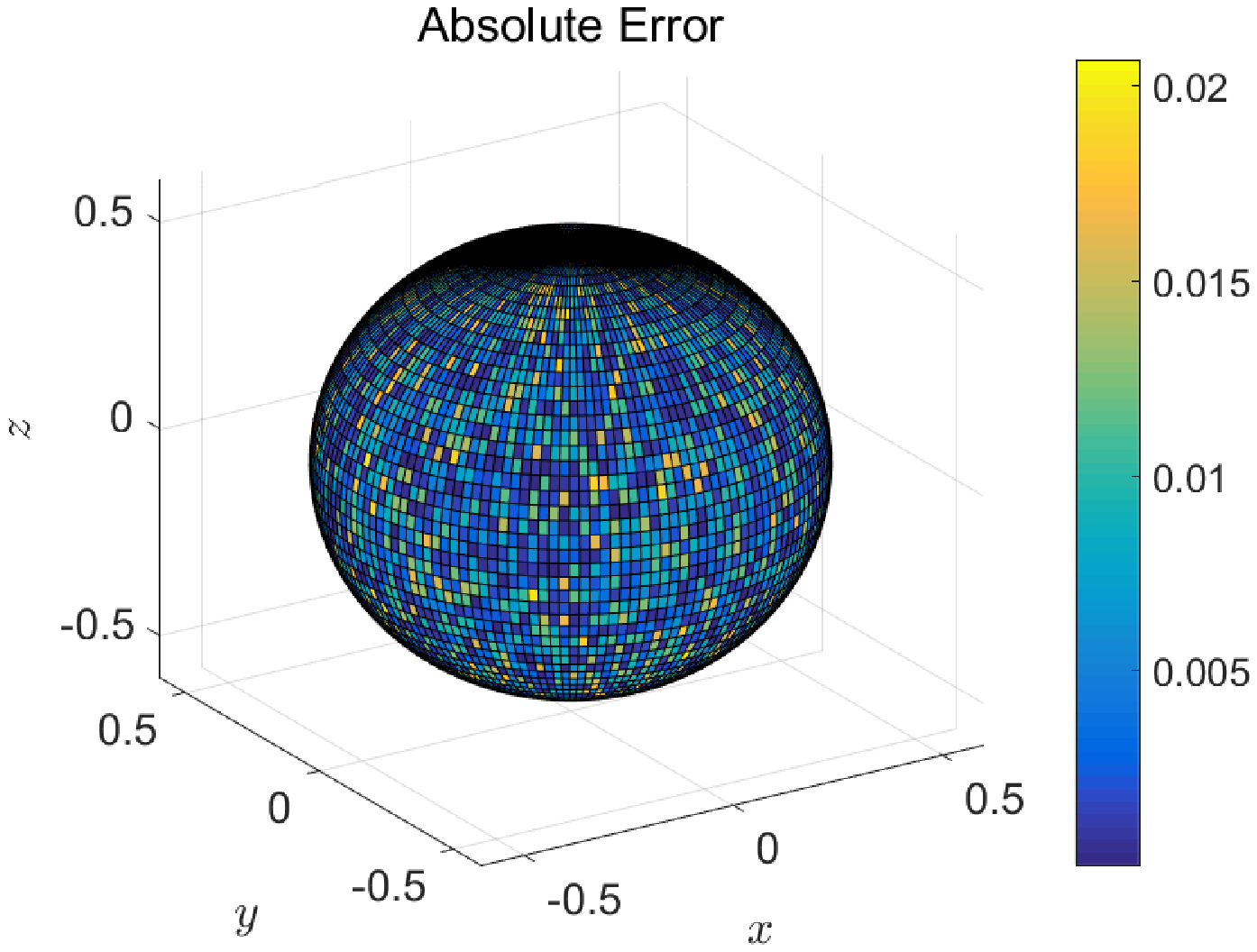}
    }
    \subfigure[MSE of PINN, SFPINN and SFHCPINN]{
        \label{Test_MSE:Multiscale_Spere3D}
        \includegraphics[scale=0.35]{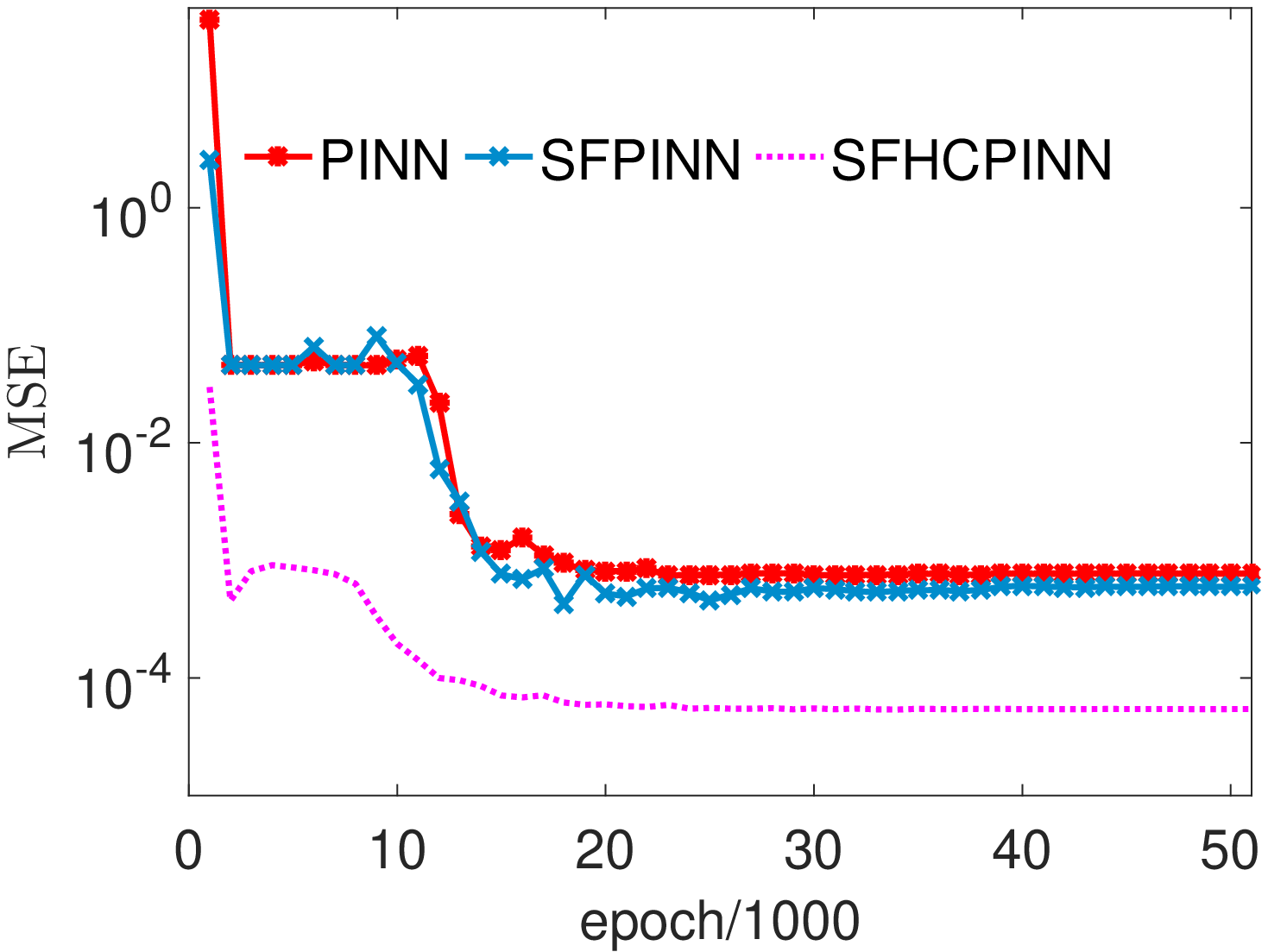}
	}
    \subfigure[REL of PINN, SFPINN and SFHCPINN]{
        \label{Test_REL:Multiscale_Spere3D}
        \includegraphics[scale=0.35]{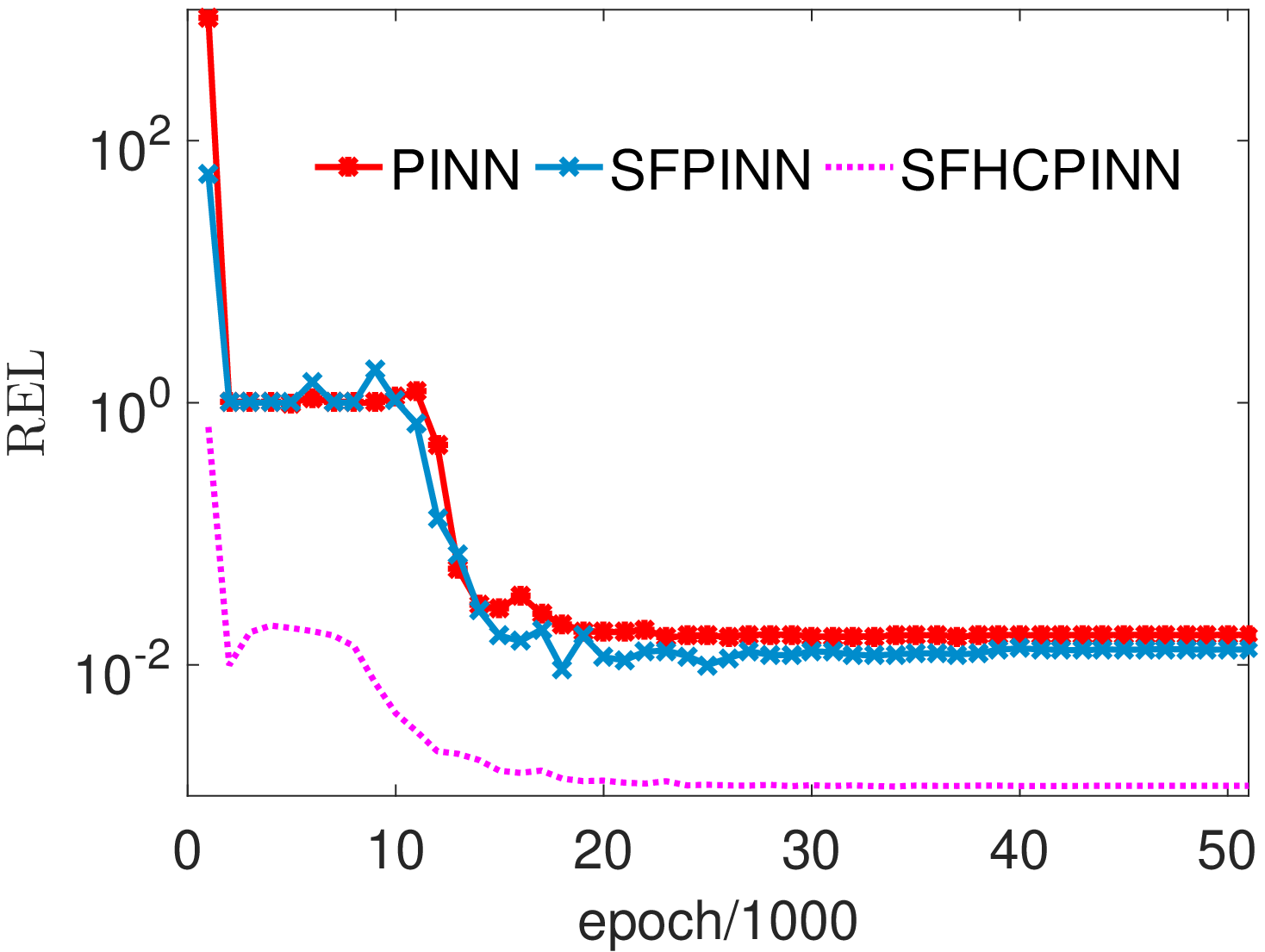}
	}
    \caption{Testing results for Example \ref{MultiscaleCase_Spere3D}.}
\end{figure}

\begin{table}[!htbp] 
\caption{MSE and REL of SFHCPINN, SFPINN, and PINN for Example \ref{MultiscaleCase_Spere3D} at $t=0.5$}
\label{Table:Multiscale_Spere3D}
\begin{center}
\setlength{\tabcolsep}{3pt}
\resizebox{!}{!}{  
    \begin{tabular}{llll} 
        \toprule 
        &constraint & MSE & REL \\
        \midrule 
         PINN      & soft &\textcolor{black}{$7.73\times 10^{-4}$}  &\textcolor{black}{$1.69\times 10^{-2}$}\\
         SFPINN    & soft &\textcolor{black}{$5.98\times 10^{-4}$}  &\textcolor{black}{$1.31\times 10^{-2}$}\\
         SFHCPINN  & hard &\textcolor{black}{$5.44\times 10^{-5}$}  &\textcolor{black}{$1.19\times 10^{-3}$}\\
        \bottomrule 
    \end{tabular}
}
\end{center}
\end{table}

\textcolor{black}{Based on the data in Table~\ref{Table:Multiscale_Spere3D} and the heatmap of three models, the performance of the normal PINN competes against that of SFPINN model, and they are inferior to the hard constraints sub-Fourier architecture when dealing with 3D issues.
In addition, we can deduce from Figs.~\ref{Test_MSE:Multiscale_Spere3D} and \ref{Test_REL:Multiscale_Spere3D} that SFHCPINN has lower initial errors, faster convergence rate, and higher precision. Additionally, Figs. \ref{PINNPWE:Multiscale_Spere3D} -- \ref{SFHCPINNPWE:Multiscale_Spere3D} and Table \ref{Table:Multiscale_Spere3D} illustrate that the point-wise errors of SFHCPINN are superior to that of PINN and SFPINN by more than two times.
In conclusion, SFHCPINN retains its high precision, convergence rate, and robust stability in solving three-dimensional ADE problems.}
\end{example}

\section{Conclusion}\label{sec:conclusion}

This study introduces SFHCPINN, a novel neural network approach that combines hard-constraint PINN with sub-networks featuring Fourier feature embedding. The purpose is to solve a specific class of advection-diffusion equations with Dirichlet and/or Neumann boundary conditions. The methodology transforms the original problem into an unconstrained optimization problem by utilizing a well-trained PINN, a distance function denoted as $D(\bm{x},t)$, and a smooth function denoted as $G(\bm{x},t)$.

To handle high-frequency modes, a Fourier activation function is employed for inputs with different frequencies, and a sub-network is designed to match the target function. The computational results demonstrate that this novel method is highly effective and efficient for solving advection-diffusion equations with Dirichlet or Neumann boundaries in one-dimensional (1D), two-dimensional (2D), and three-dimensional (3D) domains.

Importantly, SFHCPINN maintains high precision even as the dimension and/or frequency of the problem increases, unlike the soft-constraint PINN approach, which becomes degenerate in such scenarios. However, the selection of the distance function $D(\bm{x},t)$ and extension function $G(\bm{x},t)$ significantly impacts SFHCPINN's performance and may not be accurately determined in real-world engineering problems. Consequently, appropriate modifications such as employing a robust deep neural network (DNN) to fit the boundary conditions might be necessary, presenting an opportunity for future research. In addition, it shall be noted that ``the selection of the proper activation function is crucial for better accuracy and
faster convergence of the network''~\citep{jagtap2023important}, and we acknowledge the potential advantages of adaptive activation functions that outperform their classical counterparts for smooth and rough functions and problems~\citep{jagtap2020adaptive, jagtap2020locally}, such as Deep Kronecker neural networks~\citep{jagtap2022deep2}. In the future, we will attempt to explore the effectiveness of these adaptive activation functions for solving ADE problems.

\section*{Declaration of interests}
The authors declare that they have no known competing financial interests or personal relationships that could have appeared to influence the work reported in this paper.

\section*{Credit authorship contribution Statement}
Xi'an Li: Conceptualization, Methodology, Investigation, Formal analysis, Validation, Writing - Review \& Editing.
Jiaxin Deng: Investigation, Formal analysis, Validation, Writing - Original Draft. 
Jinran Wu: Writing - Original Draft,  Writing - Review \& Editing.
Shaotong Zhang: Writing - Review \& Editing.
Weide Li: Writing - Review \& Editing, Project administration.
You-Gan Wang: Writing - Review \& Editing, Project administration.

\section*{Acknowledgements}
This study was supported by the National Natural Sciences Foundation of China (No. 42130113).

\bibliographystyle{unsrtnat}
\bibliography{References}

\end{document}